\documentclass[aop,MSNbibl,seceqn,citesort,dvips]{arximspdf}
\usepackage{graphicx}

\doi{10.1214/11-AOP645}
\volume{40}
\issue{3}
\pubyear{2012}
\firstpage{1212}
\lastpage{1284}

\makeatletter

\newtheorem{theorem}{Theorem}[section]
\newtheorem{corollary}[theorem]{Corollary}
\newtheorem{lemma}[theorem]{Lemma}
\newtheorem{proposition}[theorem]{Proposition}

\newproclaim{definition}[theorem]{Definition}
\newproclaim{example}[theorem]{Example}
\newproclaim{remark}[theorem]{Remark}

\newproclaim{notation}{Notation}

\newcommand{\limfunc}[1]{\mathop{\operatorname{#1}}}
\newcommand{\func}[1]{\mathop{\operatorname{#1}}}
\newcommand{\tbigcup}{\bigcup}

\makeatother

\begin{document}
\begin{frontmatter}

\title{Two-sided estimates of heat kernels on metric measure spaces}
\runtitle{Estimates of heat kernels}

\begin{aug}
\author[A]{\fnms{Alexander} \snm{Grigor'yan}\corref{}\thanksref{t1}\ead[label=e1]{grigor@math.uni-bielefeld.de}} and
\author[B]{\fnms{Andras} \snm{Telcs}\thanksref{t2}\ead[label=e2]{telcs@szit.bme.hu}}
\runauthor{A. Grigor'yan and A. Telcs}
\affiliation{University of Bielefeld and Budapest University of
Technology and Economics}
\address[A]{Department of Mathematics\\
University of Bielefeld\\
33501 Bielefeld\\
Germany\\
\printead{e1}} %adresu isvedimo komanda gale!
\address[B]{Department of Computer Sciences\\
\quad and Information Theory\\
Budapest University of Technology\\
\quad and Economics\\
Magyar tud\'{o}sok k\"{o}r\'{u}tja 2\\
H-1117, Budapest\\
Hungary\\
\printead{e2}}
\end{aug}

\thankstext{t1}{Supported in part by SFB 701 of the German Research Council.}

\thankstext{t2}{Supported in part by a visiting grant of SFB 701 of the German Research Council.}

% HISTORY:
\received{\smonth{5} \syear{2010}}
\revised{\smonth{1} \syear{2011}}

% ABSTRACT
%
\begin{abstract}
We prove equivalent conditions for two-sided sub-Gaussian estimates of
heat kernels on
metric measure spaces.
\end{abstract}

% KEYWORDS
%
\begin{keyword}[class=AMS]
\kwd[Primary ]{47D07}
\kwd[; secondary ]{35K08}
\kwd{28A80}.
\end{keyword}
\begin{keyword}
\kwd{Heat kernel}
\kwd{heat semigroup}
\kwd{sub-Gaussian estimates}
\kwd{fractals}.
\end{keyword}

\end{frontmatter}

\tableofcontents[level=2]

%s1 ###
\section{Introduction}
\label{remmorepictures}

%s1.1 ###
\subsection{Historical background}

The notion of heat kernel has a long history. The oldest and the best-known
heat kernel is the Gauss--Weierstrass function
\[
p_{t}( x,y) =\frac{1}{( 4\pi t) ^{n/2}}\exp
\biggl( -%
\frac{\vert x-y\vert^{2}}{4t}\biggr) ,
\]
where $t>0$ and $x,y\in\mathbb{R}^{n}$, which is a fundamental
solution of
the heat equation
%
%e1.1 ###
%
\begin{equation} \label{heat}
\frac{\partial u}{\partial t}=\Delta u,
\end{equation}
where $\Delta$ is the Laplace operator in $\mathbb{R}^{n}$. A more general
parabolic equation $\frac{\partial u}{\partial t}=Lu$, where
\[
L=\sum_{i,j=1}^{n}\frac{\partial}{\partial x_{i}}\biggl( a_{ij}(
x)\, \frac{\partial}{\partial x_{j}}\biggr)
\]
is a uniformly elliptic operator with measurable coefficients
$a_{ij}=a_{ji}$%
, has also a positive fundamental solution $p_{t}( x,y) $, and
the latter admits the Gaussian bounds%
%
%e1.2 ###
%
\begin{equation} \label{Gaussian}
p_{t}( x,y) \asymp\frac{C}{t^{n/2}}\exp\biggl( -\frac
{\vert
x-y\vert^{2}}{ct}\biggr) ,
\end{equation}
where the sign $\asymp$ means that both $\leq$ and $\geq$ are true, but
the positive constants $c$ and $C$ may be different for upper and lower
bounds. The estimate~(\ref{Gaussian}) was proved by Aronson \cite{Aron}
using the parabolic Harnack inequality of Moser \cite{Moser}.

The next chapter in the history of heat kernels was opened in differential
geometry. Consider the heat equation (\ref{heat}) on a Riemannian
manifold~$M$, where $\Delta$ is now the Laplace--Beltrami operator on $M$. The heat
kernel $p_{t}( x,y) $ is defined as the minimal positive
fundamental solution of~(\ref{heat}), which always exists and is a smooth
nonnegative function of $t,x,y$; cf.
\cite{Chavbook,Grigbook,SchoenYauBook}. The question of estimating the
heat kernel on Riemannian
manifolds was addressed by many authors (see, e.g.,
\cite{Davbook,Grigbook,MolchSurvey,Varbook}). Apart from obvious
analytic and geometric motivation, a strong interest to heat kernel
estimates persists in stochastic analysis because the heat kernel coincides
with the transition density of Brownian motion on $M$ generated by the
Laplace--Beltrami operator.

One of the most powerful estimates of heat kernels was proved by Li
and Yau~\cite{LiYau}: if $M$ is a complete Riemannian manifold of
nonnegative Ricci curvature, then%
%
%e1.3 ###
%
\begin{equation}\label{LiYau}
p_{t}( x,y) \asymp\frac{C}{V( x,\sqrt{t})
}\exp
\biggl( -\frac{d^{2}( x,y) }{ct}\biggr) ,
\end{equation}
where $d( x,y) $ is the geodesic distance on $M$, and
$V(
x,r) $ is the Riemannian volume of the geodesic ball $B(
x,r) =\{ y\in M\dvtx d( x,y) <r\} $. Similar
estimates were obtained by Gushchin with coauthors \cite{Gush,GushM}
for certain unbounded domains in $\mathbb{R}^{n}$ with the Neumann boundary
condition.

An interesting question is what minimal geometric assumptions imply~(\ref{LiYau}).
The upper bound in (\ref{LiYau}) is know to be
equivalent to a
certain \textit{Faber--Krahn}-type inequality (see Section \ref
{SecEF}). The
geometric background of the lower bound in (\ref{LiYau}) is more complicated
and is closely related to the Harnack inequalities. In fact, the full
estimate (\ref{LiYau}) is equivalent, on one hand, to the parabolic
Harnack inequality of Moser (see \cite{FS}), and, on the other hand,
to the
conjunction of the volume doubling property and the Poincar\'{e} inequality
(see \cite{GrigHar,SalHar}). For a more detailed account of heat
kernel bounds on manifolds we refer the reader to the books and
surveys~\cite{Chavbook,Davbook,GrigNotes,GrigW,Grigbook,Jorg,SchoenYauBook,Varbook}.

New dimensions in the history of heat kernels were literally discovered in
analysis on fractals. Fractals are typically subsets of $\mathbb{R}^{n}$
with certain self-similarity properties, like the Sierpinski gasket
($\mathrm{SG}
$) or the Sierpinski carpet ($\mathrm{SC}$). One makes a fractal into
a metric
measure space by choosing appropriately a metric $d$ (e.g., the
extrinsic metric from the ambient~$\mathbb{R}^{n}$) and a~measure~$\mu$
(usually the Hausdorff measure). The next crucial step is introduction
of a
strongly\vspace*{1pt} local regular \textit{Dirichlet form} on a
fractal, that is, an
analogy of the Dirichlet integral $\int\vert\nabla f
\vert
^{2} $ on manifolds, which is equivalent to construction of Brownian motion
on the fractal in question; cf. \cite{FOT}. This step is highly nontrivial
and its implementation depends on a particular class of fractals. On $%
\mathrm{SG}$ Brownian motion was constructed by Goldstein~\cite{Goldstein}
and Kusuoka \cite{Kus}, on $\mathrm{SC}$, by Barlow and Bass
\cite{BarBasCon}. Kigami \cite{Kigami,Kigamibook} introduced
a~class
of \textit{post-critically finite} (p.c.f.) fractals, containing
$\mathrm{SG}$,
and constructed the Dirichlet form on such a fractal as a scaled limit of
the discrete Dirichlet forms on the graph approximations.

A strongly local regular Dirichlet form canonically leads to the notion of
the heat semigroup and the heat kernel, where the latter can be defined
either as the integral kernel of the heat semigroup or as the transition
density of Brownian motion. Surprisingly enough, the Dirichlet forms on many
families of fractals admit continuous heat kernels that satisfy the
\textit{%
sub-Gaussian} estimates%
%
%e1.4 ###
%
\begin{equation}\label{SubGauss}
p_{t}( x,y) \asymp\frac{C}{t^{\alpha/\beta}}\exp\biggl(
-c\biggl( \frac{d^{\beta}(x,y)}{t}\biggr) ^{{1}/({\beta
-1})}\biggr) ,
\end{equation}
where $\alpha>0$ and $\beta>1$ are two parameters that come from the
geometric properties of the underlying fractal. Estimate (\ref{SubGauss}
) was proved by Barlow and Perkins \cite{BarPerGas} on $\mathrm{SG}$, by
Kumagai \cite{Kumagai} on nested fractals, by Fitzsimmons, Hambly and
Kumagai \cite{FHK} on affine nested fractals and by Barlow and Bass on
$\mathrm{SC}$ \cite{BarBasTran} and on generalized Sierpinski carpets
\cite{BarBas} (see also \cite{Barlow,BBKStab,Kigamibook,KigamiVD,KusZhou}).
In fact, $\alpha$ is the Hausdorff dimension
of the space, while $\beta$ is a new quantity that is called the \textit{walk
dimension} and that can be characterized either in terms of the exit
time of
Brownian motion from balls or as the critical exponent of a family of Besov
function spaces on the fractal
(cf.~\cite{Barlow,GrigHuLau,GrigIHP,GrigHGA}).\looseness=-1

%s1.2 ###
\subsection{Description of the results}

The purpose of this paper is to find convenient equivalent conditions for
sub-Gaussian estimates of the heat kernels on abstract metric measure
spaces. Let $( M,d) $ be a locally compact separable metric
space, $\mu$~be a Radon measure on $M$ with full support and $(
\mathcal{E},\mathcal{F}) $ be a strongly local regular
Dirichlet form
on $M$ (see Section \ref{SecBasic} for the details). We are interested in
the conditions that ensure the existence of the heat kernel $p_{t}(
x,y) $ as a measurable or continuous function of $x,y$, and the
estimates of the following type%:
%
%e1.5 ###
%
\begin{equation} \label{RFi}
p_{t}( x,y) \asymp\frac{C}{V( x,\mathcal{R}
( t)
) }\exp\biggl( -ct\Phi\biggl( c\frac{d( x,y)
}{t}\biggr)
\biggr) ,
\end{equation}
where $V( x,r) =\mu( B( x,r) )
$ and $%
\mathcal{R}( t) $, $\Phi( s) $ are some nonnegative
increasing functions on $[0,\infty)$. For example, (\ref{LiYau}) has the
form (\ref{RFi}) with $\mathcal{R}( t) =\sqrt{t}$ and
$\Phi
( s) =s^{2}$, while (\ref{SubGauss}) has the form (\ref{RFi})
with $\mathcal{R}( t) =t^{1/\beta}$ and $\Phi(
s)
=s^{{\beta}/({\beta-1})}$ [assuming that\setcounter
{footnote}{2}\footnote{%
The sign $\simeq$ means that the ratio of both sides is bounded between
two positive constants.} $V( x,r) \simeq r^{\alpha}$,
which, in
fact, follows from (\ref{SubGauss})].

To describe the results of the paper, let us introduce some hypotheses.
First, we assume that the metric space $( M,d) $ is unbounded
and that all metric balls are precompact\footnote{%
The precompactness of balls implies that $( M,d) $ is a complete
metric spaces. The following partial converse is also true: if $(
M,d) $ is complete and the volume doubling property (\ref{VD})
holds, then all balls are precompact. However, since we do not always
assume (\ref{VD}), we make an independent assumption of
precompactness of
the balls.} (although these assumptions are needed only for a part of the
results). Next, define the following conditions:

$\bullet$ the volume doubling property (\ref{VD}): there
is a
constant $C$ such that%
{\renewcommand{\theequation}{\textit{VD}}
\begin{equation}\label{VD}
V( x,2r) \leq CV( x,r)
\end{equation}}
for all $x\in M$ and $r>0$;\vadjust{\goodbreak}

$\bullet$ the elliptic Harnack inequality (\ref{H}): there
is a
constant $C$ such that, for any nonnegative harmonic function $u$ in any
ball $B( x,r) \subset M$,%
{\renewcommand{\theequation}{$H$}
\begin{equation}\label{H}
\limfunc{esup}_{B( x,r/2) }u\leq C\limfunc{einf}_{B(
x,r/2) }u,
\end{equation}}

\vspace*{-8pt}

\noindent
where $\limfunc{esup}$ and $\limfunc{einf}$ are the essential
supremum and
infimum, respectively, (see Section \ref{SecHarnack} for more details);

$\bullet$ the estimate of the mean exit time (\ref{EF}),
{\renewcommand{\theequation}{${E}_{F}$}
\begin{equation}\label{EF}
\mathbb{E}_{x}\tau_{B( x,r) }\simeq F( r) ,
\end{equation}}

\vspace*{-\baselineskip}

\noindent
where $\tau_{B( x,r) }$ is the first exist time from ball
$%
B( x,r) $ of the associated diffusion process, started at the
center $x$, and $F( r) $ is a given function with a certain
regularity (see Section \ref{SecEF} for more details). A typical
example is $%
F( r) =r^{\beta}$ for some constant $\beta>1$.

The conditions $\mbox{(\ref{H})} + \mbox{(\ref{VD})} +
\mbox{(\ref{EF})} $ are
known to be true on p.c.f. fractals (see
\cite{Kigamibook,HamblyKum}) as well as on generalized Sierpinski carpets
(see \cite{BarBas,BBKT}) so that our results apply to such fractals.
Another situation
where $%
\mbox{(\ref{H})} +\mbox{(\ref{VD})} +\mbox{(\ref{EF})} $ are
satisfied is
the setting of \textit{resistance forms} introduced by Kigami~\cite{KigamiRes}.
A resistance form is a specific Dirichlet form that
corresponds to a~strongly recurrent Brownian motion. Kigami showed that, in the setting of
resistance forms on self-similar sets, (\ref{VD}) alone
implies (\ref{H}) and (\ref{EF}) with $F( r)
=r^{\beta}$, for a suitable choice of a distance function. Examples with
more general functions $F( r) $ appear in \cite{BarHam}
and \cite{Telcsbook}.

Let us emphasize in this connection that our results do not depend on the
recurrence or transience hypotheses and apply to both cases, which partly
explains the complexity of the proofs. A transient case occurs, for example,
for some generalized Sierpinski carpets. Another point worth mentioning is
that we do not assume specific properties of the metric $d$ such as being
geodesic; the latter is quite a common assumption in the fractal literature.
This level of generality enables applications to resistance forms where the
distance function is usually the resistance metric that is not geodesic.

Our first main result, which is stated in Theorem \ref{Tmain} and
which, in
fact, is a combination of Theorems \ref{TG=>FK}, \ref{TDUE}, \ref
{THolder}, %
\ref{TNLE}, says the following: if the hypotheses $\mbox{(\ref{VD})}
+\mbox{(\ref{H})} +\mbox{(\ref{EF})}$ are satisfied, then the heat
kernel $
p_{t}( x,y) $ exists, is H\"{o}lder continuous in $x,y$ and
satisfies the following upper estimate:%
{\renewcommand{\theequation}{\textit{UE}}
\begin{equation}\label{UE}
p_{t}( x,y) \leq\frac{C}{V( x,\mathcal{R}(
t)
) }\exp\biggl( -\frac{1}{2}t\Phi\biggl( c\frac{d(
x,y) }{t}%
\biggr) \biggr),
\end{equation}}

\vspace*{-8pt}

\noindent
where $\mathcal{R} = F^{-1}$ and
\[
\Phi( s) :=\sup_{r>0}\biggl\{ \frac{s}{r}-\frac
{1}{F(
r) }\biggr\} ,
\]
and the \textit{near-diagonal} lower estimate%
{\renewcommand{\theequation}{\textit{NLE}}
\begin{equation}\label{NLE}
p_{t}( x,y) \geq\frac{c}{V( x,\mathcal{R}(
t)
) } \qquad\mbox{provided }d( x,y) \leq\eta\mathcal
{R}(t),
\end{equation}}

\vspace*{-13pt}
\noindent where $\eta>0$ is a small enough constant. Furthermore, assuming that
(\ref{VD}) holds a priori, we have the equivalence\footnote{%
For comparison, let us observe that, under the same standing
assumptions, it
was proved in \cite{BarGrigKumHar} that
\[
\mbox{(\ref{UE})} + \mbox{(\ref{NLE})} \Leftrightarrow(
\mathit{PHI}_{F}),
\]
where $(\mathit{PHI}_{F}) $ stands for the \textit{parabolic} Harnack
inequality for caloric functions. Hence, we see that the
``difference'' between $( \mathit{PHI}_{F}) $ and
(\ref{H}) is the condition (\ref{EF}), that in particular
provides a necessary space/time scaling for $( \mathit{PHI}_{F}) $.}
%
%e1.6 ###
\setcounter{equation}{5}
\begin{equation} \label{eq}
\mbox{(\ref{UE})} + \mbox{(\ref{NLE})} \Leftrightarrow\mbox{(\ref{H})}
+\mbox{(\ref{EF})}
\end{equation}
(Theorem \ref{Tconv}).

For example, if $F( r) =r^{\beta}$ for some $\beta>1$,
then $%
\mathcal{R}( t) =t^{1/\beta}$ and $\Phi( s
) =\func{%
const}s^{{\beta}/({\beta-1})}$. Hence, (\ref{UE}) and
(\ref{NLE}) become as follows:%
%
%e1.7 ###
%
\begin{equation}\label{subu}
p_{t}( x,y) \leq\frac{C}{V( x,t^{1/\beta})
}\exp\biggl( -c\biggl( \frac{d^{\beta}(x,y)}{t}\biggr) ^{{1}/({\beta-1})}\biggr)
\end{equation}
and%
\[
p_{t}( x,y) \geq\frac{c}{V( x,t^{1/\beta})}\qquad
\mbox{provided }d( x,y) \leq\eta t^{1/\beta}.
\]
It is desirable to have a lower bound of $p_{t}( x,y) $
for all $%
x,y$ that would match the upper bound (\ref{subu}). However, such a lower
bound fails in general. The reason for that is the lack of \textit{chaining
properties} of the distance function, where by chaining properties we
loosely mean a possibility to connect any two points $x,y\in M$ by a chain
of balls of controllable radii so that the number of balls in this
chain is
also under control. More precisely, this property can be stated in
terms of
the modified distance $d_{\varepsilon}( x,y) $ where $%
\varepsilon>0$ is a~parameter. The exact definition of $d_{\varepsilon}$
is given in Section \ref{Secde}, where it is also shown that
\[
d_{\varepsilon}( x,y) \simeq\varepsilon N_{\varepsilon
}(x,y) ,
\]
where $N_{\varepsilon}( x,y) $ is the smallest number of balls
in a chain of balls of radii $\varepsilon$ connecting $x$ and $y$. As $
\varepsilon$ goes to $0$, $d_{\varepsilon}( x,y) $ increases
and can go to $\infty$ or even become equal to $\infty$. If the distance
function $d$ is geodesic then $d_{\varepsilon}\equiv d$, which corresponds
to the best possible chaining property. In general, the rate of growth
of $%
d_{\varepsilon}( x,y) $ as $\varepsilon\rightarrow0$
can be
regarded as a quantitative description of the chaining properties of~$d$.
For this part of our work, we assume that%
%
%e1.8 ###
%
\begin{equation}\label{ad}
\frac{F( \varepsilon) }{\varepsilon}d_{\varepsilon}
(x,y) \rightarrow0 \qquad\mbox{as }\varepsilon\rightarrow0,
\end{equation}
which allows to define a function $\varepsilon( t,x,y) $ from
the identity%
%
%e1.9 ###
%
\begin{equation} \label{ede}
\frac{F( \varepsilon) }{\varepsilon}d_{\varepsilon}
(x,y) =t.
\end{equation}
Our second main result states the following: if (\ref{ad}) and $\mbox
{(\ref{VD})}
+\mbox{(\ref{H})} + \mbox{(\ref{EF})}$ are satisfied,
then%
%
%e1.11 ###
%e1.10 ###
%
\begin{eqnarray}
\label{twoe}
p_{t}( x,y) &\asymp&\frac{C}{V(x,\mathcal{R}(
t) )}%
\exp\biggl( -ct\Phi\biggl( c\frac{d_{\varepsilon}( x,y
) }{t}%
\biggr) \biggr) \\
\label{twohk}
&\asymp&\frac{C}{V(x,\mathcal{R}( t) )}\exp(
-cN_{\varepsilon}) ,
\end{eqnarray}
where $\varepsilon=\varepsilon( ct,x,y) $ (Theorem \ref{Ttwo}).
For example, the above hypotheses and, hence, the estimates (\ref
{twoe}) and (%
\ref{twohk}) hold on connected p.c.f. fractals endowed with resistance
distance, where one has $V( x,r) \simeq r^{\alpha}$ and
$F(
r) =r^{\alpha+1}$ for some constant $\alpha$. The estimate
(\ref{twohk}) on p.c.f. fractals was first proved by Hambly and Kumagai
\cite{HamblyKum}. In fact, we use the argument from \cite{HamblyKum} to verify
our hypotheses (see Remark \ref{ExHK}).

Note that the dependence on $t,x,y$ in the estimates (\ref{twoe}) and (\ref
{twohk}) in very implicit and is hidden in $\varepsilon(
ct,x,y)
$. One can loosely interpret the use of this function in (\ref
{twoe}) and (\ref{twohk}) as follows. In order to find a most probable path
for Brownian
motion to go from $x$ to $y$ in time $t$, one determines the optimal
size $%
\varepsilon=\varepsilon( ct,x,y) $ of balls and then the
optimal chain of balls of radii~$\varepsilon$ connecting $x$ and $y$, and
this chain provides an optimal route between~$x$ and~$y$. This phenomenon
was discovered by Hambly and Kumagai in the setting of p.c.f. fractals,
where they used instead of balls the construction cells of the fractal. As
it follows from our results, this phenomenon is generic and independent of
self-similar structures.

If the distance function satisfies the \textit{chain condition} $%
d_{\varepsilon}\leq Cd$, which is stronger than (\ref{ad}), then one can
replace in (\ref{twoe}) $d_{\varepsilon}$ by $d$ and obtain (\ref{RFi})
(Corollary \ref{Cortwo}). In fact, in this case we have the equivalence
%
%e1.12 ###
%
\begin{equation} \label{3+2}
\mbox{(\ref{VD})} +\mbox{(\ref{H})} +\mbox{(\ref{EF})}
\Leftrightarrow\mbox{(\ref{RFi})}
\end{equation}
(Corollary \ref{Corunbound}).

In the setting of random walks on infinite graphs, the equivalence
(\ref{3+2}%
) was proved by the authors in \cite{GrigTelTran,GrigTelLoc}. Of
course, in this case all the conditions have to be adjusted to the discrete
setting.

For the sake of applications (cf., e.g, \cite{BBKT}), it is desirable
to replace the probabilistic condition (\ref{EF}) in all the
above results by an analytic condition, namely, by a certain estimate
of the
capacity between two concentric balls. This type of result requires
different techniques and will be treated elsewhere.

%s1.3 ###
\subsection{Structure of the paper and interconnection of the results}

In Section~\ref{SecHH} we revise the basic properties of the heat semigroups
and heat kernels and prove the criterion for the existence of the heat
kernel in terms of local ultracontractivity of the heat semigroup
(Theorem %
\ref{TptOm}).

In Section \ref{Secaux} we prove two preparatory results:

\begin{longlist}[(2)]
\item[(1)]
$\mbox{(\ref{VD})} +\mbox{(\ref{H})} +\mbox{(\ref{EF})}
\Rightarrow
\mbox{(\ref{FK})}$ where (\ref{FK}) stands for a certain
\textit{%
Faber--Krahn inequality}, which provides a lower bound for the bottom
eigenvalue in any bounded open set $\Omega\subset M$ via its measure
(Theorem \ref{TG=>FK}). In turn,~(\ref{FK}) implies the local
ultracontractivity of the heat semigroup, which by Theorem~\ref{TptOm}
ensures the existence of the heat kernel.

\item[(2)] $( E_{F}) $ implies the following estimate of the tail
of the
exit time from balls:%
%
%e1.13 ###
%
\begin{equation} \label{Psi2}
\mathbb{P}_{x}\bigl( \tau_{B( x,R) }\leq t\bigr) \leq
C\exp
\biggl( -t\Phi\biggl( c\frac{R}{t}\biggr) \biggr)
\end{equation}
(Theorem \ref{TEF}).
\end{longlist}

In Section \ref{SecDUE} we prove the upper estimate of the heat
kernel, more
precisely, the implication%
\[
\mbox{(\ref{VD})} + \mbox{(\ref{FK})} +\mbox{(\ref{EF})}
\Rightarrow\mbox{(\ref{UE})}
\]
(Theorem \ref{TDUE}). The main difficulty lies already in the proof of the
diagonal upper bound%
{\renewcommand{\theequation}{\textit{DUE}}
\begin{equation}\label{DUE}
p_{t}( x,x) \leq\frac{C}{V( x,\mathcal{R}(
t)
) }.
\end{equation}}

\vspace*{-8pt}

\noindent
Using (\ref{FK}), we obtain first some diagonal upper bound
for the
Dirichlet heat kernels in balls, and then use Kigami's iteration argument
and (\ref{Psi2}) to pass to (\ref{DUE}). The latter
argument is
borrowed from \cite{GrigHuUpper}. The full upper estimate (\ref{UE})
follows from (\ref{DUE}) and (\ref{Psi2}).

In Section \ref{SecLow} we prove the lower bounds of the heat kernel. The
diagonal lower bound%
{\renewcommand{\theequation}{\textit{DLE}}
\begin{equation}\label{DLE}
p_{t}( x,x) \geq\frac{C}{V( x,\mathcal{R}(
t)
) }
\end{equation}}

\vspace*{-8pt}

\noindent
follows directly from (\ref{Psi2}) (Lemma \ref{LemDLE}). To obtain
the near
diagonal lower estimate (\ref{NLE}), one estimates from
above the
difference%
%
%e1.14 ###
%
\setcounter{equation}{13}
\begin{equation}\label{di}
\vert p_{t}( x,x) -p_{t}( x,y)
\vert,
\end{equation}
where $y$ is close to $x$, which requires the following two ingredients:

\begin{longlist}[(2)]
\item[(1)]
the oscillation inequalities that are consequences of the elliptic
Harnack inequality (\ref{H}) (Lemma \ref{Lemosc} and
Proposition \ref{Posc});

\item[(2)] the upper estimate of the time derivative $\partial_{t}p_{t}(
x,y) $ (Corollary \ref{Cdtpt}).
\end{longlist}

Combining them with (\ref{UE}), one obtains an upper bound
for (\ref{di}), which together with (\ref{DLE}) yields (\ref{NLE})
(Theorem \ref{TNLE}).

The same method gives also the H\"{o}lder continuity of the heat kernel
(Theorem~\ref{THolder}).

In Section \ref{SecTwo} we prove two-sided estimates (\ref{twoe}) and (\ref
{twohk}) (Theorem~\ref{Ttwo}). For the upper bound, we basically
repeat the
proof of (\ref{UE}) by tracing the use of the distance
function $d$
and replacing it by $d_{\varepsilon}$. The lower bound for large
$d(
x,y) $ is obtained from (\ref{NLE}) by a standard chaining
argument using the semigroup property of the heat kernel and the chaining
property of the distance function.

In Section \ref{Secconv} we prove the converse Theorem \ref{Tconv}, which
essentially consists of the equivalence (\ref{eq}).\vspace*{-3pt}
\begin{notation}
We use the letters $C,c,C^{\prime},c^{\prime}$ etc. to denote positive
constant whose value is unimportant and can change at each occurrence. Note
that the value of such constants in the conclusions depend on the
values of
the constants in the hypotheses (and, perhaps, on some other explicit
parameters). In this sense, all our results are quantitative.

The relation $f\simeq g$ means that $C^{-1}g\leq f\leq Cg$ for some positive
constant~$C$ and for a specified range of the arguments of functions
$f$ and
$g$. The relation $f\asymp g$ means that both inequalities $f\leq g$
and $%
f\geq g$ hold but possibly with different values of constants $c,C$
that are
involved in the expressions~$f$ and/or $g$.\vspace*{-3pt}
\end{notation}

%s2 ###
\section{Heat semigroups and heat kernels}\vspace*{-3pt}
\label{SecHH}

%s2.1 ###
\subsection{Basic setup}

\label{SecBasic}Throughout the paper, we assume that $(
M,d) $ is a~locally compact separable metric space, and $\mu$ is a Radon measure
on $M$
with full support. As usual, denote by $L^{q}( M) $ where
$q\in%
[ 1,+\infty] $ the Lebesgue function space with respect
measure $%
\mu$, and by $\Vert\cdot\Vert_{q}$ the norm in
$L^{q}(
M) $. The inner product in $L^{2}( M) $ is denoted
by $%
( \cdot,\cdot) $. All functions on $M$ are supposed to
be real
valued. Denote by $C_{0}( M) $ the space of all continuous
functions on $M$ with compact supports, equipped with the $\sup$-norm.

Let $( \mathcal{E},\mathcal{F}) $ be Dirichlet form in $%
L^{2}( M) $. This means that $\mathcal{F}$ is a dense
subspace of
$L^{2}( M) $, and $\mathcal{E}( f,g) $ is a bilinear,
nonnegative definite, closed\footnote{The form $(\mathcal{E},\mathcal{F})$ is called
closed if $\mathcal{F}$ is a Hilbert space with respect to the following inner product:
\[
\mathcal{E}_{1}(f,g)=\mathcal{E}(f,g)+(f,g).
\]}
form defined for functions $f,g\in\mathcal{F}$, which satisfies, in
addition, the Markovian property.\footnote{%
The Markovian property (which could be also called the Beurling--Deny
property) means that if $f\in\mathcal{F}$, then also the function
$\hat{f}%
=f_{+}\wedge1$ belongs to $\mathcal{F}$ and $\mathcal{E(}\hat
{f},\hat{f}%
)\leq\mathcal{E}( f,f)$.} The Dirichlet form $(
\mathcal{E%
},\mathcal{F}) $ is called \textit{regular} if $\mathcal{F}\cap
C_{0}( M) $ is dense both in $\mathcal{F}$ and in
$C_{0}(
M) $. The Dirichlet form is called \textit{strongly local} if
$\mathcal{E%
}( f,g) =0$ for all functions $f,g\in\mathcal{F}$ such
that $g$
has a compact support and $f\equiv\func{const}$ in a neighborhood of $
\limfunc{supp}g$. In this paper, we assume by default that $(
\mathcal{E%
},\mathcal{F}) $ is a regular, strongly local Dirichlet form. A
general theory of Dirichlet forms can be found in \cite{FOT}.\vadjust{\goodbreak}

Let $\mathcal{L}$ be the generator of $( \mathcal{E},\mathcal
{F})
$; that is, $\mathcal{L}$ is a self-adjoint nonnegative definite operator
in $L^{2}( M) $ with the domain $\func{dom}(
\mathcal{L}%
) $ that is a dense subset of $\mathcal{F}$ and such that, for
all $%
f\in\func{dom}( \mathcal{L}) $ and $g\in\mathcal{F}$%
\[
\mathcal{E}( f,g) =( \mathcal{L}f,g) .
\]
The associated \textit{heat semigroup}
\[
P_{t}=e^{-t\mathcal{L}}, \qquad t\geq0,
\]
is a family of bounded self-adjoint operators in $L^{2}( M
) $. The
Markovian properties allow the extension of $P_{t}$ to a bounded
operator in
$L^{q}( M) $, with the norm $\leq1$, for any $q\in
\lbrack
1,+\infty]$.

Denote by $\mathcal{B}( M) $ the class of all Borel
functions on $%
M$, by $\mathcal{B}_{b}$ the class of bounded Borel functions, by
$\mathcal{B%
}_{+}( M) $ the class of nonnegative Borel functions and
by $%
\mathcal{B}L^{q}( M) $ the class of Borel functions that belong
to $L^{q}( M) $.

By \cite{FOT}, Theorem 7.2.1, for any local Dirichlet form, there
exists a
diffusion process $\{ \{ X_{t}\} _{t\geq0},\{
\mathbb{P%
}_{x}\} _{x\in M\setminus\mathcal{N}_{0}}\} $ with the initial
point $x$ outside some properly exceptional set\footnote{%
A set $\mathcal{N}\subset M$ is called properly exceptional if it is
Borel, $%
\mu( \mathcal{N}) =0$ and%
\[
\mathbb{P}_{x}( X_{t}\in\mathcal{N}\mbox{ for some }t\geq
0) =0
\]
for all $x\in M\setminus\mathcal{N}$ (see \cite{FOT},
page 134 and Theorem 4.1.1 on page 137).} $\mathcal{N}_{0}\subset M$,
which is associated
with the heat semigroup $\{ P_{t}\} $ as follows: for any
$f\in
\mathcal{B}L^{q}( M) $, $1\leq q\leq\infty$,%
%
%e2.1 ###
%
\begin{equation}\label{Ee}
\mathbb{E}_{x}f( X_{t}) =P_{t}f( x)
\qquad\mbox{for }\mu\mbox{-a.a. } x\in M.
\end{equation}
Consider the family of operators $\{ \mathcal{P}_{t}\}
_{t\geq0}$
defined by%
%
%e2.2 ###
%
\begin{equation}\label{PtEx}
\mathcal{P}_{t}f( x) :=\mathbb{E}_{x}f
(X_{t}),\qquad x\in M\setminus\mathcal{N}_{0},
\end{equation}
for all functions $f\in\mathcal{B}_{b}( M) $ (if $X_{t}$
has a
finite lifetime, then $f$ is to be extended by $0$ at the cemetery). The
function $\mathcal{P}_{t}f( x) $ is a bounded Borel
function on $%
M\setminus\mathcal{N}_{0}$. \label{remdoweneedthisextension?}
It is convenient to extend it to all $x\in M$ by setting%
%
%e2.3 ###
%
\begin{equation}\label{PtEx0}
\mathcal{P}_{t}f( x) =0,\qquad x\in\mathcal{N}_{0},
\end{equation}
so that $\mathcal{P}_{t}$ can be considered as an operator in
$\mathcal{B}%
_{b}( M) $. Obviously, $\mathcal{P}_{t}f\geq0$ if $f\geq
0$ and $%
\mathcal{P}_{t}1\leq1$. Moreover, the family $\{ \mathcal{P}%
_{t}\} _{t\geq0}$ satisfies the semigroup identity
\[
\mathcal{P}_{t}\mathcal{P}_{s}=\mathcal{P}_{t+s}.
\]
Indeed, if $x\in M\setminus\mathcal{N}_{0}$, then we have by the Markov
property, for any $f\in\mathcal{B}_{b}( M) $,%
\[
\mathcal{P}_{t+s}f( x) =\mathbb{E}_{x}( f(
X_{t+s}) ) =\mathbb{E}_{x}( \mathbb
{E}_{X_{t}}(
f( X_{s}) ) ) =\mathbb{E}_{x}(
\mathcal{P}%
_{s}f( X_{t}) ) =\mathcal{P}_{t}( \mathcal
{P}%
_{s}f) ( x),
\]
where we have used that $X_{t}\in M\setminus\mathcal{N}_{0}$ with
$\mathbb{P%
}_{x}$-probability $1$. If $x\in\mathcal{N}_{0}$, then we have again%
\[
\mathcal{P}_{t+s}f( x) =\mathcal{P}_{t}( \mathcal
{P}%
_{s}f) ( x),
\]
because the both sides are $0$.

By considering increasing sequences of bounded functions, one extends the
definition of $\mathcal{P}_{t}f$ to all $f\in\mathcal{B}_{+}(
M)
$ so that the defining identities (\ref{PtEx}) and (\ref{PtEx0}) remain
valid also for $f\in\mathcal{B}_{+}( M) $ [allowing
value $%
+\infty$ for $\mathcal{P}_{t}f( x) $]. For a signed
function $%
f\in\mathcal{B}( M) $, define $\mathcal{P}_{t}f$ by%
\[
\mathcal{P}_{t}f( x) =\mathcal{P}_{t}( f_{+}
) (x) -\mathcal{P}_{t}( f_{-}) ( x) ,
\]
provided at least one of the functions $\mathcal{P}_{t}(
f_{+}) $%
, $\mathcal{P}_{t}( f_{-}) $ is finite. Obviously, identities
(\ref{PtEx}), (\ref{PtEx0}) are satisfied for such functions as well.

If follows from the comparison of (\ref{Ee}) and (\ref{PtEx}) that,
for all $%
f\in\mathcal{B}L^{q}( M) $,
\[
\mathcal{P}_{t}f( x) =P_{t}f( x) \qquad\mbox
{for }\mu
\mbox{-a.a. }x\in M.
\]
It particular, $\mathcal{P}_{t}f$ is finite almost everywhere.

The set of the above assumptions will be referred to as the \textit{basic
hypotheses}, and they are assumed by default in all parts of this paper.
Sometimes we need also the following property.
\begin{definition}
The Dirichlet form $( \mathcal{E},\mathcal{F}) $ is called
\textit{conservative} (or \textit{stochastically complete}) if $\mathcal
{P}%
_{t}1\equiv1$ for all $t>0$.
\end{definition}
\begin{example}
Let $M$ be a connected Riemannian manifold, $d$ be the geodesic distance
on $M$, $\mu$ be the Riemannian volume. Define the Sobolev space%
\[
W^{1}=\{ f\in L^{2}( M) \dvtx\nabla f\in L^{2}(
M)\},
\]
where $\nabla f$ is the Riemannian gradient of $f$ understood in the weak
sense. For all $f,g\in W^{1}$, one defines the energy form%
\[
\mathcal{E}( f,g) =\int_{M}( \nabla f,\nabla
g) \,d\mu.
\]
Let $\mathcal{F}$ be the closure of $C_{0}^{\infty}( M)
$ in $%
W^{1}$. Then $( \mathcal{E},\mathcal{F}) $ is a regular strongly
local Dirichlet form in $L^{2}( M) $. \label
{remexamplesofDirichletformonfractals}
\end{example}

%s2.2 ###
\subsection{The heat kernel and the transition semigroup}
\label{SecHeat}

\begin{definition}
\label{DefHK}
The \textit{heat kernel} (or the \textit{transition density})
of the transition semigroup $\{ \mathcal{P}_{t}\} $ is a function
$p_{t}( x,y) $ defined for all $t>0$ and $x,y\in
D:=M\setminus
\mathcal{N}$, where $\mathcal{N}$ is a properly exceptional set
containing $%
\mathcal{N}_{0}$, and such that the following properties are satisfied:

\begin{longlist}[(2)]
\item[(1)] for any $t>0$, the function $p_{t}( x,y) $ is measurable
jointly in $x,y$;\vadjust{\goodbreak}

\item[(2)] for all $f\in\mathcal{B}_{+}( M) $,
$t>0$ and $%
x\in D$,%
%
%e2.4 ###
%
\begin{equation}\label{Pt=pt}
\mathcal{P}_{t}f( x) =\int_{D}p_{t}( x,y)
f(
y) \,d\mu( y) ;
\end{equation}

\item[(3)] for all $t>0$ and $x,y\in D$,
%
%e2.5 ###
%
\begin{equation}\label{sym}
p_{t}( x,y) =p_{t}( y,x) ;
\end{equation}

\item[(4)] for all $t,s>0$ and $x,y\in D$,%
%
%e2.6 ###
%
\begin{equation} \label{semi}
p_{t+s}( x,y) =\int_{D}p_{t}( x,z)
p_{s}(
z,y) \,d\mu( z) .
\end{equation}
The set $D$ is called the domain of the heat kernel.
\end{longlist}
\end{definition}

Let us extend $p_{t}( x,y) $ to all $x,y\in M$ by setting $
p_{t}( x,y) =0$ if $x$ or $y$ is outside~$D$. Then (\ref{sym})
and (\ref{semi}) hold for all $x,y\in M$, and the domain of
integration in (%
\ref{Pt=pt}) and (\ref{semi}) can be extended to $M$. The existence
of the
heat kernel allows us to extend the definition of $\mathcal{P}_{t}f$
to all
measurable functions~$f$ by choosing a Borel measurable version of~$f$ and
noticing that the integral~(\ref{Pt=pt}) does not change if function
$f$ is
changed on a set of measure~$0$.

It follows from (\ref{Ee}) and (\ref{Pt=pt}) that, for any $f\in
L^{2}(
M) $,
%
%e2.7 ###
%
\begin{equation}\label{ept}
P_{t}f( x) =\int_{M}p_{t}( x,y) f(
y) \,d\mu
( y)
\end{equation}
for all $t>0$ and $\mu$-a.a. $x\in M$. A measurable
function $%
p_{t}( x,y) $ that satisfies~(\ref{ept}) is called the heat
kernel of the semigroup $P_{t}$. It is well known that the heat kernel
of $%
P_{t}$ satisfies (\ref{sym}) and (\ref{semi}) although for \textit
{almost} all
$x,y\in M$ (see \cite{GrigHuUpper}, Section~3.3).

Hence, the relation between the heat kernels of $\mathcal{P}_{t}$ and $P_{t}$
is as follows: the former is defined as a pointwise function of $x,y$, while
the latter is defined almost everywhere, and the former is a pointwise
realization of the latter, where the defining identities (\ref{Pt=pt}),
(\ref{sym}), (\ref{ept}) must be satisfied pointwise. In this paper the heat
kernel is understood exclusively in the sense of Definition~\ref{DefHK}.

The existence of the heat kernel is not obvious at all and will be
given a
special treatment. Those who are interested in the settings where the
pointwise existence of the heat kernel is known otherwise, can skip the rest
of this section and go to Section \ref{Secaux}.
\begin{lemma}
\label{Lemuniq}Let $p_{t}$ be the heat kernel of $\mathcal{P}_{t}$.

\begin{longlist}[(a)]
\item[(a)]
The function $p_{t}( x,\cdot) $
belongs to $%
\mathcal{B}L^{2}( M) $ for all $t>0$ and $x\in M$.

\item[(b)] For all $t>0$, $x,y\in M$, we have $p_{t}(
x,y) \geq0$ and
%
%e2.8 ###
%
\begin{equation}\label{int1}
\int_{M}p_{t}( x,z) \,d\mu( z) \leq1.
\end{equation}
Consequently, $p_{t}( x,\cdot) \in\mathcal{B}L^{1}(
M) $.\vadjust{\goodbreak}

\item[(c)] If $q_{t}$ is another heat kernel, then
$p_{t}=q_{t}$ in
the common part of their domains.
\end{longlist}
\end{lemma}
\begin{pf}
(a) Set $f=p_{t/2}( x,\cdot) $ and
observe that,
by (\ref{sym}) and (\ref{semi}),%
%
%e2.9 ###
%
\begin{equation}\label{t2}
p_{t}( x,y) =\int_{M}p_{t/2}( x,\cdot)
p_{t/2}(
y,\cdot) \,d\mu=\mathcal{P}_{t/2}f( y)
\end{equation}
for all $t>0$ and $x,y\in D$. Since $\mathcal{P}_{t/2}f$ is a Borel
function, we obtain that $p_{t}( x,\cdot) $ is Borel. The latter
is true also if $x\in\mathcal{N}$ since in this case $p_{t}(
x,\cdot
) =0$. Setting in (\ref{t2}) $x=y$, we obtain%
%
%e2.10 ###
%
\begin{equation}\label{pt2}
\int_{M}p_{t/2}( x.\cdot) ^{2}\,d\mu=p_{t}(
x,x)
<\infty,
\end{equation}
whence $p_{t/2}( x,\cdot) \in L^{2}( M) $.

(b) By (\ref{PtEx}), (\ref{PtEx0}) we have $\mathcal
{P}%
_{t}f( x) \geq0$ for all $t>0$, $x\in M$ provided $f\geq0$.
Setting $f=[ p_{t}( x,\cdot) ] _{-}$, we
obtain%
\[
0\leq\mathcal{P}_{t}f( x) =\int_{M}p_{t}( x,\cdot
) %
[ p_{t}( x,\cdot) ] _{-}\,d\mu=-\int
_{M}[
p_{t}( x,\cdot) ] _{-}^{2}\,d\mu,
\]
whence it follows that $[ p_{t}( x,\cdot) ]
_{-}=0$ a.e., that is, $p_{t}( x,\cdot) \geq0$
a.e.
on $M$. It follows from (\ref{t2}) that, for all $x,y\in M$,%
\[
p_{t}( x,y) =\int_{M}p_{t/2}( x,\cdot)
p_{t/2}(
y,\cdot) \,d\mu\geq0.
\]
Inequality (\ref{int1}) is trivial if $x\in\mathcal{N}$, and if
$x\in D$
then it follows from%
\[
\int_{M}p_{t}( x,\cdot) \,d\mu=\mathcal{P}_{t}1(
x) =%
\mathbb{E}_{x}1\leq1.
\]

(c) Let $D$ be the intersection of the domains of
$p_{t}$ and
$q_{t}$. For all $f\in\mathcal{B}_{+}( M) $
and $t>0$,
$x\in D$, we have

\[
\int_{D}p_{t}( x,\cdot) f\,d\mu=\mathcal{P}_{t}f(
x)
=\int_{D}q_{t}( x,\cdot) f\,d\mu.
\]
Applying this identity to function $f=p_{t}( y,\cdot) $
where $%
y\in D$, and using (\ref{t2}), we obtain%
\[
p_{2t}( x,y) =\int_{D}q_{t}( x,\cdot)
p_{t}(
y,\cdot) \,d\mu.
\]
Similarly, we have%
\[
q_{2t}( x,y) =\int_{D}p_{t}( y,\cdot)
q_{t}(
x,\cdot) \,d\mu,
\]
whence $p_{2t}( x,y) =q_{2t}( x,y) $.
\end{pf}

Following \cite{FOT}, page 67, a sequence $\{ F_{n}\}
_{n=1}^{\infty
}$ of subsets of $M$ will be called a \textit{regular nest} if:

\begin{longlist}[(2)]
\item[(1)] each $F_{n}$ is closed;\vadjust{\goodbreak}

\item[(2)] $F_{n}\subset F_{n+1}$ for all $n\geq1$;

\item[(3)] $\limfunc{Cap}(M\setminus F_{n})\rightarrow0$ as $%
n\rightarrow\infty$ (see \cite{FOT} for the definition of capacity);

\item[(4)] measure $\mu|_{F_{n}}$ has full support in $F_{n}$ (in the induced
topology of $F_{n}$).
\end{longlist}
\begin{definition}
\label{DefN}
A set $\mathcal{N}\subset M$ is called \textit{truly
exceptional} if:

\begin{longlist}[(2)]
\item[(1)] $\mathcal{N}$ is properly exceptional;

\item[(2)] $\mathcal{N}\supset\mathcal{N}_{0}$;

\item[(3)] there is a regular nest $\{ F_{n}\} $ in $M$ such
that $%
M\setminus\mathcal{N}=\tbigcup_{n=1}^{\infty}F_{n}$ and that the
function $%
\mathcal{P}_{t}f\vert_{F_{n}}$ is continuous for all
$f\in
\mathcal{B}L^{1}( M) $, $t>0$, and $n\in\mathbb{N}$.
\end{longlist}
\end{definition}

The conditions under which a truly exceptional set exists, will be discussed
later on. Let us mention some important consequences of the existence of
such a set.
\begin{lemma}
\label{LemPt<fi}Let $\mathcal{N}$ be a truly exceptional set. If, for
some $%
f\in\mathcal{B}L^{1}( M) $, \mbox{$t>0$}, and for an upper
semicontinuous function $\varphi\dvtx M\rightarrow(-\infty,+\infty]$, the
inequality%
\[
\mathcal{P}_{t}f( x) \leq\varphi( x)
\]
holds for $\mu$-a.a. $x\in M$, then it is true for all
$x\in
M\setminus\mathcal{N}$. Similarly, if $\psi\dvtx M\rightarrow\lbrack
-\infty
,+\infty)$ is a lower semicontinuous function and%
\[
\mathcal{P}_{t}f( x) \geq\psi( x)
\]
holds for $\mu$-a.a. $x\in M$, then it is true for all
$x\in
M\setminus\mathcal{N}$.
\end{lemma}
\begin{pf}
This proof is essentially the same as in \cite{FOT}, Theorem 2.1.2(ii).
Assume that $\mathcal{P}_{t}f( x_{0}) >\varphi(
x_{0}) $ for some $x_{0}\in M\setminus\mathcal{N}$. By
Definition \ref{DefN}, $x_{0}$ belongs to one of the sets $F_{n}$.
Since $\mathcal{P}%
_{t}f|_{F_{n}}$ is continuous and, hence, $( \mathcal
{P}_{t}f-\varphi
) |_{F_{n}}$ is lower semicontinuous, the condition $(
\mathcal{P}%
_{t}f-\varphi) ( x_{0}) >0$ implies that $(
\mathcal{P%
}_{t}f-\varphi) ( x) >0$ for all $x$ in some open
neighborhood $U$ of $x_{0}$ in $F_{n}$. Since measure $\mu$ has full
support in $F_{n}$, we have $\mu( U) >0$ so that
$\mathcal{P}%
_{t}f( x) >\varphi( x) $ in a set of positive
measure, that contradicts the hypothesis.

The second claim follows from the first one with
$\varphi=-\psi$.
\end{pf}

Denote by $\limfunc{esup}_{A}f$ the $\mu$-essential supremum of a
function $%
f$ on a set $A\subset M$, and by $\limfunc{einf}_{A}f$---the $\mu$%
-essential infimum.
\begin{corollary}
\label{Cesup=sup1}Let $\mathcal{N}$ be a truly exceptional set. Then, for
any $f\in\mathcal{B}L^{1}( M) $, $t>0$, and an open set $
X\subset M$,%
%
%e2.11 ###
%
\begin{equation} \label{esi}
\limfunc{esup}_{X}\mathcal{P}_{t}f=\sup_{X\setminus\mathcal
{N}}\mathcal{P}%
_{t}f \quad\mbox{and}\quad\limfunc{einf}_{X}\mathcal{P}_{t}f=\inf
_{X\setminus
\mathcal{N}}\mathcal{P}_{t}f.
\end{equation}
\end{corollary}
\begin{pf}
Function%
\[
\varphi( x) =\cases{
\displaystyle \limfunc{esup}_{X}\mathcal{P}_{t}f, &\quad$x\in X$, \vspace*{2pt}\cr
+\infty, &\quad$x\notin X$,}\vadjust{\goodbreak}
\]
is upper semicontinuous. Since $\mathcal{P}_{t}f( x) \leq
\varphi
( x) $ for $\mu$-a.a. $x\in M$, we conclude
by Lemma %
\ref{LemPt<fi} that this inequality is true for all $x\in M\setminus
\mathcal{N}$, whence%
\[
\sup_{X\setminus\mathcal{N}}\mathcal{P}_{t}f\leq\limfunc
{esup}_{X}\mathcal{%
P}_{t}f.
\]
The opposite inequality follows trivially from the definition of the
essential supremum.

The second identity in (\ref{esi}) follows from the first one by
changing $f$
to~$-f$.
\end{pf}

Note that if $p_{t}( x,y) $ is the heat kernel with domain
$%
D=M\setminus\mathcal{N}$, then we have by (\ref{semi}) that, for all
$x,y\in
D$, $0<s<t$,%
%
%e2.12 ###
%
\begin{equation}\label{ptPt}
p_{t}( x,y) =\mathcal{P}_{s}f( x) ,
\end{equation}
where $f=p_{t-s}( \cdot,y) $. Hence, if $\mathcal{N}$ is truly
exceptional, then the claims of Lem\-ma~\ref{LemPt<fi} and Corollary
\ref{Cesup=sup1} apply to function $p_{t}( x,y) $ in place of $
\mathcal{P}_{t}f( x) $, with any fixed $y\in D$.
\begin{lemma}
\label{LemBGK}
Let $p_{t}(x,y)$ be the heat kernel with the
domain $%
D=M\setminus\mathcal{N}$ such that $\mathcal{N}$ is a truly exceptional
set. Let $\varphi\dvtx D\times D\rightarrow[ 0,+\infty] $ be an
upper semicontinuous function and $\psi\dvtx D\times D\rightarrow\lbrack
0,+\infty)$ be a lower semicontinuous function. If, for some fixed $t>0$,
the following inequality:
%
%e2.13 ###
%
\begin{equation}\label{ptA}
\psi( x,y) \leq p_{t}(x,y)\leq\varphi(x,y)
\end{equation}
holds for $\mu\times\mu$-almost all $x,y\in D$, then (\ref{ptA})
holds for all $x,y\in D$.
\end{lemma}

This lemma is a generalization of
\cite{BarGrigKum}, Lemma 2.2, and
the proof follows the argument in \cite{BarGrigKum}.
\begin{pf*}{Proof of Lemma \ref{LemBGK}}
Consider the set%
\[
D^{\prime}=\{ y\in D\dvt\mbox{(\ref{ptA}) holds for }\mu\mbox
{-a.a. }x\in D\} .
\]
If $y\in D^{\prime}$ then applying Lemma \ref{LemPt<fi} to the
function $%
p_{t}( \cdot,y) $, we obtain that~(\ref{ptA}) holds for
all $%
x\in D$.

Now fix $x\in D$. Since by Fubini's theorem $\mu(D\setminus D^{\prime
})=0$,~(\ref{ptA}) holds for $\mu$-a.a. $y\in M$. Applying
Lemma~\ref{LemPt<fi} to the function $p_{t}( x,\cdot) $, we
conclude that~(\ref{ptA}) holds for all $y\in D$.
\end{pf*}
\begin{corollary}
\label{Cesup=sup}Under the hypotheses of Lemma \ref{LemBGK},
if $X,Y$
are two open subsets of $M$ then%
%
%e2.14 ###
%
\begin{equation}\label{esup=sup}
\mathop{\limfunc{esup}_{x\in X}}_{y\in Y}p_{t}( x,y)
=\mathop{\sup
_{x\in X\setminus\mathcal{N}}}_{y\in Y\setminus\mathcal{N}}%
p_{t}( x,y)
\end{equation}
and%
%
%e2.15 ###
%
\begin{equation}\label{einf=inf}
\mathop{\limfunc{einf}_{x\in X}}_{y\in Y}p_{t}( x,y)
=\mathop{\inf
_{x\in X\setminus\mathcal{N}}}_{y\in Y\setminus\mathcal{N}}%
p_{t}( x,y) .
\end{equation}
\end{corollary}
\begin{pf}
This follows from Lemma \ref{LemBGK} with functions
\[
\varphi( x,y) =\cases{
\func{const}, &\quad$x\in X,y\in Y$, \cr
+\infty, &\quad otherwise,}%
\]
and%
\[
\psi( x,y) =\cases{
\func{const}, &\quad$x\in X,y\in Y$, \cr
0, &\quad otherwise.}%
\]
\upqed\end{pf}

In conclusion of this section, let us state a result that ensures the
existence of the heat kernel outside a truly exceptional set.
\begin{theorem}[(\cite{BBCK}, Theorem 2.1)]
\label{TBBCK} Assume that there is a
positive left-continuous function $\gamma( t) $ such that for
all $f\in L^{1}\cap L^{2}( M) $ and $t>0$,%
%
%e2.16 ###
%
\begin{equation} \label{ultra}
\Vert{}P_{t}f\Vert_{\infty}\leq\gamma( t) \Vert
{}f\Vert
_{1}.
\end{equation}
Then the transition semigroup $\mathcal{P}_{t}$ possesses the heat
kernel $%
p_{t}( x,y) $ with domain $D=M\setminus\mathcal{N}$ for some
truly exceptional set $\mathcal{N}$, and $p_{t}( x,y)
\leq\gamma
( t) $ for all $x,y\in D$ and $t>0$.
\end{theorem}

If the semigroup $\{ P_{t}\} $ satisfies (\ref{ultra}),
then it is
called \textit{ultracontractive} (cf.~\cite{Davbook}). It was proved in
\cite{BBCK} that the ultracontractivity implies the existence of a
function $
p_{t}( x,y) $ that satisfies all the requirements of
Definition~\ref{DefHK} except for the joint measurability in $x,y$. Let us prove the
latter so that $p_{t}( x,y) $ is indeed the heat kernel in our
strict sense. Given that $p_{t}( x,y) $ satisfies conditions
(2)--(4) of Definition~\ref{DefHK}, we see that the statement of Lemma~\ref
{Lemuniq} remains true because the proof of that lemma does not use
the
joint measurability.\vspace*{1pt} In particular, for any $t>0$, $x\in D$, the
function $%
p_{t}( x,\cdot) $ is in $L^{2}( M) $. Also, the
mapping $x\mapsto p_{t}( x,\cdot) $ is weakly measurable
as a
mapping from $D$ to $L^{2}( M) $ because for any $f\in
L^{2}( M) $, the function $x\mapsto(
p_{t,x},f) =%
\mathcal{P}_{t}f(x)$ is measurable. Since
$L^{2}(
M) $ is separable, by Pettis's measurability theorem (see
\cite{Yosida}, Chapter V,\vspace*{1pt} Section 4) the mapping $x\mapsto p_{t}(
x,\cdot)
$ is
strongly measurable in $L^{2}( M)$. It follows that the
function%
\[
p_{2t}( x,y) =( p_{t}( x,\cdot)
,p_{t}(
y,\cdot) )
\]
is jointly measurable in $x,y\in D$ as the composition of two strongly
measurable mappings $D\rightarrow L^{2}( M) $ and a continuous
mapping $ f,g\mapsto( f,g)$.

%s2.3 ###
\subsection{Restricted heat semigroup and local ultracontractivity}

Any open subset $\Omega$ of $M$ can be considered as a metric measure space
$( \Omega,d,\mu) $. Let us identify\vadjust{\goodbreak} $L^{2}( \Omega
)
$ as a subspace in $L^{2}( M) $ by extending functions
outside~$\Omega$ by~$0$. Define~$\mathcal{F}( \Omega) $ as the closure
of $\mathcal{F}\cap C_{0}( \Omega) $ in $\mathcal{F}$
so that $%
\mathcal{F}( \Omega) $ is a~subspace of both $\mathcal
{F}$ and $L^{2}( \Omega) $. Then $( \mathcal
{E},\mathcal{F}%
( \Omega) ) $ is a regular strongly local
Dirichlet form
in $L^{2}( \Omega) $, which is called the restriction of
$(
\mathcal{E},\mathcal{F}) $ to $\Omega$. Let~$\mathcal{L}^{\Omega}$
be the generator of the form $( \mathcal{E},\mathcal{F}(
\Omega
) ) $ and $P_{t}^{\Omega}=e^{-t\mathcal{L}^{\Omega}}$,
$t\geq
0 $, be the \textit{restricted} heat semigroup.

Define the \textit{first exit time} from $\Omega$ by%
\[
\tau_{\Omega}=\inf\{ t>0\dvtx X_{t}\notin\Omega\} .
\]
The diffusion process associated with the restricted Dirichlet form,
can be
canonically obtained from $\{ X_{t}\} $ by killing the latter
outside $\Omega$, that is, by restricting the life time of $X_{t}$ by
$\tau
_{\Omega}$ (see \cite{FOT}). It follows that the transition operator $
\mathcal{P}_{t}^{\Omega}$ of the killed diffusion is given by%
%
%e2.17 ###
%
\begin{equation}\label{PtOm}
\mathcal{P}_{t}^{\Omega}f( x) =\mathbb{E}_{x}\bigl(
\mathbf{1}%
_{\{ t<\tau_{\Omega}\} }f( X_{t}) \bigr)
\qquad \mbox{for all }x\in\Omega\setminus\mathcal{N}_{0},
\end{equation}
for all $f\in\mathcal{B}_{+}( \Omega) $. Then $\mathcal
{P}%
_{t}^{\Omega}f$ is defined for $f$ from other function classes in the same
way as $\mathcal{P}_{t}$. Also, extend $\mathcal{P}_{t}^{\Omega
}f(
x) $ to all $x\in\Omega$ by setting it to be $0$ if $x\in
\mathcal{N}%
_{0}$.
\begin{definition}
We say that the semigroup $P_{t}$ is \textit{locally ultracontractive}
if the
restricted heat semigroup $P_{t}^{B}$ is ultracontractive for any metric
ball $B$ of $( M,d) $.
\end{definition}
\begin{theorem}
\label{TptOm}Let the semigroup $P_{t}$ be locally ultracontractive.
Then the
following is true.

\begin{longlist}[(a)]
\item[(a)] There exists a properly exceptional set
$\mathcal{N%
}\subset M$ such that, for any open subset $\Omega\subset M$, the semigroup
$\mathcal{P}_{t}^{\Omega}$ possesses the heat kernel $p_{t}^{\Omega
}(
x,y) $ with the domain $\Omega\setminus\mathcal{N}$.

\item[(b)] If $\Omega_{1}\subset\Omega_{2}$ are open
subsets of $M$, then $p_{t}^{\Omega_{1}}( x,y) \leq
p_{t}^{\Omega
_{2}}( x,y) $ for all $t>0$, $x,y\in\Omega_{1}\setminus
\mathcal{N}$.

\item[(c)] If $\{ \Omega_{k}\}
_{k=1}^{\infty}$
is an increasing sequence of open subsets of $M$ and $\Omega
=\tbigcup_{k}\Omega_{k}$, then $p_{t}^{\Omega_{k}}( x,y)
\rightarrow p_{t}^{\Omega}( x,y) $ as $k\rightarrow
\infty$ for
all $t>0$, $x,y\in\Omega\setminus\mathcal{N}$.

\item[(d)] Set $D=M\setminus\mathcal{N}$. Let
$\varphi
( x,y) \dvtx D\times D\rightarrow[ 0,+\infty] $
be an
upper semi-continuous function such that, for some open set $\Omega
\subset
M $ and for some $t>0$,%
%
%e2.18 ###
%
\begin{equation}\label{ptfi}
p_{t}^{\Omega}( x,y) \leq\varphi( x,y)
\end{equation}
for almost all $x,y\in\Omega$. Then (\ref{ptfi}) holds for
all $%
x,y\in\Omega\setminus\mathcal{N}$.
\end{longlist}
\end{theorem}

For simplicity of notation, set $p_{t}^{\Omega}( x,y) $
to be $0$
for all $x,y$ outside $\Omega$ (which, however, does not mean the extension
of the domain of $p_{t}^{\Omega}$).
\begin{pf*}{Proof of Theorem \ref{TptOm}}
(a) Since the metric space $( M,d) $ is separable,
there is a countable family of balls that form a base. Let $\mathcal
{U}$ be
the family of all finite unions of such balls so that $\mathcal{U}$ is
countable and any open set $\Omega\subset M$ can be represented\vadjust{\goodbreak} as an
increasing union of sets of $\mathcal{U}$. Since any set $U\in
\mathcal{U}$
is contained in a metric ball, the semigroup $P_{t}^{U}$ is dominated
by~$P_{t}^{B}$ and, hence, is ultracontractive. By Theorem \ref{TBBCK},
there is
a truly exceptional set $\mathcal{N}_{U}\subset U$ such that the
$\mathcal{P}%
_{t}^{U}$ has the heat kernel $p_{t}^{U}$ in the domain $U\setminus
\mathcal{%
N}_{U}$. Since the family $\mathcal{U}$ is countable, the set%
%
%e2.19 ###
%
\begin{equation}\label{N=}
\mathcal{N}=\tbigcup_{U\in\mathcal{U}}\mathcal{N}_{U}
\end{equation}
is properly exceptional.

Let us first show that if $U_{1}$, $U_{2}$ are the sets from $\mathcal{U}$
and $U_{1}\subset U_{2}$, then
%
%e2.20 ###
%
\begin{equation} \label{pt12}
p_{t}^{U_{1}}( x,y) \leq p_{t}^{U_{2}}( x,y)
\qquad\mbox{for all }t>0,x,y\in U_{1}\setminus\mathcal{N}.
\end{equation}
It follows from (\ref{PtOm}) that, for any $f\in\mathcal{B}_{+}(
U_{1}) $,%
\[
\mathcal{P}_{t}^{U_{1}}f( x) \leq\mathcal
{P}_{t}^{U_{2}}f(
x) \qquad\mbox{for all }t>0\mbox{ and }x\in U_{1},
\]
that is,%
%
%e2.21 ###
%
\begin{equation}\label{Om12}
\int_{U_{1}}p_{t}^{U_{1}}( x,\cdot) f\,d\mu\leq
\int_{U_{2}}p_{t}^{U_{2}}( x,\cdot) f\,d\mu.
\end{equation}
Setting here $f=P_{t}^{U_{1}}( y,\cdot) $ where $y\in
U_{1}\setminus\mathcal{N}$, we obtain%
\[
p_{2t}^{U_{1}}( x,y) \leq\int_{U_{1}}p_{t}^{U_{2}}
( x,\cdot
) p_{t}^{U_{1}}( y,\cdot) \,d\mu.
\]
Setting in (\ref{Om12}) $f=P_{t}^{U_{2}}( y,\cdot) $, we
obtain%
\[
\int_{U_{1}}p_{t}^{U_{1}}( x,\cdot) P_{t}^{U_{2}}(
y,\cdot
) \,d\mu\leq p_{2t}^{U_{2}}( x,y) .
\]
Combining the above two lines gives (\ref{pt12}).

Let $\Omega$ be any open subset of $M$ and $\{ U_{n}\}
_{n=1}^{\infty}$ be an increasing sequence of sets from $\mathcal{U}$ such
that $\Omega=\tbigcup_{n=1}^{\infty}U_{n}$. Let us set%
%
%e2.22 ###
%
\begin{equation}\label{ptOmlim}
p_{t}^{\Omega}( x,y) =\lim_{n\rightarrow\infty
}p_{t}^{U_{n}}( x,y) \qquad\mbox{for all }t>0\mbox{ and
}x,y\in
\Omega\setminus\mathcal{N}.
\end{equation}
This limit exists (finite or infinite) by the monotonicity of the
sequence $%
\{ p_{t}^{U_{n}}\hspace*{-0.2pt}( x,y) \} $. It follows from
(\ref{PtOm}) that, for any $f\in\mathcal{B}_{+}( \Omega) $,%
\[
\mathcal{P}_{t}^{U_{n}}f( x) \uparrow\mathcal
{P}_{t}^{\Omega
}f( x) \qquad\mbox{for all }t>0\mbox{ and }x\in\Omega
\setminus
\mathcal{N}.
\]
By the monotone convergence theorem, we obtain%
\[
\mathcal{P}_{t}^{U_{n}}f( x) =\int_{\Omega
}p_{t}^{U_{n}}(
x,y) f( y) \,d\mu( y) \rightarrow\int
_{\Omega
}p_{t}^{\Omega}( x,y) f( y) \,d\mu(
y)
\]
for all $t>0$ and $x\in\Omega\setminus\mathcal{N}$. Comparing the above
two lines, we obtain%
\[
\mathcal{P}_{t}^{\Omega}f( x) =\int_{\Omega
}p_{t}^{\Omega
}( x,y) f( y) \,d\mu( y) \qquad\mbox
{for all }t>0%
\mbox{ and }x\in\Omega\setminus\mathcal{N}.
\]
The symmetry of $p_{t}^{\Omega}( x,y) $ is obvious from
(\ref{ptOmlim}), and the semigroup property of $p_{t}^{\Omega}$ follows from
that of $p_{t}^{U_{n}}$ by the monotone convergence theorem. Note that $
p_{t}^{\Omega}$ does not depend on the choice of $\{ U_{n}
\} $
by the uniqueness of the heat kernel (Lemma \ref{Lemuniq}).

(b) For two arbitrary open sets $\Omega_{1}\subset
\Omega
_{2}$ let $\{U_{n}\}_{n=1}^{\infty}$ and $\{W_{n}\}_{n=1}^{\infty}$ be
increasing sequences of sets from $\mathcal{U}$ that exhaust $\Omega_{1}$
and $\Omega_{2}$, respectively. Set $V_{n}=U_{n}\cup W_{n}$ so that $%
V_{n}\in\mathcal{U}$ and $\Omega_{2}$ is the\vspace*{1pt} increasing union of
sets $%
V_{n}$ (see Figure \ref{pic6}). Then $U_{n}\subset V_{n}$ and, hence, $
p_{t}^{U_{n}}\leq p_{t}^{V_{n}}$, which implies as $n\rightarrow\infty$
that $p_{t}^{\Omega_{1}}\leq p_{t}^{\Omega_{2}}$.

%
%f1 ###
%
\begin{figure}

\includegraphics{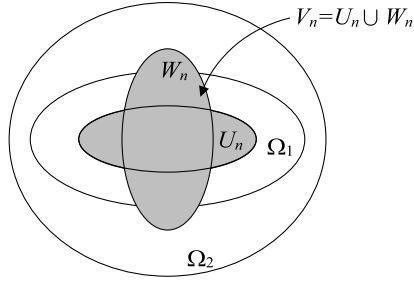}

\caption{Sets $U_{n},W_{n},V_{n}$.}\label{pic6}
\end{figure}

(c) Let $\{ \Omega_{k}\} _{k=1}^{\infty
}$ be an
increasing sequence of open sets whose union is~$\Omega$. Let $%
\{U_{n}^{( k) }\}_{n=1}^{\infty}$ be an increasing
sequence of
sets from $\mathcal{U}$ that exhausts~$\Omega_{k}$. As in the previous
argument,\vspace*{1pt} we can replace $U_{n}^{( 2) }$ by $V_{n}^{(
2) }=U_{n}^{( 1) }\cup U_{n}^{( 2) }$
so that $%
U_{n}^{( 1) }\subset V_{n}^{( 2) }$. Rename $%
V_{n}^{( 2) }$ back to $U_{n}^{( 2) }$, and
assume\vspace*{1pt} in
the sequel that $U_{n}^{( 1) }\subset U_{n}^{(
2) }$.
Similarly, replace $U_{n}^{( 3) }$ by $U_{n}^{(
1)
}\cup U_{n}^{( 2) }\cup U_{n}^{( 3) }$ and
assume in
the sequel that $U_{n}^{( 2) }\subset U_{n}^{(
3) }$.
Arguing by induction, we redefine the double sequence $U_{n}^{(
k) }$ in the way that it is monotone increasing not only in $n$ but
also in $k$. Then we claim that
\[
\Omega=\tbigcup_{m=1}^{\infty}U_{m}^{( m) }.
\]
Indeed, if $x\in\Omega, $ then $x\in\Omega_{k}$ for some $k$ and,
hence, $%
x\in U_{n}^{( k) }$ for some $n$, which implies $x\in
U_{m}^{( m) }$ for $m=\max( k,n) $. Finally,
we have $%
p_{t}^{\Omega}\geq p_{t}^{\Omega_{m}}$ and%
\[
p_{t}^{\Omega}=\lim_{m\rightarrow\infty}p_{t}^{U_{m}^{(
m)
}}\leq\lim_{m\rightarrow\infty}p^{\Omega_{m}},
\]
whence it follows that%
\[
p_{t}^{\Omega}=\lim_{m\rightarrow\infty}p^{\Omega_{m}}.
\]

(d) Let $U\in\mathcal{U}$ be subset of $\Omega$.
Then the
semigroup $P_{t}^{U}$ is ultracontractive and possesses the heat kernel
$%
p_{t}^{U}$ with the domain $U\setminus\mathcal{N}_{U}$ where
$\mathcal{N}%
_{U}$ is a~truly exceptional\vadjust{\goodbreak} set as in part (a). Note
that $%
\mathcal{N}_{U}\subset\mathcal{N}$. Since $p_{t}^{U}\leq
p_{t}^{\Omega}$
in $U\setminus\mathcal{N}$, we obtain by hypothesis that
\[
p_{t}^{U}( x,y) \leq\varphi( x,y)
\]
for almost all $x,y\in U$. By Lemma \ref{LemBGK}, we conclude that this
inequality is true for all $x,y\in U\setminus\mathcal{N}$. Exhausting
$%
\Omega$ be a sequence of subsets $U\in\mathcal{U}$ and using (\ref
{ptOmlim}%
), we obtain (\ref{ptfi}).
\end{pf*}

%s3 ###
\section{Some preparatory results}
\label{Secaux}

%s3.1 ###
\subsection{Green operator}
\label{SecFK}

A priori we assume here only the basic hypotheses. All
necessary additional assumptions are explicitly stated. The main result of
this section is Theorem \ref{TG=>FK}.

Given an open set $\Omega\subset M$, define the Green operator
$G^{\Omega}$
first for all $f\in\mathcal{B}_{+}( \Omega) $ by%
%
%e3.1 ###
%
\begin{equation}\label{Gdef}
G^{\Omega}f( x) =\int_{0}^{\infty}\mathcal
{P}_{t}^{\Omega
}f ( x) \,dt
\end{equation}
for all $x\in M\setminus\mathcal{N}_{0}$, where we admit infinite
values of
the integral. If $f$ $\in\mathcal{B}( \Omega) $ and
$G^{\Omega
}\vert f\vert<\infty, $ then $G^{\Omega}f$ is also
defined by
\[
G^{\Omega}f=G^{\Omega}f_{+}-G^{\Omega}f_{-} .
\]

\begin{lemma}
We have, for any open $\Omega\subset M$ and all $f\in\mathcal
{B}_{+}(
\Omega) $,%
%
%e3.2 ###
%
\begin{equation} \label{GOm}
G^{\Omega}f( x) =\mathbb{E}_{x}\biggl( \int_{0}^{\tau
_{\Omega
}}f( X_{t}) \,dt\biggr)
\end{equation}
for any $x\in\Omega\setminus\mathcal{N}_{0}$. In particular,%
%
%e3.3 ###
%
\begin{equation}\label{Gtau}
G^{\Omega}1( x) =\mathbb{E}_{x}\tau_{\Omega}.
\end{equation}
\end{lemma}
\begin{pf}
Indeed, integrating (\ref{PtOm}) in $t$, we obtain%
\begin{eqnarray*}
G^{\Omega}f( x) &=&\int_{0}^{\infty}\mathcal
{P}_{t}^{\Omega
}f( x) \,dt \\
&=&\int_{0}^{\infty}\mathbb{E}_{x}\bigl( \mathbf{1}_{\{
t<\tau
_{\Omega}\} }f( X_{t}) \bigr) \,dt \\
&=&\mathbb{E}_{x}\int_{0}^{\infty}\bigl( \mathbf{1}_{\{
t<\tau
_{\Omega}\} }f( X_{t}) \bigr) \,dt \\
&=&\mathbb{E}_{x}\biggl( \int_{0}^{\tau_{\Omega}}f(
X_{t})
\,dt\biggr) .
\end{eqnarray*}
Obviously, (\ref{Gtau}) follows from (\ref{GOm}) for
$f\equiv1$.
\end{pf}

Denote by $\lambda_{\min}( \Omega) $ the bottom of the
spectrum of $\mathcal{L}^{\Omega}$ in $L^{2}( \Omega)
$, that is,%
%
%e3.4 ###
%
\begin{equation}\label{ladef}
\lambda_{\min}( \Omega) :=\inf\func{spec}\mathcal
{L}^{\Omega
}=\inf_{f\in\mathcal{F}( \Omega) \setminus\{
0\} }%
\frac{\mathcal{E}( f,f) }{( f,f) }.\vadjust{\goodbreak}
\end{equation}
For any open set $\Omega\subset M$, we will consider the \textit{mean exit
time} $\mathbb{E}_{x}\tau_{\Omega}$ from $\Omega$ as a function of
$x\in
\Omega\setminus\mathcal{N}_{0}$. Also, set%
%
%e3.5 ###
%
\begin{equation}\label{Ebar}
\widetilde{E}( \Omega) :=\limfunc{esup}_{x\in\Omega
}\mathbb{E}%
_{x}\tau_{\Omega}.
\end{equation}

\begin{lemma}
\label{LG1-1}If $\widetilde{E}( \Omega) <\infty, $
then $%
G^{\Omega}$ is a bounded operator on $\mathcal{B}_{b}( \Omega
), $
and it uniquely extends to each of the spaces $L^{\infty}(
\Omega
) $, $L^{1}( \Omega) $, $L^{2}( \Omega
) $,
with the following norm estimates:%
%
%e3.8 ###
%e3.7 ###
%e3.6 ###
%
\begin{eqnarray}\label{Gi-i}
\Vert G^{\Omega}\Vert_{L^{\infty}\rightarrow L^{\infty}}&\leq&
\widetilde{E%
}( \Omega) ,
\\[-1pt]
\label{G1-1}
\Vert G^{\Omega}\Vert_{L^{1}\rightarrow L^{1}}&\leq&\widetilde
{E}(
\Omega) ,
\\[-1pt]
\label{G2-2}
\Vert G^{\Omega}\Vert_{L^{2}\rightarrow L^{2}}&\leq&\widetilde
{E}(
\Omega) .
\end{eqnarray}
Moreover,
%
%e3.9 ###
%
\begin{equation}\label{lamin<}
\lambda_{\min}( \Omega) ^{-1}\leq\widetilde{E}(
\Omega) ,
\end{equation}
and $G^{\Omega}$ is the inverse in $L^{2}( \Omega) $ to the
operator $\mathcal{L}^{\Omega}$.
\end{lemma}

\begin{pf}
It follows from (\ref{Gtau}) that%
%
%e3.10 ###
%
\begin{equation}\label{GOm1}
\Vert G^{\Omega}1\Vert_{\infty}=\widetilde{E}( \Omega) ,
\end{equation}
which implies that for any $f\in\mathcal{B}_{b}( \Omega
) $,
\[
\Vert G^{\Omega}f\Vert_{\infty}\leq\widetilde{E}( \Omega
)
\Vert f\Vert_{\infty}.
\]
Hence, $G^{\Omega}$ can be considered as a bounded operator in
$L^{\infty}$
with the norm estimate (\ref{Gi-i}).

Estimate (\ref{G1-1}) follows from (\ref{Gi-i}) by duality. Indeed, for
any two functions $f,h\in\mathcal{B}_{+}( \Omega) $, we
have%
%
%e3.11 ###
%
\begin{equation}\label{Gfh}
\int_{\Omega}( G^{\Omega}f ) h\,d\mu=\int_{\Omega
}fG^{\Omega
}h \,d\mu,
\end{equation}
which follows from (\ref{Gdef}) and the symmetry of $\mathcal
{P}_{t}^{\Omega
}$. By linearity, (\ref{Gfh}) extends to all $f,h\in\mathcal
{B}_{b}(
\Omega) $. Then, for any $f\in C_{0}( \Omega) $,
we have%
\begin{eqnarray*}
\Vert G^{\Omega}f\Vert_{1} &=&\sup_{h\in\mathcal{B}_{b}(
\Omega
) \setminus\{ 0\} }\frac{\int_{\Omega}(
G^{\Omega
}f) h\,d\mu}{\Vert h\Vert_{\infty}} \\[-1pt]
&=&\sup_{h\in\mathcal{B}_{b}( \Omega) \setminus
\{
0\} }\frac{\int_{\Omega}fG^{\Omega}h\,d\mu}{\Vert h\Vert
_{\infty}}
\\[-1pt]
&\leq&\sup_{h\in\mathcal{B}_{b}( \Omega) \setminus
\{
0\} }\frac{\Vert G^{\Omega}h\Vert_{\infty}\Vert f\Vert
_{1}}{\Vert
h\Vert_{\infty}} \\[-1pt]
&\leq&\widetilde{E}( \Omega) \Vert f\Vert_{1},
\end{eqnarray*}
whence it follows that $G^{\Omega}$ uniquely extends to a bounded operator
in $L^{1}$ with the norm estimate (\ref{G1-1}).\vadjust{\goodbreak}

The estimate (\ref{G2-2}) [as well as a similar estimate for $
\Vert
G\Vert_{L^{p}\rightarrow L^{p}}$ for any $p\in( 1,\infty
) $] follows from (\ref{Gi-i}) and (\ref{G1-1}) by the Riesz--Thorin
interpolation theorem.

To prove (\ref{lamin<}), let us consider the following
``cut-down'' version of the Green operator:%
\[
G_{T}^{\Omega}f=\int_{0}^{T}\mathcal{P}_{t}^{\Omega}f \,dt,
\]
where $T\in(0,+\infty)$. The same argument as above shows that $%
G_{T}^{\Omega}$ can be considered as an operator in $L^{2}$ with the same
norm bound%
\[
\Vert G_{T}^{\Omega}\Vert_{L^{2}\rightarrow L^{2}}\leq\widetilde
{E}(
\Omega) .
\]
On the other hand, using the spectral resolution $\{ E_{\lambda
}\} _{\lambda\geq0}$ of the generator~$\mathcal{L}^{\Omega
}$, we
obtain, for any $f\in C_{0}( \Omega) $,%
%
%e3.12 ###
%
\begin{eqnarray} \label{fiT}
G_{T}^{\Omega}f &=&\int_{0}^{T}\biggl( \int_{0}^{\infty}e^{-\lambda
t}\,dE_{\lambda}f\biggr) \,dt \nonumber\\
&=&\int_{0}^{\infty}\biggl( \int_{0}^{T}e^{-\lambda t}\,dt\biggr)\,
dE_{\lambda
}f \nonumber\\[-8pt]\\[-8pt]
&=&\int_{0}^{\infty}\varphi_{T}( \lambda) \,dE_{\lambda}f
\nonumber\\
&=&\varphi_{T}( \mathcal{L}^{\Omega}) f,\nonumber
\end{eqnarray}
where%
\[
\varphi_{T}( \lambda) =\int_{0}^{T}e^{-\lambda
t}\,dt=\frac{%
1-e^{-T\lambda}}{\lambda}.
\]
Since $\varphi_{T}$ is a bounded function on $[0,+\infty)$, the
operator $%
\varphi_{T}( \mathcal{L}^{\Omega}) $ is a~bounded
operator in $%
L^{2}$. By the spectral mapping theorem, we obtain%
\begin{eqnarray*}
\sup\varphi_{T}( \func{spec}\mathcal{L}^{\Omega})
&=&\sup
\func{spec}\varphi_{T}( \mathcal{L}^{\Omega}) \\
&=&\Vert\varphi_{T}( \mathcal{L}^{\Omega}) \Vert
_{L^{2}\rightarrow L^{2}} \\
&=&\Vert G_{T}^{\Omega}\Vert_{L^{2}\rightarrow L^{2}} \\
&\leq&\widetilde{E}( \Omega) .
\end{eqnarray*}
On the other hand, since $\varphi_{T}( \lambda) $ is decreasing
in $\lambda$,%
\[
\sup\varphi_{T}( \func{spec}\mathcal{L}^{\Omega})
=\varphi
_{T}( \lambda_{\min}( \Omega) ) ,
\]
whence%
\[
\varphi_{T}( \lambda_{\min}( \Omega) )
\leq
\widetilde{E}( \Omega) .
\]
By letting $T\rightarrow\infty$ and observing that $\varphi_{T}(
\lambda) \rightarrow\frac{1}{\lambda}$, we obtain%
\[
\lambda_{\min}( \Omega) ^{-1}\leq\widetilde{E}(
\Omega
) ,\vadjust{\goodbreak}
\]
which in particular means that $\lambda_{\min}( \Omega) >0$.
Consequently, the operator $\mathcal{L}^{\Omega}$ has a bounded inverse.
Passing in (\ref{fiT}) to the limit as $T\rightarrow\infty$, we
obtain $%
G^{\Omega}=( \mathcal{L}^{\Omega}) ^{-1}$.
\end{pf}

%s3.2 ###
\subsection{Harmonic functions and Harnack inequality}
\label{SecHarnack}

Let $\Omega$ be an open subset of~$M$.
\begin{definition}
We say that a function $u\in\mathcal{F}$ is \textit{harmonic}
in $\Omega$ if%
\[
\mathcal{E}( u,v) =0 \qquad\mbox{for any }v\in\mathcal
{F}(
\Omega) .
\]
\end{definition}
\begin{lemma}
\label{LemG-G}Let $\Omega$ be an open subset of $M$ such that
$\widetilde{E}%
( \Omega) <\infty$, and let $U$ be an open subset of
$\Omega$.

\begin{longlist}[(a)]
\item[(a)]
For any $f\in L^{2}( \Omega\setminus
U) $,
the function $G^{\Omega}f$ is harmonic in $U$.\label{remdoweneeda?}

\item[(b)] For any $f\in L^{2}( \Omega) $, the
function $%
G^{\Omega}f-G^{U}f$ is harmonic in $U$.
\end{longlist}
\end{lemma}
\begin{remark}
If $f\in L^{2}( \Omega), $ then $G^{U}f$ is defined
as the
extension of $G^{U}( f|_{U}) $ to $\Omega$ by setting it
to be
equal to $0$ in $\Omega\setminus U$.
\end{remark}
\begin{pf*}{Proof of Lemma \ref{LemG-G}}
(a) Set $u=G^{\Omega}f$. To prove that $u$ is
harmonic in~$U$%
, we need to show that $\mathcal{E}( u,v) =0$, for any
$v\in
\mathcal{F}( U) $. Since by Lemma~\ref{LG1-1} $G^{\Omega
}=(
\mathcal{L}^{\Omega}) ^{-1}$, we have $u\in\func{dom}(
\mathcal{%
L}^{\Omega}) $. Therefore, by the definition of $\mathcal
{L}^{\Omega
} $,%
\[
\mathcal{E}( u,v) =( \mathcal{L}^{\Omega
}u,v) =(
f,v) =0.
\]

(b) Set $u=G^{\Omega}f-G^{U}f$. Any function $v\in
\mathcal{F%
}( U) $ can be considered as an element of $\mathcal
{F}(
\Omega) $ by setting it to be $0$ in $\Omega\setminus U$. Then
both $%
u$ and $v$ are in~$\mathcal{F}( \Omega) $ whence%
\begin{eqnarray*}
\mathcal{E}( u,v) &=&\mathcal{E}( G^{\Omega
}f,v) -%
\mathcal{E}( G^{U}f,v) \\
&=&( f,v) _{L^{2}( \Omega) }-(
f,v)
_{L^{2}( U) } \\
&=&0.
\end{eqnarray*}
\upqed\end{pf*}

Denote by
\[
B( x,r) =\{ y\in M\dvtx d( x,y) <r\}
\]
the open metric ball of radius $r>0$ centered at a point $x\in M$, and
set%
\[
V( x,r) =\mu( B( x,r) ) .
\]
That $\mu$ has full support implies $V( x,r) >0$ whenever $r>0$.
Whenever we use the function $V( x,r) $, we always assume that
\[
V( x,r) <\infty\qquad\mbox{for all }x\in M\mbox{ and }r>0.
\]
For example, this condition is automatically satisfied if all balls are
precompact. However, we do not assume precompactness of all balls unless
otherwise explicitly stated.\vadjust{\goodbreak}
\begin{definition}
We say that the \textit{elliptic Harnack inequality} (\ref{condH})
holds on~$M$, if there exist constants $C>1$ and $\delta\in(
0,1) $ such that, for any ball~$B( x,r) $ in $M$
and for
any function $u\in\mathcal{F}$ that is nonnegative and harmonic~$B(x,r) $,
{\renewcommand{\theequation}{$H$}
\begin{equation}\label{condH}
\limfunc{esup}_{B( x,\delta r) }u\leq C \mathop{\func
{einf}}_{B(
x,\delta r) }u .
\end{equation}}
\end{definition}

\begin{definition}
We say that the \textit{volume doubling} property
(\ref{VD})
holds if there exists a constant $C$ such that, for all $x\in M$ and
$r>0$%
\label{condVD}%
{\renewcommand{\theequation}{\textit{VD}}
\begin{equation}
V( x,2r) \leq CV( x,r) .
\end{equation}}
\end{definition}

It is known that (\ref{VD}) implies that, for all $x,y\in M$
and $%
0<r<R$,%
%
%e3.13 ###
%
\setcounter{equation}{12}
\begin{equation}\label{Va}
\frac{V( x,R) }{V( y,r) }\leq C\biggl( \frac
{R+d(
x,y) }{r}\biggr) ^{\alpha}
\end{equation}
for some $\alpha>0$ (see \cite{GrigHuUpper}).
\begin{lemma}
Assume that $\mbox{(\ref{VD})} +\mbox{(\ref{condH})} $ hold. Let $\Omega
$ be an
open subset of~$M$ such that $\widetilde{E}( \Omega)
<\infty, $
and let $B=B( x,r) $ be a ball contained in $\Omega$.

\begin{longlist}[(a)]
\item[(a)]
\label{LgOm}
For any nonnegative function
$\varphi\in
L^{1}( \Omega\setminus B) $, \label{remdoweneeda?copy1}%
%
%e3.14 ###
%
\begin{equation}\label{gOm<E}
\limfunc{esup}_{B( x,\delta r) }G^{\Omega}\varphi\leq
C\frac{%
\widetilde{E}( \Omega) }{V( x,r) }\Vert
{}\varphi
\Vert_{1}.
\end{equation}

\item[(b)]\label{LemGB-GB}For and any nonnegative function
$%
\varphi\in L^{1}( \Omega) $,%
%
%e3.15 ###
%
\begin{equation}\label{gBB}
\limfunc{esup}_{B( x,\delta r) }( G^{\Omega
}\varphi
-G^{B}\varphi) \leq\frac{C\widetilde{E}( \Omega
) }{%
V( x,r) }\Vert{}\varphi\Vert_{1}.
\end{equation}
\end{longlist}
\end{lemma}

\begin{pf}
(a) Since identity (\ref{gOm<E}) survives monotone
increasing limits of sequences of functions $\varphi$, it suffices to prove
(\ref{gOm<E}) for any nonnegative function $\varphi\in L^{1}\cap
L^{2}( \Omega\setminus B) $. Then, by Lemma \ref{LemG-G}, the
function $u=G^{\Omega}\varphi$ is harmonic in $B( x,r)$.\vadjust{\goodbreak} Since
$u\geq0$, we can use the Harnack inequality (\ref{condH}) in
ball $B$,
which yields%
%
%e3.16 ###
%
\begin{eqnarray}\label{fi1}
\limfunc{esup}_{B( x,\delta r) }u( x) &\leq
&C\limfunc{%
einf}_{B( x,\delta r) }u\leq\frac{C}{V( x,r
) }\Vert
u\Vert_{1} \nonumber\\
&\leq&\frac{C}{V( x,r) }\Vert G^{\Omega}\Vert
_{L^{1}\rightarrow L^{1}}\Vert\varphi\Vert_{1} \\
&\leq&\frac{C\widetilde{E}( \Omega) }{V(
x,r) }\Vert
\varphi\Vert_{1}.\nonumber
\end{eqnarray}

(b) Assume first that $\varphi\in L^{1}\cap
L^{2}(
\Omega) $. By Lemma \ref{LemG-G}, the function
$u=G^{\Omega
}\varphi-G^{B}\varphi$ is harmonic in $B( x,r) $. Since
$u\geq
0 $, applying for this function argument (\ref{fi1}), we obtain (\ref
{gBB}). An arbitrary nonnegative function $\varphi\in L^{1}(
\Omega
) $ can be represented as a sum in $L^{1}( \Omega)
$%
\[
\varphi=\sum_{k=0}^{\infty}\varphi_{k},
\]
where $\varphi_{k}:=( \varphi-k) _{+}\wedge1\in
L^{1}\cap
L^{\infty}( \Omega) $. Applying (\ref{gBB}) to each
$\varphi
_{k}$ and summing up, we obtain (\ref{gBB}) for $\varphi$.
\end{pf}

%s3.3 ###
\subsection{Faber--Krahn inequality and mean exit time}
\label{SecEF}

A classical theorem of Faber and Krahn says that for any
bounded open set $\Omega\subset\mathbb{R}^{n}$,
\[
\lambda_{\min}( \Omega) \geq\lambda_{\min}(B),
\]
where $B$ is a ball in $\mathbb{R}^{n}$ of the same volume as $\Omega
$. If
the radius of $B$ is $r$, then $\lambda_{\min}( B)
=\frac{c}{%
r^{2}}$ where $c$ is a positive constant depending only on $n$, which
implies that
%
%e3.17 ###
%
\begin{equation}\label{FKRn}
\lambda_{\min}( \Omega) \geq c\mu( \Omega
) ^{-2/n};
\end{equation}
cf. \cite{Chavbook,Chavnotes}. We refer to lower estimates of
$%
\lambda_{\min}( \Omega) $ via a function of $\mu
( \Omega
) $ as \textit{Faber--Krahn inequalities}. A more general type of a
Faber--Krahn inequality holds on a complete $n$-dimensional Riemannian
manifold $M$ of nonnegative Ricci curvature: for any bounded open set $
\Omega\subset M$ and for any ball $B$ of radius $r$ containing~$\Omega
$,%
\label{condFK0}%
%
%e3.18 ###
%
\begin{equation}\label{FKrel}
\lambda_{\min}( \Omega) \geq\frac{c}{r^{2}}\biggl(
\frac{\mu
( B) }{\mu( \Omega) }\biggr) ^{\nu},
\end{equation}
where $\nu=2/n$ and $c=c( n) >0$ (see \cite{GrigHar}).
Note that
(\ref{FKRn}) follows from~(\ref{FKrel}) (apart from the sharp value
of the
constant $c$) because in $\mathbb{R}^{n}$ we have $\mu( B
) =%
\func{const}r^{n}$.

It was proved in \cite{GrigHeat} that, on any complete Riemannian manifold,
\[
\mbox{(\ref{FKrel})}\Leftrightarrow\mbox{(\ref{VD})} + \mbox{(\textit{UE})},
\]
where (\textit{UE}) is here the upper bound of the heat kernel
in the
Li--Yau estimate (\ref{LiYau}). In Section \ref{SecDUE} we will derive a
general upper bound (\ref{UEE}) from a~set of hypotheses containing
a suitable version of (\ref{FKrel}). In this section, we will deduce a
Faber--Krahn inequality from the main hypotheses.

We fix from now on a function $F\dvtx( 0,\infty) \rightarrow
(
0,\infty) $ that is a continuous increasing bijection of
$(0,\infty)$
onto itself, such that, for all $0<r\leq R$,%
%
%e3.19 ###
%
\begin{equation}\label{Fb}
C^{-1}\biggl( \frac{R}{r}\biggr) ^{\beta}\leq\frac{F(
R) }{%
F( r) }\leq C\biggl( \frac{R}{r}\biggr) ^{\beta^{\prime}}
\end{equation}
for some constants $1<\beta\leq\beta^{\prime}$, $C>1$. In the
sequel we
will use the inverse function $\mathcal{R}=F^{-1}$. It follows from
(\ref{Fb}%
) that%
%
%e3.20 ###
%
\begin{equation} \label{Rb}
C^{-1}\biggl( \frac{T}{t}\biggr) ^{1/\beta^{\prime}}\leq\frac
{\mathcal{R}%
( T) }{\mathcal{R}( t) }\leq C\biggl( \frac
{T}{t}%
\biggr) ^{1/\beta}
\end{equation}
for all $0<t\leq T$.
\begin{definition}
We say that the Faber--Krahn inequality (\ref{FK}) holds
if, for
any ball $B$ in $M$ of radius $r$ and any open set $\Omega\subset
B$,\label{condFK}%
{\renewcommand{\theequation}{\textit{FK}}
\begin{equation}\label{FK}
\lambda_{\min}( \Omega) \geq\frac{c}{F( r
) }\biggl(
\frac{\mu( B) }{\mu( \Omega) }\biggr)
^{\nu}
\end{equation}}

\vspace*{-8pt}

\noindent with some positive constants $c,\nu$.
\end{definition}
\begin{definition}
We say that the mean exit time estimate (\ref{EF})
holds if,
for all $x\in M\setminus\mathcal{N}_{0}$ and $r>0$,\label{condEF}%
{\renewcommand{\theequation}{${E}_{F}$}
\begin{equation}\label{EFF}
C^{-1}F( r) \leq\mathbb{E}_{x}\tau_{B( x,r
) }\leq
CF( r)
\end{equation}}

\vspace*{-8pt}

\noindent
with some constant $C>1$.
\end{definition}

We denote by $({E}_{F}\mbox{$\leq$}) $ and $({E}_{F}\mbox{$\geq$}
) $ the
upper and lower bounds of $\mathbb{E}_{x}\tau_{B( x,r)
}$ in~(\ref{EFF}), respectively.
\begin{theorem}
\label{TG=>FK}The hypotheses $\mbox{(\ref{VD})} +\mbox{(\ref{condH})}
+({E}_{F}\mbox{$\leq$}) $ imply (\ref{FK}).
\end{theorem}
\begin{pf}
We have by (\ref{lamin<}) and (\ref{Gtau})%
%
%e3.21 ###
%
\setcounter{equation}{20}
\begin{equation}\label{laG}
\lambda_{\min}( \Omega) ^{-1}\leq\widetilde{E}(
\Omega
) =\func{esup}_{x\in\Omega}G^{\Omega}1_{\Omega}.
\end{equation}
It will be convenient to rename $R$ to $R/2$ and let the original ball $B$
be $B( z,R/2) $ and $\Omega\subset B( z,R/2)
$. Fix a
point $x\in\Omega$ so that $\Omega\subset B( x,R) $, consider
a numerical sequence $R_{k}=\delta^{k}R$, $k=0,1,2,\ldots,$ where
$\delta
$ is
the parameter from (\ref{condH}), and the balls $B_{k}=B(
x,R_{k}) $. We have%
\[
G^{\Omega}1_{\Omega}\leq G^{B_{0}}1_{\Omega}=\sum_{k=0}^{n-1}(
G^{B_{k}}-G^{B_{k+1}}) 1_{\Omega}+G^{B_{n}}1_{\Omega},
\]
where $n$ is to be chosen (see Figure \ref{pic3}), whence%
\[
\limfunc{esup}_{B_{n+1}}G^{\Omega}1_{\Omega}\leq\sum
_{k=0}^{n-1}\limfunc{%
esup}_{B_{k+2}}( G^{B_{k}}-G^{B_{k+1}}) 1_{\Omega
}+\limfunc{esup}%
_{B_{n}}G^{B_{n}}1_{\Omega}.
\]

%
%f2 ###
%
\begin{figure}

\includegraphics{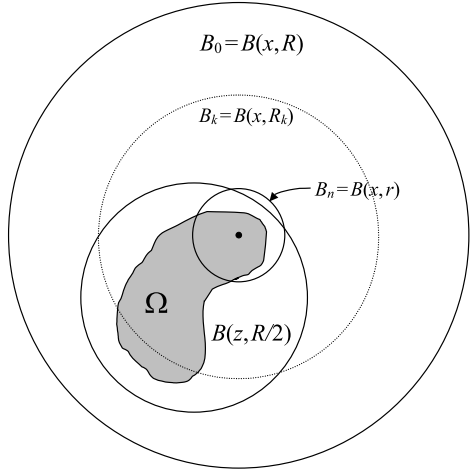}

\caption{Balls $B_{k}$.}\label{pic3}
\end{figure}

Setting $V( r) =V( x,r) $ and using
$\widetilde{E}%
( B_{k}) \leq F( R_{k}) $, we obtain, by
Lemma~\ref{LemGB-GB},%
\[
\limfunc{esup}_{B_{k+2}}( G^{B_{k}}-G^{B_{k+1}})
1_{\Omega}\leq
\frac{CF( R_{k}) }{V( R_{k}) }\mu(
\Omega).
\]
Also, by (\ref{GOm1}),%
\[
\limfunc{esup}_{B_{n}}G^{B_{n}}1_{\Omega}\leq\limfunc{esup}%
_{B_{n}}G^{B_{n}}1=\widetilde{E}( B_{n}) \leq CF(
R_{n}) .
\]
Hence, collecting together the previous lines, we obtain%
\[
\limfunc{esup}_{B_{n+1}}G^{\Omega}1_{\Omega}\leq C\sum
_{k=0}^{n-1}\frac{%
F( R_{k}) }{V( R_{k}) }\mu( \Omega
)
+CF( R_{n}) .
\]
Choose any $\nu\in( 0,1) $ so that $\nu<\beta/\alpha
$. Using
(\ref{Fb}), (\ref{Va}) and the monotonicity of $V( s) $,
we obtain%
\begin{eqnarray*}
\sum_{k=0}^{n-1}\frac{F( R_{k}) }{V( R_{k})
} &=&\frac{%
F( R) }{V( R_{n}) ^{1-\nu}V( R)
^{\nu}}%
\sum_{k=0}^{n-1}\frac{F( R_{k}) }{F( R)
}\biggl( \frac{%
V( R) }{V( R_{k}) }\biggr) ^{\nu}\biggl(
\frac{V(
R_{n}) }{V( R_{k}) }\biggr) ^{1-\nu} \\
&\leq&\frac{CF( R) }{V( R_{n}) ^{1-\nu
}V(
R) ^{\nu}}\sum_{k=0}^{n-1}\biggl( \frac{R_{k}}{R}\biggr)
^{\beta
}\biggl( \frac{R}{R_{k}}\biggr) ^{\alpha\nu} \\
&=&\frac{CF( R) }{V( R_{n}) ^{1-\nu}V
( R)
^{\nu}}\sum_{k=0}^{n-1}\delta^{k( \beta-\alpha\nu) }
\\
&\leq&\frac{CF( R) }{V( R_{n}) ^{1-\nu
}V(
R) ^{\nu}}.
\end{eqnarray*}
Now choose $n$ from the condition
\[
V( R_{n+1}) <\mu( \Omega) \leq V(
R_{n}),
\]
and set $r=R_{n}$. We obtain then%
%
%e3.22 ###
%
\begin{eqnarray} \label{esupG2}
\limfunc{esup}_{B( x,\delta r) }G^{\Omega}1_{\Omega}
&\leq&C%
\frac{F( R) }{V( r) ^{1-\nu}V(
R) ^{\nu}}%
\mu( \Omega) +CF( r) \nonumber\\[-8pt]\\[-8pt]
&\leq&CF( R) \biggl( \frac{V( r) }{V(
R) }%
\biggr) ^{\nu}+CF( r) .\nonumber
\end{eqnarray}
Using again (\ref{Va}), (\ref{Fb}) and $\alpha\nu<\beta$, we obtain
\[
\frac{F( r) }{F( R) }\leq C\biggl( \frac
{r}{R}\biggr)
^{\beta}\leq C\biggl( \frac{r}{R}\biggr) ^{a\nu}\leq C\biggl( \frac
{V(
r) }{V( R) }\biggr) ^{\nu},
\]
which implies that the second term in (\ref{esupG2}) can be absorbed
by the
first one, thus giving
\[
\limfunc{esup}_{B( x,\delta r) }G^{\Omega}1_{\Omega
}\leq
CF( R) \biggl( \frac{V( r) }{V(
R) }\biggr)
^{\nu}\leq CF( R) \biggl( \frac{\mu( \Omega
) }{%
V( R) }\biggr) ^{\nu}.
\]
Since the point $x\in\Omega$ was arbitrary, covering $\Omega$ by a
countable family of balls like $B( x,\delta r) $, we obtain
\[
\limfunc{esup}_{\Omega}G^{\Omega}1_{\Omega}\leq CF( R)
\biggl(
\frac{\mu( \Omega) }{V( R) }\biggr) ^{\nu},
\]
which together with (\ref{laG}) finishes the proof.
\end{pf}

%s3.4 ###
\subsection{Estimates of the exit time}

Our main result in this section is Theorem \ref{TEF} saying that the
condition (\ref{EFF}) implies a certain upper bound for
the tail $%
\mathbb{P}_{x}( \tau_{B}\leq t) $ of the exit time from balls.
The results of this type go back to Barlow \cite{Barlow}, Theorem~3.11. Here
we give a self-contained proof in the present setting, which is based
on the
ideas of \cite{Barlow}. An alternative analytic approach can be found in
\cite{GrigHuUpper}.

For any open set $\Omega\subset M$, set%
%
%e3.23 ###
%
\begin{equation}\label{Ebarsup}
\overline{E}( \Omega) =\sup_{\Omega\setminus\mathcal
{N}_{0}}%
\mathbb{E}_{x}\tau_{\Omega}.
\end{equation}

\begin{lemma}
\label{lTtail}For any open $\Omega\subset M$ such that $\overline
{E}(
\Omega) <\infty$, we have, for all $t>0$ and $x\in\Omega
\setminus
\mathcal{N}_{0}$,%
%
%e3.24 ###
%
\begin{equation}\label{Px<}
\mathbb{P}_{x}( \tau_{\Omega}<t) \leq1-\frac{\mathbb
{E}%
_{x}( \tau_{\Omega}) }{\overline{E}( \Omega
) }+%
\frac{t}{\overline{E}( \Omega) }.
\end{equation}
\end{lemma}
\begin{pf}
Denote $\tau=\tau_{\Omega}$, and observe that%
\[
\tau\leq t+( \tau-t) \mathbf{1}_{\{ \tau\geq
t\}
}=t+( \tau\circ\Theta_{t}) \mathbf{1}_{\{ \tau
\geq
t\} },\vadjust{\goodbreak}
\]
where $\Theta_{t}$ is the time shift of trajectories. Using the Markov
property, we obtain, for any $x\in\Omega\setminus\mathcal{N}_{0}$,%
\[
\mathbb{E}_{x}\tau\leq t+\mathbb{E}_{x}\bigl( ( \tau\circ
\Theta
_{t}) \mathbf{1}_{\{ \tau\geq t\} }\bigr)
=t+\mathbb{E}%
_{x}\bigl( \mathbb{E}_{X_{t}}( \tau) \mathbf{1}_{
\{ \tau
\geq t\} }\bigr),
\]
whence
\[
\mathbb{E}_{x}\tau\leq t+\mathbb{P}_{x}( \tau\geq t)
\sup_{y\in
\Omega\setminus\mathcal{N}_{0}}\mathbb{E}_{y}\tau=t+\mathbb
{P}_{x}(
\tau\geq t) \overline{E}( \Omega)
\]
(see Figure \ref{pic7}),
%
%f3 ###
%
\begin{figure}

\includegraphics{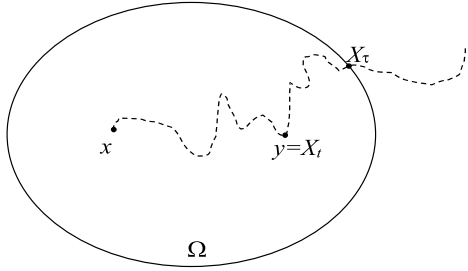}

\caption{Illustration to the proof of
Lemma \protect\ref{lTtail}.}\label{pic7}
\end{figure}
and (\ref{Px<}) follows.
\end{pf}
\begin{lemma}
\label{lE<1-e}Assume that the condition (\ref{EFF}) is satisfied.
Then there are constants $\varepsilon,\sigma>0$ such that, for all $%
x\in M\setminus\mathcal{N}_{0}$, $R>0$, and $\lambda\geq\frac
{\sigma}{%
F( R) }$,%
%
%e3.25 ###
%
\begin{equation}\label{Exe}
\mathbb{E}_{x}\bigl( e^{-\lambda\tau_{B( x,R) }}
\bigr) \leq
1-\varepsilon.
\end{equation}
\end{lemma}
\begin{pf}
Denoting $B=B( x,R) $ and using Lemma \ref{lTtail}, we
have, for
any \mbox{$t>0$},
\begin{eqnarray*}
\mathbb{E}_{x}( e^{-\lambda\tau_{B}}) &=&\mathbb
{E}_{x}\bigl(
e^{-\lambda\tau_{B}}\mathbf{1}_{\{ \tau_{B}<t\}
}\bigr) +%
\mathbb{E}_{x}\bigl( e^{-\lambda\tau_{B}}\mathbf{1}_{\{ \tau
_{B}\geq
t\} }\bigr) \\
&\leq&\mathbb{P}_{x}( \tau_{B}<t) +e^{-\lambda t} \\
&\leq&1-\frac{\mathbb{E}_{x}\tau_{B}}{\overline{E}( B
) }+\frac{t%
}{\overline{E}( B) }+e^{-\lambda t}.
\end{eqnarray*}
The condition (\ref{EFF}) implies that%
\[
\overline{E}( B) =\sup_{z\in B( x,R)
\setminus
\mathcal{N}_{0}}\mathbb{E}_{z}\tau_{B( x,R) }\leq\sup
_{z\in
M\setminus\mathcal{N}_{0}}\mathbb{E}_{z}\tau_{B( z,2R)
}\leq
CF( 2R) ,
\]
whence%
%
%e3.26 ###
%
\begin{equation}\label{EbarF}
\overline{E}( B) \leq C\mathbb{E}_{x}\tau_{B}.
\end{equation}
Using these two estimates of $\overline{E}( B) $, we
obtain%
\[
\mathbb{E}_{x}( e^{-\lambda\tau_{B}}) \leq1-\frac
{1}{C}+\frac{%
Ct}{F( R) }+e^{-\lambda t}.
\]
Setting $\varepsilon=\frac{1}{3C}$ and choosing $t=\frac{\varepsilon
}{C}%
F( R) $, we obtain%
\[
\mathbb{E}_{x}( e^{-\lambda\tau_{B}}) \leq
1-3\varepsilon
+\varepsilon+e^{-\lambda t}.
\]
If also $e^{-\lambda t}\leq\varepsilon$, then we obtain (\ref{Exe}).
Clearly, the former condition will be satisfied provided%
\[
\lambda\geq\frac{\log({1/\varepsilon})}{t}=\frac{({C}/{%
\varepsilon})\log({1/\varepsilon})}{F( R) },
\]
which finishes the proof.
\end{pf}
\begin{lemma}
\label{pf<exp}Assume that the condition (\ref{EFF}) is satisfied.
Then there exists constant $\gamma>0$ such that, for all precompact
balls $%
B( x,R) $ with $x\in M\setminus\mathcal{N}_{0}$ and for
all $%
\lambda>0$,%
%
%e3.27 ###
%
\begin{equation}\label{Ex0}
\mathbb{E}_{x}\bigl( e^{-\lambda\tau_{B( x,R) }}
\bigr) \leq
C\exp\biggl( -\gamma\frac{R}{\mathcal{R}( 1/\lambda)}\biggr),
\end{equation}
where $\mathcal{R}=F^{-1}$.
\end{lemma}
\begin{pf}
Rename the center $x$ of the ball to $z$ so that the letter $x$ will be used
to denote a variable point. Fix some $0<r<R$ to be specified later, and
set $%
n=[ \frac{R}{r}] $. Set also $\tau=\tau_{B(
z,R) }$,
\[
u( x) =\mathbb{E}_{x}( e^{-\lambda\tau})
\]
and%
\[
m_{k}=\sup_{\overline{B}( z,kr) \setminus\mathcal{N}_{0}}u,
\]
where $k=1,2,\ldots,n$. Note that all $m_{k}$ are bounded by $1$.
Choose $%
0<\varepsilon^{\prime}<\varepsilon$ where $\varepsilon$ is the constant
from Lemma \ref{lE<1-e}, and let $x_{k}$ be a point in $\overline
{B}(
z,kr) \setminus\mathcal{N}_{0}$ for which
\[
( 1-\varepsilon^{\prime}) m_{k}\leq u( x_{k}
) \leq
m_{k}.
\]
Fix $k\leq n-1$, observe that
\[
B( x_{k},r) \subset B\bigl( z,( k+1) r
\bigr) \subset
B( z,R)
\]
and consider the following function in $B( x_{k},r)$:%
\[
v_{k}( x) =\mathbb{E}_{x}( e^{-\lambda\tau
_{k}}) ,
\]
where $\tau_{k}=\tau_{B( x_{k},r) }$ (see Figure \ref{pic4}).
Since the ball $B( x_{k},r) $ is precompact, we have
$X_{\tau
_{k}}\in\overline{B}( x_{k},r) $ (while for noncompact balls
the exit point could have been at the cemetery).

%
%f4 ###
%
\begin{figure}

\includegraphics{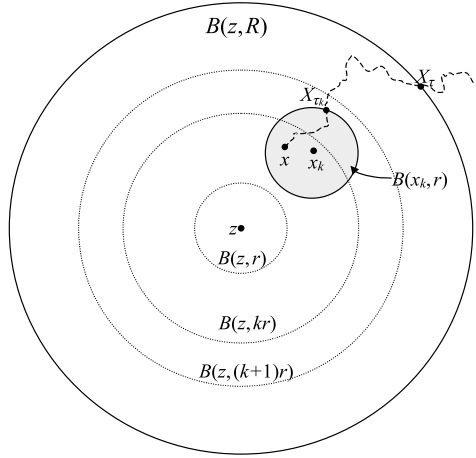}

\caption{Exit times from $B(
x_{k},r) $ and $B( z,R) $.}\label{pic4}
\end{figure}

Let us show that, for all $x\in B( x_{k},r) \setminus
\mathcal{N}%
_{0}$,%
%
%e3.28 ###
%
\begin{equation}\label{u<uu}
u( x) \leq v_{k}( x) \sup_{\overline
{B}(
x_{k},r) \setminus\mathcal{N}_{0}}u.\vadjust{\goodbreak}
\end{equation}
Indeed, we have by the strong Markov property%
\begin{eqnarray*}
u( x) &=&\mathbb{E}_{x}( e^{-\lambda\tau
_{k}}) =%
\mathbb{E}_{x}\bigl( e^{-\lambda\tau_{k}}e^{-\lambda( \tau
-\tau
_{k}) }\bigr) \\
&=&\mathbb{E}_{x}\bigl( e^{-\lambda\tau_{k}}( e^{-\lambda
\tau}\circ
\Theta_{\tau_{k}}) \bigr) \\
&=&\mathbb{E}_{x}( e^{-\lambda\tau_{k}}\mathbb{E}_{X_{\tau
_{k}}}( e^{-\lambda\tau}) ) \\
&=&\mathbb{E}_{x}( e^{-\lambda\tau_{k}}u( X_{\tau
_{k}})
) \\
&\leq&\mathbb{E}_{x}( e^{-\lambda\tau_{k}}) \sup
_{\overline{B}%
( x_{k},r) \setminus\mathcal{N}_{0}}u,
\end{eqnarray*}
which proves (\ref{u<uu}). It follows from (\ref{u<uu}) that%
\[
u( x_{k}) \leq v_{k}( x_{k}) \sup_{\overline
{B}(
z,( k+1) r) \setminus\mathcal{N}_{0}}u=v_{k}(
x_{k}) m_{k+1} ,
\]
whence%
\[
( 1-\varepsilon^{\prime}) m_{k}\leq v_{k}(
x_{k})
m_{k+1}.
\]
By Lemma \ref{lE<1-e}, if
%
%e3.29 ###
%
\begin{equation}\label{lar}
\lambda\geq\frac{\sigma}{F( r) } ,
\end{equation}
then $v_{k}( x_{k}) \leq1-\varepsilon$. Therefore, under
hypothesis (\ref{lar}), we have%
\[
( 1-\varepsilon^{\prime}) m_{k}\leq(
1-\varepsilon
) m_{k+1},
\]
whence it follows by iteration that%
%
%e3.30 ###
%
\begin{equation}\label{n}
u( z) \leq m_{1}\leq\biggl( \frac{1-\varepsilon
}{1-\varepsilon
^{\prime}}\biggr) ^{n-1}m_{n}\leq C\exp\biggl( -c\frac{R}{r}\biggr) ,
\end{equation}
where in the last inequality we have used that $n\geq\frac{R}{r}-1$
and $%
c:=\log\frac{1-\varepsilon^{\prime}}{1-\varepsilon}>0$.\vadjust{\goodbreak}

Condition (\ref{lar}) can be satisfied by choosing
%
%e3.31 ###
%
\begin{equation} \label{n=}
r=\mathcal{R}\biggl( \frac{\sigma}{\lambda}\biggr) .
\end{equation}
This value of $r$ is legitimate only if $r<R$, that is, if%
%
%e3.32 ###
%
\begin{equation} \label{R>}
R>\mathcal{R}\biggl( \frac{\sigma}{\lambda}\biggr) .
\end{equation}
If (\ref{R>}) is not fulfilled, then (\ref{Ex0}) is trivially
satisfied by
choosing the constant $C$ large enough. Assuming that (\ref{R>}) is satisfied
and defining $r$ by~(\ref{n=}) we obtain from (\ref{n}) that%
\[
u( z) \leq C\exp\biggl( -c\frac{R}{\mathcal{R}(
\sigma
/\lambda) }\biggr) ,
\]
whence (\ref{Ex0}) follows.
\end{pf}
\begin{theorem}
\label{TEF}Assume that (\ref{EF}) holds. Then, for any
precompact ball~$B( x,R) $ with $x\in M\setminus\mathcal{N}_{0}$
and for any $t>0$,
%
%e3.33 ###
%
\begin{equation} \label{Psi}
\mathbb{P}_{x}\bigl( \tau_{B( x,R) }\leq t\bigr) \leq
C\exp
( -\Phi( \gamma R,t) ),
\end{equation}
where $\gamma>0$ is the constant from Lemma \ref{pf<exp} and
%
%e3.34 ###
%
\begin{equation} \label{Fidef}
\Phi( R,t) =\sup_{r>0}\biggl\{ \frac{R}{r}-\frac
{t}{F( r) }\biggr\}.
\end{equation}
\end{theorem}

Changing the variable $r$ in (\ref{Fidef}), we obtain the following
equivalent definitions of $\Phi$:
%
%e3.35 ###
%
\begin{equation}\label{Fidef2}
\Phi( R,t) =\sup_{\xi>0}\biggl\{ \frac{R}{\mathcal
{R}( \xi
) }-\frac{t}{\xi}\biggr\} =\sup_{\lambda>0}\biggl\{ \frac
{R}{\mathcal{%
R}( 1/\lambda) }-\lambda t\biggr\} ,
\end{equation}
where $\mathcal{R}=F^{-1}$.
\begin{pf*}{Proof of Theorem \ref{TEF}}
Denoting $B=B( x,R) $ and using Lem-\break ma~\ref{pf<exp}, we obtain
that, for any $\lambda>0$,%
%
%e3.36 ###
%
\begin{eqnarray} \label{3a}
\mathbb{P}_{x}( \tau_{B}\leq t) &=&\mathbb{P}_{x}(
e^{-\lambda\tau_{B}}\geq e^{-\lambda t}) \nonumber\\
&\leq&e^{\lambda t}\mathbb{E}_{x}( e^{-\lambda\tau_{B}})
\\
&\leq&C\exp\biggl( -\gamma\frac{R}{\mathcal{R}( 1/\lambda
) }%
+\lambda t\biggr) .\nonumber
\end{eqnarray}
Taking the supremum in $\lambda$ and using (\ref{Fidef2}), we obtain
(\ref{Psi}).
\end{pf*}
\begin{remark}
\label{RemFi}
It is clear from (\ref{Fidef}) that function $\Phi
(
R,t) $ is increasing in $R$ and decreasing in $t$. Also, we
have, for
any constants $a,b>0$,%
%
%e3.37 ###
%
\begin{equation}\label{Fiab}
\Phi( aR,bt) =ab\Phi\biggl( \frac{R}{b},
\frac{t}{a}\biggr) .\vadjust{\goodbreak}
\end{equation}
In particular, it follows that
%
%e3.38 ###
%
\begin{equation} \label{FiFi}
\Phi( R,t) =t\Phi\biggl( \frac{R}{t},1\biggr) =t\Phi
\biggl(\frac{R}{t}\biggr) ,
\end{equation}
where%
%
%e3.39 ###
%
\begin{equation} \label{Fidef1}
\Phi( s) :=\Phi( s,1) =\sup_{r>0}\biggl\{
\frac{s}{r}-%
\frac{1}{F( r) }\biggr\}.
\end{equation}
Hence, (\ref{Psi}) can be written also in the form%
%
%e3.40 ###
%
\begin{equation} \label{Psi0}
\mathbb{P}_{x}\bigl( \tau_{B( x,R) }\leq t\bigr) \leq
C\exp
\biggl( -t\Phi\biggl( \gamma\frac{R}{t}\biggr) \biggr) .
\end{equation}

Clearly, $\Phi( 0) =0$. Let us show that $0<\Phi(
s)
<\infty$ for all $s>0$. Since%
\[
\lim_{r\rightarrow\infty}\biggl( \frac{s}{r}-\frac{1}{F(
r) }%
\biggr) =0,
\]
we see from (\ref{Fidef1}) that $\Phi( s) \geq0$. It follows
from (\ref{Fb}) and $\beta>1$ that
%
%e3.41 ###
%
\begin{equation} \label{Frr}
\lim_{r\rightarrow0}\frac{r}{F( r) }=\infty
\quad\mbox{and}\quad
\lim_{r\rightarrow+\infty}\frac{r}{F( r) }=0.
\end{equation}
Given $s>0$, choose $r$ so big that $\frac{r}{F( r) }<s$
[such $r$
exists by the second condition in (\ref{Frr})]. Then%
\[
\Phi( s) \geq\frac{s}{r}-\frac{1}{F( r) }>0.
\]
In order to prove that $\Phi( s) <\infty$, it suffices
to show
that%
\[
\lim_{r\rightarrow0}\biggl( \frac{s}{r}-\frac{1}{F( r)
}\biggr)
\leq0.
\]
Indeed, if $r$ is sufficiently small, then by the first condition in
(\ref{Frr}), $\frac{r}{F( r) }>s$ whence $\frac{s}{r}<\frac
{1}{%
F( r) }$.\vspace*{2pt}

Another useful property of function $\Phi( s) $ is the
inequality%
\label{remdoweuseit?}%
%
%e3.42 ###
%
\begin{equation} \label{Fi2s}
\Phi( as) \geq a\Phi( s) \qquad\mbox{for all
}s\geq0
\mbox{ and }a\geq1.
\end{equation}
Indeed, we have for any $r>0$
\[
\frac{as}{r}-\frac{1}{F( r) }\geq a\biggl( \frac
{s}{r}-\frac{1}{%
F( r) }\biggr) ,
\]
whence (\ref{Fi2s}) follows by taking $\sup$ in $r$.
\end{remark}
\begin{example}
If $F( r) $ is differentiable then the supremum in
(\ref{Fidef1}) is attained at the value of $r$ that solves the equation%
\[
\frac{r^{2}F^{\prime}( r) }{F^{2}( r) }=s.\vadjust{\goodbreak}
\]
For example, $F( r) =r^{\beta}$ then we obtain $r=(
\frac{%
\beta}{s}) ^{{1}/({\beta-1})}$ whence $\Phi( s
) =cs^{%
{\beta}/({\beta-1})}$ and%
\[
\Phi( R,t) =c\biggl( \frac{R^{\beta}}{t}\biggr)
^{{1}/({\beta
-1})}.
\]
Consequently, (\ref{Psi}) becomes%
\[
\mathbb{P}_{x}\bigl( \tau_{B( x,R) }\leq t\bigr) \leq
C\exp
\biggl( -c\biggl( \frac{R^{\beta}}{t}\biggr) ^{{1}/({\beta
-1})}\biggr) .
\]
\end{example}
\begin{example}
\label{ExF120}
Consider the following example of function $F$:%
%
%e3.43 ###
%
\begin{equation}\label{F120}
F( r) =\cases{
r^{\beta_{1}},&\quad$r<1$, \cr
r^{\beta_{2}},&\quad$r\geq1$,}
\end{equation}
where $\beta_{1},\beta_{2}>1$. It is easy to see that (\ref{Fb}) is
satisfied with $\beta=\beta_{1}\wedge\beta_{2}$ and $\beta^{\prime
}=\beta_{1}\vee\beta_{2}$. Similarly to the previous example,
one obtains that%
%
%e3.44 ###
%
\begin{equation} \label{Fi120}
\Phi( s) \simeq\cases{
s^{{\beta_{1}}/({\beta_{1}-1})}, &\quad$s>1$, \cr
s^{{\beta_{2}}/({\beta_{2}-1})}, &\quad$s\leq1$,}
\end{equation}
so that (\ref{Psi}) becomes%
\[
\mathbb{P}_{x}\bigl( \tau_{B( x,R) }\leq t\bigr) \leq
C\cases{
\displaystyle \exp\biggl( -c\biggl( \frac{R^{\beta_{1}}}{t}\biggr) ^{
{1}/({\beta_{1}-1})%
}\biggr) , &\quad$t<R$, \vspace*{2pt}\cr
\displaystyle \exp\biggl( -c\biggl( \frac{R^{\beta_{2}}}{t}\biggr) ^{
{1}/({\beta_{2}-1})%
}\biggr) , &\quad$t\geq R$.}
\]
\end{example}
\begin{lemma}
\label{LemtiFi}The function $\Phi( R,t) $ satisfies the
following inequality:%
%
%e3.45 ###
%
\begin{equation} \label{tFi}
\Phi( R,t) \geq c\min\biggl\{ \biggl( \frac{F(
R) }{t}%
\biggr) ^{{1}/({\beta^{\prime}-1})},\biggl( \frac{F(
R) }{t}%
\biggr) ^{{1}/({\beta-1})}\biggr\}
\end{equation}
for all $R,t>0$.
\end{lemma}
\begin{pf}
By (\ref{Fidef}), we have, for any $r>0$,%
\[
\Phi( R,t) \geq\frac{R}{r}-\frac{t}{F( r) }.
\]
We claim that there exists $r>0$ such that%
%
%e3.46 ###
%
\begin{equation} \label{rt}
\frac{t}{F( r) }=\frac{1}{2}\frac{R}{r}.
\end{equation}
Indeed, the function $\frac{F( r) }{r}$ is continuous on
$(
0,+\infty) $, tends to $0$ as $r\rightarrow0$ and to $\infty$
as $%
r\rightarrow\infty$ so that $\frac{F( r) }{r}$ takes all\vspace*{1pt}
positive values, whence the claim follows. With the value of $r$ as in
(\ref{rt}), we have
%
%e3.47 ###
%
\begin{equation} \label{lati}
\Phi( R,t) \geq\frac{t}{F( r) }.\vadjust{\goodbreak}
\end{equation}
If $r\leq R$ then using the left-hand side inequality of (\ref{Fb}), we
obtain%
\[
\frac{R}{r}\geq c\biggl( \frac{F( R) }{F( r
) }\biggr)
^{1/\beta},
\]
which together with (\ref{rt}) yields
\[
F( r) \leq C\biggl( \frac{t^{\beta}}{F( R)
}\biggr) ^{%
{1}/({\beta-1})}.
\]
Substituting into (\ref{lati}), we obtain%
\[
\Phi( R,t) \geq c\biggl( \frac{F( R)
}{t}\biggr) ^{%
{1}/({\beta-1})}.
\]
Similarly, it $r>R$ then using the right-hand side inequality in (\ref{Fb})
we obtain%
\[
\frac{R}{r}\geq c\biggl( \frac{F( R) }{F( r
) }\biggr)
^{1/\beta^{\prime}},
\]
whence it follows that
\[
\Phi( R,t) \geq c\biggl( \frac{F( R)
}{t}\biggr) ^{%
{1}/({\beta^{\prime}-1})}.
\]
\upqed\end{pf}
\begin{corollary}
\label{CEF}Under the hypotheses of Theorem \ref{TEF}, we
have, for
any $x\in M\setminus\mathcal{N}_{0}$, $R>0$, $t>0$,
%
%e3.48 ###
%
\begin{equation} \label{Ptau}
\mathbb{P}_{x}\bigl( \tau_{B( x,R) }\leq t\bigr) \leq
C\exp
\biggl( -c\biggl( \frac{F( R) }{t}\biggr) ^{
{1}/({\beta^{\prime
}-1})}\biggr) .
\end{equation}
\end{corollary}
\begin{pf}
Indeed, if $\frac{F( R) }{t}\geq1$, then (\ref{Ptau}) follows
from Theorem \ref{TEF}, Lem\-ma~\ref{LemtiFi} and (\ref{Fb}). If
$\frac{%
F( R) }{t}<1$, then (\ref{Ptau}) is trivial.
\end{pf}

%s4 ###
\section{Upper bounds of heat kernel}
\label{SecDUE}

The following result will be
used in
the proof of Theorem \ref{TDUE} below.
\begin{proposition}[(\cite{GrigHuUpper}, Lemma 5.5)]
\label{Pdirhk} Let $U$ be an open
subset of $M$, and assume that, for any nonempty open set $\Omega
\subset U$,%
\[
\lambda_{\min}( \Omega) \geq a\mu( \Omega
) ^{-\nu}
\]
for some positive constants $a,\nu$. Then the semigroup $\{P_{t}^{B}\}
$ is
ultracontractive with the following estimate:%
%
%e4.1 ###
%
\begin{equation} \label{Uultra}
\Vert{}P_{t}^{B}f\Vert_{\infty}\leq C( at) ^{-1/\nu
}\Vert
{}f\Vert_{1}
\end{equation}
for any $f\in L^{1}( B) $.
\end{proposition}

The next theorem provides pointwise upper bounds for the heat
kernel.\vadjust{\goodbreak}
\begin{theorem}
\label{TDUE}\label{TUE}If the conditions $\mbox{(\ref{VD})} +
\mbox{(\ref{FK})}
+\mbox{(\ref{EFF})}$ are satisfied and all balls are precompact
then the
heat kernel exists with the domain $M\setminus\mathcal{N}$ for some
properly exceptional set $\mathcal{N}$, and satisfies the upper
bound\label{condUE}%
{\renewcommand{\theequation}{\textit{UE}}
\begin{equation}\label{UEE}
p_{t}( x,y) \leq\frac{C}{V( x,\mathcal{R}(
t))}\exp\biggl( -\frac{1}{2}\Phi( cd( x,y)
,t)\biggr)
\end{equation}}

\vspace*{-8pt}

\noindent
for all $t>0$ and $x,y\in M\setminus\mathcal{N}$, where
$\mathcal{%
R}=F^{-1}$ and $\Phi$ is defined by (\ref{Fidef}).
\end{theorem}
\begin{remark}
As it follows from Theorem \ref{TG=>FK}, the hypotheses $\mbox{(\ref{VD})}
+ \mbox{(\ref{FK})} + \mbox{(\ref{EFF})}$ here can be
replaced by $%
\mbox{(\ref{VD})} +\mbox{(\ref{condH})} + \mbox{(\ref{EFF})}$. Also,
using (\ref{FiFi}), one can write (\ref{UEE}) in the form%
\[
p_{t}( x,y) \leq\frac{C}{V( x,\mathcal{R}(
t)
) }\exp\biggl( -\frac{1}{2}t\Phi\biggl( c\frac{d(
x,y) }{t}%
\biggr) \biggr)
\]
as it was stated in the \hyperref[remmorepictures]{Introduction}.
\end{remark}
\begin{remark}
A version of Theorem \ref{TDUE} was proved by Kigami \cite{KigamiNash}
under additional assumptions that the heat kernel is a priori
continuous and
ultracontractive, and using instead of (\ref{FK}) a local Nash
inequality. In the case $F( r) =r^{\beta}$, another
version of
Theorem \ref{TDUE} was proved in \cite{GrigHuUpper}, where the upper
bound (\ref{UEE}) was understood for \textit{almost all} $x,y$. The proof
below uses a combination of techniques from \cite{GrigHuUpper} and
\cite{KigamiNash}.
\end{remark}
\begin{example}
If function $F( r) $ is given by (\ref{F120}) as in
Example %
\ref{ExF120}, then $\Phi( s) $ is given by (\ref{Fi120}) and (\ref{UEE})
becomes%
\[
p_{t}( x,y) \leq\frac{C}{V( x,\mathcal{R}(
t)
) }\cases{
\displaystyle \exp\biggl( -c\biggl( \frac{r^{\beta_{1}}}{t}\biggr) ^{
{1}/({\beta_{1}-1})%
}\biggr) , &\quad$t<r$, \vspace*{2pt}\cr
\displaystyle \exp\biggl( -c\biggl( \frac{r^{\beta_{2}}}{t}\biggr) ^{
{1}/({\beta_{2}-1})%
}\biggr) , &\quad$t\geq r$,}
\]
where $r=d( x,y) $.
\end{example}
\begin{pf*}{Proof of Theorem \protect\ref{TDUE}}
The hypothesis (\ref{FK}) can be stated as follows: for any
ball $%
B=B( x,r) $ where $x\in M$ and $r>0$, and for any
nonempty open
set $\Omega\subset B$, we have%
%
%e4.2 ###
%
\setcounter{equation}{1}
\begin{equation} \label{laa}
\lambda_{\min}( \Omega) \geq a( B) \mu
(
\Omega) ^{-\nu},
\end{equation}
where%
%
%e4.3 ###
%
\begin{equation}\label{aB}
a( B) =\frac{c}{F( r) }\mu( B)
^{\nu}
\end{equation}
and $\nu,c$ are positive constants. Hence, (\ref{FK})
implies by
Proposition \ref{Pdirhk} that%
%
%e4.4 ###
%
\begin{equation} \label{Bultra}
\Vert{}P_{t}^{B}f\Vert_{L^{1}\rightarrow L^{\infty}}\leq\frac
{C}{(
a( B) t) ^{1/\nu}}.\vadjust{\goodbreak}
\end{equation}
In particular, the semigroup $\{ P_{t}^{B}\} $ is
ultracontractive and $\{ P_{t}\} $ is locally ultracontractive.
By Theorem \ref{TptOm}, there exists a properly exceptional set
$\mathcal{N}%
\subset M$ (containing $\mathcal{N}_{0}$) such that, for any open
subset $%
\Omega\subset M$, the semigroup $\mathcal{P}_{t}^{\Omega}$ possesses the
heat kernel $p_{t}^{\Omega}( x,y) $ with the\vspace*{1pt} domain
$\Omega
\setminus\mathcal{N}$. Fix this set $\mathcal{N}$ for what follows. By
Theorem \ref{TBBCK}, (\ref{aB}) and (\ref{Bultra}) imply the following
estimate:
%
%e4.5 ###
%
\begin{equation} \label{ptB}
p_{t}^{B}( x,y) \leq\frac{C}{( a( B)
t)
^{1/\nu}}=\frac{C}{\mu( B) }\biggl( \frac{F(
r) }{t}%
\biggr) ^{1/\nu}
\end{equation}
for any ball $B$ of radius $r$, and for all $t>0$, $x,y\in
B\setminus\mathcal{N}$.

Our next step is to prove the on-diagonal estimate%
{\renewcommand{\theequation}{\textit{DUE}}
\begin{equation}
p_{t}( x,x) \leq\frac{C}{V( x,\mathcal{R}(
t)
) }
\end{equation}}

\vspace*{-8pt}

\noindent
for all $x\in M\setminus\mathcal{N}$ and $t>0$. To understand the
difficulties, let us first consider a particular case when the volume
function satisfies the following estimate:%
%
%e4.6 ###
%
\setcounter{equation}{5}
\begin{equation} \label{VF}
V( x,R) \simeq F( R) ^{1/\nu}
\end{equation}
for all $x\in M$ and $R>0$, where $\nu$ is the exponent in (\ref{FK})
[e.g., (\ref{VF}) holds, if $V( x,R) \simeq
R^{\alpha}$, $F( R) =R^{\beta}$ and $\nu=\beta/\alpha$].
In this case, the value $F( R) $ in~(\ref{FK}) cancels
our, and we obtain
%
%e4.7 ###
%
\begin{equation} \label{FKc}
\lambda_{\min}( \Omega) \geq c\mu( \Omega
) ^{-\nu}.
\end{equation}
Hence, by Proposition \ref{Pdirhk}, the semigroup $\{
P_{t}\} $
is ultracontractive, and by Theorem~\ref{TBBCK} we obtain the estimate%
%
%e4.8 ###
%
\begin{equation} \label{t1nu}
p_{t}( x,x) \leq Ct^{-1/\nu}
\end{equation}
for all $x\in M\setminus\mathcal{N}$ and $t>0$. Observing that%
\[
V( x,\mathcal{R}( t) ) \simeq F(
\mathcal{R}%
( t) ) ^{1/\nu}=t^{1/\nu},
\]
we see that (\ref{t1nu}) is equivalent to (\ref{DUE}). Although
this argument works only under restriction (\ref{VF}), it has an
advantage that it can be localized as follows. Assuming that (\ref{VF}) is
satisfied for all $R<R_{0}$ with some fixed constant $R_{0}$,~(\ref
{FKc}) is
satisfied for all open sets $\Omega$ with a bounded value of~$\mu
(
\Omega) $, and~(\ref{EFF}) is satisfied for all
balls with
a bounded value of the radius, one can prove in the same way that (\ref
{t1nu}) is true for $t<t_{0}$ for some $t_{0}>0$. The proof below
does not
allow such a localization in the general case.\looseness=-1

In the general case, without the hypotheses (\ref{VF}), the heat
semigroup~$\{ P_{t}\} $ is not necessarily ultracontractive, which requires
other tools for obtaining~(\ref{DUE}). In the case of Riemannian
manifolds, one can obtain~(\ref{DUE}) from~(\ref{FK})
using a certain mean value inequality (see \cite{GrigNotes,CouGgraph})
but this method heavily relies on a specific property of the distance
function that \mbox{$\vert\nabla d\vert\leq1$},\vadjust{\goodbreak} which is not
available in our generality. We will use Kigami's iteration argument that
allows us to obtain (\ref{DUE}) from (\ref{ptB}) using, in
addition, the
hypothesis~(\ref{EFF}). This argument is presented in an abstract
form in~\cite{GrigHuUpper}, Lemma 5.6, that says the following. Assume that
the following two conditions are satisfied:

\begin{longlist}[(2)]
\item[(1)] for any ball $B=B( x,r) $,%
%
%e4.9 ###
%
\begin{equation}\label{Psit}
\limfunc{esup}_{B}p_{t}^{B}\leq\Psi_{t}( x,r),
\end{equation}
where function $\Psi_{t}( x,t) $ satisfies certain
conditions%
\footnote{%
Function $\Psi_{t}( x,r) $ should be monotone decreasing
in $t$
and should satisfy the following doubling condition: if $r\leq
r^{\prime
}\leq2r$ and $t^{\prime}\geq t/2$, then
\[
\Psi_{t^{\prime}}( x,r^{\prime}) \leq K\Psi_{t}(
x,r)
\]
for some constant $K$. This is obviously satisfied for the function
$\Psi$
given by (\ref{Psidef}).};

\item[(2)] for all $x\in M\setminus\mathcal{N}_{0}$, $t>0$, and $r\geq
\varphi
( t) $,%
%
%e4.10 ###
%
\begin{equation}\label{Pxtau}
\mathbb{P}_{x}( \tau_{B}\leq t) \leq\varepsilon,
\end{equation}
where $\varepsilon>0$ is a sufficiently small\footnote{%
More precisely, this means that $\varepsilon\leq\frac{1}{2K}$ where
$K$ is
the constant from the conditions for $\Psi$.} constant and $\varphi$
is a
positive increasing function on $( 0,+\infty) $ such that%
%
%e4.11 ###
%
\begin{equation}\label{intfi}
\int_{0}\varphi( t) \,\frac{dt}{t}<\infty.
\end{equation}
\end{longlist}

Then the heat kernel on $M$ satisfies the estimate%
%
%e4.12 ###
%
\begin{equation} \label{ptxx}
\limfunc{esup}_{B( x,\varphi( t) )
}p_{t}\leq C\Psi
_{t}( x,\varphi( t) ) .
\end{equation}
Obviously, (\ref{ptB}) implies (\ref{Psit}) with the function
%
%e4.13 ###
%
\begin{equation}\label{Psidef}
\Psi_{t}( x,r) =\frac{C}{V( x,r) }\biggl(
\frac{%
F( r) }{t}\biggr) ^{1/\nu}.
\end{equation}
By Corollary \ref{CEF}, (\ref{EFF}) implies (\ref{Ptau}), which
means that (\ref{Pxtau}) is satisfied provided $Ct\leq F(
r) $
for a sufficiently large constant $C$; hence, the function~$\varphi
(
t) $ can be chosen as follows:
\[
\varphi( t) =\mathcal{R}( Ct) ,
\]
which clearly satisfies (\ref{intfi}) [indeed, by (\ref{Rb}), we have
$%
\mathcal{R}( t) \leq Ct^{1/\beta^{\prime}}$ for all $0<t<1$,
whence (\ref{intfi}) follows].

By (\ref{ptxx}), we obtain%
\[
\limfunc{esup}_{B( x,\varphi( t) )
}p_{t}\leq C\Psi
_{t}( x,\mathcal{R}( Ct) ) \leq\frac
{C}{V( x,%
\mathcal{R}( Ct) ) }\leq\frac{C}{V(
x,\mathcal{R}%
( t) ) },
\]
where we have also used (\ref{Rb}) and (\ref{Va}). By Theorem \ref
{TptOm}(d), we can replace here $\limfunc{esup}p_{t}$ by
$\sup p_{t}$
outside $\mathcal{N}$, whence (\ref{DUE}) follows.\vadjust{\goodbreak}

Now we prove the full upper estimate (\ref{UEE}). Fix two disjoint
open subsets~$U,V$ of $M$ and use the following inequality proved in
\cite{BarGrigKum}, Lemma 2.1: for all functions $f,g\in\mathcal{B}_{+}(
M) $,%
%
%e4.14 ###
%
\begin{eqnarray} \label{Ptfgb}
(\mathcal{P}_{t}f,g)&\leq&\bigl(\mathbb{E}_{\cdot}\bigl(\mathbf{1}_{\{\tau
_{U}\leq
t/2\}}\mathcal{P}_{t-\tau_{U}}f(X_{\tau_{U}})\bigr),g\bigr)\nonumber\\[-9pt]\\[-9pt]
&&{} +\bigl(\mathbb
{E}_{\cdot}\bigl(%
\mathbf{1}_{\{\tau_{V}\leq t/2\}}\mathcal{P}_{t-\tau_{V}}g(X_{\tau
_{V}})\bigr),f\bigr)\nonumber
\end{eqnarray}
(see Figure \ref{pic5}).

%
%f5 ###
%
\begin{figure}

\includegraphics{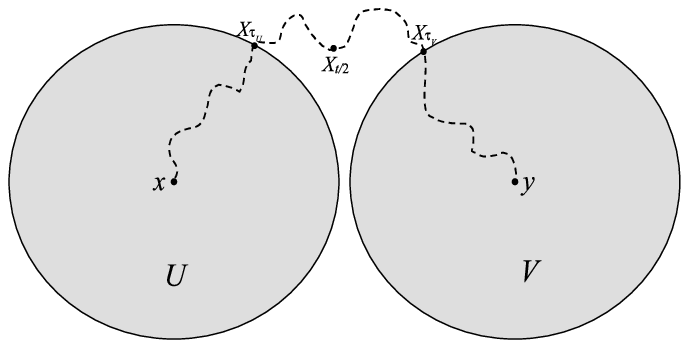}

\caption{Illustration to (\protect\ref{Ptfgb}).}\label{pic5}
\vspace*{-3pt}
\end{figure}

Assume in addition that $f\in\mathcal{B}L^{1}( V) $ and
$g\in
\mathcal{B}L^{1}( U) $. Then, under the condition $\tau
_{U}\leq
t/2$, we have%
\[
\mathcal{P}_{t-\tau_{U}}f(X_{\tau_{U}})=\int_{V\setminus\mathcal
{N}%
}p_{t-\tau_{U}}( X_{\tau_{U}},y) f( y) \,d\mu
(
y) \leq S\Vert{}f\Vert_{1}
\]
almost surely, where%
%
%e4.15 ###
%
\begin{equation} \label{S}
S:=\sup_{t/2\leq s\leq t}\mathop{\sup_{u\in\overline{U}\setminus
\mathcal{N}}}_{v\in\overline{V}\setminus\mathcal{N}}p_{s}(
u,v) .
\end{equation}
Here we have used that $X_{\tau_{U}}\in\overline{U}\setminus
\mathcal{N}$
almost surely, which is due to the fact that $\{ X_{t}\} $
is a
diffusion and the set $\mathcal{N}$ is properly exceptional. It
follows~that%
\[
\bigl(\mathbb{E}_{\cdot}\bigl(\mathbf{1}_{\{\tau\leq t/2\}}\mathcal
{P}_{t-\tau
}f(X_{\tau})\bigr),g\bigr)\leq S\Vert{}f\Vert_{1}\int_{U}\mathbb{P}_{x}
\biggl( \tau
_{U}\leq\frac{t}{2}\biggr) g( x) \,d\mu( x
) .
\]
Estimating similarly the second term in (\ref{Ptfgb}), we obtain from
(\ref{Ptfgb})%
\begin{eqnarray*}
&&\int_{U}\int_{V}p_{t}( x,y) f( y) g(
x)
\,d\mu( x) \,d\mu( y) \\[-2pt]
&&\qquad\leq S\Vert
f\Vert
_{1}\int_{U}\mathbb{P}_{x}\biggl( \tau_{U}\leq\frac{t}{2}\biggr)
g(
x) \,d\mu( x)\\[-2pt]
&&\qquad\quad{} +S\Vert g\Vert_{1}\int_{V}\mathbb{P}_{y}\biggl( \tau
_{V}\leq
\frac{t}{2}\biggr) f( y) \,d\mu( y) .
\end{eqnarray*}
By \cite{GrigHuUpper}, Lemma 3.4, we conclude that, for $\mu
$-a.a. $x\in V$ and $y\in U$,
%
%e4.16 ###
%
\begin{equation}\label{p<pB2}
p_{t}( x,y) \leq S\mathbb{P}_{x}\biggl( \tau_{V}\leq
\frac{t}{2}%
\biggr) +S\mathbb{P}_{y}\biggl( \tau_{U}\leq\frac{t}{2}\biggr)
.\vadjust{\goodbreak}
\end{equation}
A slightly different inequality (i.e., also enough for our purposes) was
proved in~\cite{GrigHuComp}. For the case of heat kernels on Riemannian
manifolds, (\ref{p<pB2}) was proved in \cite{GrigSal}, Lemma 3.3.

By Theorem \ref{TEF}, we have%
%
%e4.17 ###
%
\begin{equation} \label{Psi1}
\mathbb{P}_{x}\bigl( \tau_{B( x,R) }\leq t\bigr) \leq
C\exp
( -\Phi( \gamma R,t) )
\end{equation}
for all $x\in M\setminus\mathcal{N}$ and $t,R>0$. Let
\[
V_{R}=\{ x\in V\dvtx d( x,V^{c}) >R\} .
\]
Then, for any $x\in V_{R}\setminus\mathcal{N}$, we obtain by (\ref
{Psi1})%
\[
\mathbb{P}_{x}\biggl( \tau_{V}\leq\frac{t}{2}\biggr) \leq\mathbb
{P}%
_{x}\bigl( \tau_{B( x,R) }\leq t\bigr) \leq C\exp
( -\Phi
( \gamma R,t) ) .
\]
Using a similar estimate for $y\in U_{R}$, we obtain from (\ref{p<pB2})
that, for $\mu$-a.a. $x\in V_{R}$ and $y\in U_{R}$,%
%
%e4.18 ###
%
\begin{equation} \label{ptro}
p_{t}( x,y) \leq CS\exp( -\Phi( \gamma
R,t)
) .
\end{equation}
Since the right-hand side here is a constant in $x,y$, we conclude by
Theorem~\ref{TptOm}(d) that (\ref{ptro}) holds for
all $x\in
V_{R}\setminus\mathcal{N}$ and $y\in U_{R}\setminus\mathcal{N}$.

Now fix two distinct points $x,y\in M\setminus\mathcal{N}$, set
%
%e4.19 ###
%
\begin{equation}\label{14}
R=\tfrac{1}{4}d( x,y)
\end{equation}
and observe that the balls $V=B( x,2R) $ and $U=B(
y,2R) $ are disjoint. Since $x\in V_{R}$ and $y\in U_{R}$, we conclude
that (\ref{ptro}) is satisfied for these points $x,y$ with the above value
of $R$. Let us estimate the quantity $S$ defined by (\ref{S}). Using the
semigroup property and (\ref{DUE}), we obtain, for all
$u,v\in
M\setminus\mathcal{N}$,%
\[
p_{s}( u,v) \leq\sqrt{p_{s}( u,u)
p_{s}(
v,v) }\leq\frac{C}{\sqrt{V( u,\mathcal{R}(
s) )
V( v,\mathcal{R}( s) ) }}.
\]
Observe that by (\ref{Va}), for all $z\in M$,%
%
%e4.20 ###
%
\begin{equation}\label{VRR}
\frac{V( z,\mathcal{R}( s) ) }{V(
x,\mathcal{R}%
( s) ) }\leq C\biggl( 1+\frac{d( x,z)
}{\mathcal{R%
}( s) }\biggr) ^{\alpha}.
\end{equation}
Applying this for $z=u\in\overline{U}$ and $z=v\in\overline{V}$ so
that $%
d( x,u) \leq2R$ and $d( x,v) \leq6R$, and
substituting to the above estimate of $p_{s}( u,v) $ we
obtain%
\[
p_{s}( u,v) \leq\frac{C}{V( x,\mathcal{R}(
s)
) }\biggl( 1+\frac{R}{\mathcal{R}( s) }\biggr)
^{\alpha}.
\]
Using that $s\in[ t/2,t] $ as well as (\ref{Rb}) and
(\ref{Va}),
we obtain
\[
S\leq\frac{C}{V( x,\mathcal{R}( t) )
}\biggl( 1+\frac{R%
}{\mathcal{R}( t) }\biggr) ^{\alpha},\vadjust{\goodbreak}
\]
which together with (\ref{ptro}) yields
%
%e4.21 ###
%
\begin{equation} \label{ptS}
p_{t}( x,y) \leq\frac{C}{V( x,\mathcal{R}(
t)
) }\biggl( 1+\frac{R}{\mathcal{R}( t) }\biggr)
^{\alpha
}\exp( -\Phi( \gamma R,t) ) .
\end{equation}
On the other hand, we have by (\ref{Fidef2})
\[
\Phi( R,t) =\sup_{\xi>0}\biggl\{ \frac{R}{\mathcal
{R}( \xi
) }-\frac{t}{\xi}\biggr\} \geq\frac{R}{\mathcal{R}(
t) }%
-1,
\]
where we have chosen $\xi=t$. Using the elementary estimate%
\[
1+z\leq\frac{1}{a}\exp( az) ,\qquad z>0, 0<a\leq1,
\]
and its consequence%
\[
2+z\leq\frac{2}{a}\exp( az) ,
\]
we obtain%
\[
1+\frac{R}{\mathcal{R}( t) }\leq2+\Phi(
R,t) \leq
\frac{2}{a}\exp\biggl( \frac{a}{\gamma}\Phi( \gamma
R,t) \biggr),
\]
whence%
%
%e4.22 ###
%
\begin{equation}\label{1+}
\biggl( 1+\frac{R}{\mathcal{R}( t) }\biggr) ^{\alpha
}\leq\biggl(
\frac{2}{a}\biggr) ^{\alpha}\exp\biggl( \frac{\alpha a}{\gamma
}\Phi(
\gamma R,t) \biggr) .
\end{equation}
Choosing $a$ small enough and substituting this estimate to (\ref
{ptS}), we obtain~(\ref{UEE}).
\end{pf*}
\begin{remark}
\label{Rembdd}It is desirable to have a localized version of
Theorem \ref{TDUE} when the hypotheses are assumed for balls of bounded
radii and the
conclusions are proved for a bounded range of time. As was already mentioned
in the proof, Kigami's argument requires the ultracontractivity of $%
P_{t}^{B} $ for \textit{all} balls, and (\ref{EFF}) should
also be
satisfied for all balls, because, loosely speaking, one deals with the
estimate of $p_{t}^{B_{k+1}}-p_{t}^{B_{k}}$ for an exhausting sequence of
balls $\{ B_{k}\} $ (see \cite{GrigHuUpper} or \cite{KigamiNash}).
As we will see in Section \ref{SecRVD}, the hypotheses $\mbox{(\ref{VD})}
+\mbox{(\ref{condH})} + \mbox{(\ref{EFF})}$ for all balls imply that
the space $%
( M,d) $ is unbounded. Note that all other arguments used
in this
paper do admit localization.
\end{remark}

%s5 ###
\section{Lower bounds of heat kernel}
\label{SecLow}

%s5.1 ###
\subsection{Oscillation inequalities}

The Harnack inequality (\ref{condH}) is a standing assumption in this
subsection. The main result is Proposition \ref{Posc} that is heavily based
on the oscillation inequality of Lemma \ref{Lemosc}. The latter is
considered as a standard consequence of (\ref{condH}), but we still
provide a full proof to emphasize the use of the precompactness of balls.
\begin{lemma}
\label{LemH}Let $B$ be a precompact ball of radius $r$ in $M$. If
$u\in
\mathcal{F}$ is harmonic in $B$, and if $u\geq a$ in $B$ for some real
constant $a$, then%
%
%e5.1 ###
%
\begin{equation}\label{H-a}
\limfunc{esup}_{\delta B}( u-a) \leq C\limfunc
{einf}_{\delta
B}( u-a) ,
\end{equation}
where $C$ and $\delta$ are the same constants as in (\ref{condH}).
\end{lemma}
\begin{pf}
Let $\psi$ be a cutoff function of $B$ in $M$, that is, $\psi\in
\mathcal{F%
}\cap C_{0}( M) $ and $\psi\equiv1$ in an open
neighborhood of $%
\overline{B}$. The function $v=u-a\psi$ belongs to $\mathcal{F}$ and is
equal to $u-a$ in $B$. Let us show that $v$ is harmonic in $B$. Indeed, for
any $\varphi\in\mathcal{F}( B) $, we have%
%
%e5.2 ###
%
\begin{equation}\label{EEE}
\mathcal{E}( v,\varphi) =\mathcal{E}( u-a\psi
,\varphi
) =\mathcal{E}( u,\varphi) -a\mathcal{E}(
\psi
,\varphi) =0,
\end{equation}
because $\mathcal{E}( u,\varphi) =0$ by the harmonicity
of $u$
in $B$, and $\mathcal{E}( \psi,\varphi) =0$ by the strong
locality as $\psi\equiv1$ in a neighborhood of $\limfunc
{supp}\varphi$.
Applying (\ref{condH}) to $v$, we obtain (\ref{H-a}).
\end{pf}

For any function $f$ on any set $S\subset M$, define the essential
oscillation of $f$ in $S$ by%
\[
\limfunc{eosc}_{S}f=\limfunc{esup}_{S}f-\limfunc{einf}_{S}f.
\]
The following statement is well known for functions in $\mathbb{R}^{n}$
(see, e.g., \cite{MoserEl} and~\cite{Salbook}, Lemma 2.3.2).
\begin{lemma}
\label{Lemosc}There exists $\theta>0$ such that, for any precompact
ball $%
B( x,\allowbreak r) \subset M$, for any nonnegative harmonic function
$u$ in
$B( x,r) $, and any $0<\rho\leq r$,%
%
%e5.3 ###
%
\begin{equation} \label{eosc}
\limfunc{eosc}_{B( x,\rho) }u\leq2\biggl( \frac{\rho
}{r}\biggr)
^{\theta}\limfunc{eosc}_{B( x,r) }u,
\end{equation}
where $\theta$ is a positive constant that depends on the constants in
(\ref{condH}).
\end{lemma}
\begin{pf}
Fix $x\in M$, and write for simplicity $B_{r}=B( x,r) $. Consider
first the case when $\rho=\delta r$ [where $\delta$ is a parameter
from (\ref{condH})] and set%
\[
a=\limfunc{esup}_{B_{r}}u,\qquad b=\limfunc{einf}_{B_{r}}u
\]
and%
\[
a^{\prime}=\limfunc{esup}_{B_{\rho}}u,\qquad b^{\prime}=\limfunc{einf}
_{B_{\rho}}u.
\]
Clearly, $b\leq b^{\prime}\leq a^{\prime}\leq a$. By Lemma \ref{LemH}, we
have%
%
%e5.4 ###
%
\begin{equation}\label{esiu}
\limfunc{esup}_{B_{\rho}}( u-b) \leq C\limfunc
{einf}_{B_{\rho
}}( u-b),
\end{equation}
that is,%
\[
a^{\prime}-b\leq C( b^{\prime}-b) .\vadjust{\goodbreak}
\]
Similarly, applying Lemma \ref{LemH} to function $-u$, we obtain%
\[
\limfunc{esup}_{B_{\rho}}( a-u) \leq C\limfunc
{einf}_{B_{\rho
}}( a-u) ,
\]
whence%
\[
a-b^{\prime}\leq C( a-a^{\prime}) .
\]
Adding up the two inequalities yields%
\[
( a-b) +( a^{\prime}-b^{\prime}) \leq
C(
a-b) -C( a^{\prime}-b^{\prime}),
\]
whence%
\[
a^{\prime}-b^{\prime}\leq\frac{C-1}{C+1}( a-b) ,
\]
that is,%
%
%e5.5 ###
%
\begin{equation}\label{der}
\limfunc{eosc}_{B_{\delta r}}u\leq\gamma\limfunc{eosc}_{B_{r}}u,
\end{equation}
where $\gamma:=\frac{C-1}{C+1}<1$. For an arbitrary $\rho\leq r$,
find a
nonnegative integer $k$ such that
\[
\delta^{k+1}r<\rho\leq\delta^{k}r.
\]
Iterating (\ref{der}), we obtain%
\begin{eqnarray*}
\limfunc{eosc}_{B_{\rho}}u&\leq&\limfunc{eosc}_{B_{\delta
^{k}r}}u\leq
\gamma^{k}\limfunc{eosc}_{B_{r}}u\leq\gamma^{{\log
({r/\rho})}{%
\log({1}/{\delta})}-1}\limfunc{eosc}_{B_{r}}u\\
&=&\frac{1}{\gamma
}\biggl(
\frac{\rho}{r}\biggr) ^{({\log({1/\gamma})})/({\log
({1/\delta})%
})}\limfunc{eosc}_{B_{r}}u.
\end{eqnarray*}
Note that the constant $C$ in (\ref{esiu}) can be assumed to be big enough,
say $C>3$. Then $\gamma>1/2$ and (\ref{eosc}) follows from the previous
line with $\theta={({\log({1/\gamma})})/({\log
({1/\delta})%
})}$.
\end{pf}
\begin{proposition}
\label{Posc}Let $\Omega$ be an open subset of $M$ such that
$\widetilde{E}%
( \Omega) <\infty$. Fix a function $f\in\mathcal
{B}_{b}(
\Omega), $ and set $u=G^{\Omega}f$. Then, for any precompact
ball $%
B( x,r) \subset\Omega$ and all $\rho\in(0,r]$,%
%
%e5.6 ###
%
\begin{equation} \label{oscu}
\limfunc{eosc}_{B( x,\rho) }u\leq2\widetilde{E}(
x,r) \limfunc{esup}_{B( x,r) }\vert f
\vert
+4\biggl( \frac{\rho}{r}\biggr) ^{\theta}\limfunc{esup}_{B(
x,r)
}\vert u\vert,
\end{equation}
where $\theta$ is the same constant as in Lemma \ref{Lemosc}
and%
\[
\widetilde{E}( x,r) :=\widetilde{E}( B(
x,r)
) .
\]
\end{proposition}
\begin{pf}
Write for simplicity $B_{r}:=B( x,r) $. Let us first prove that
if $f\geq0$, then
%
%e5.7 ###
%
\begin{equation} \label{osc+}
\limfunc{eosc}_{B_{\rho}}u\leq\widetilde{E}( x,r)
\limfunc{esup}%
_{B_{r}}f+2\biggl( \frac{\rho}{r}\biggr) ^{\theta}\limfunc
{esup}_{B_{r}}u.\vadjust{\goodbreak}
\end{equation}
By Lemma \ref{LG1-1}, we have for the function $v=G^{B_{r}}f$ that%
\[
\limfunc{eosc}_{B_{\rho}}v\leq\limfunc{esup}_{B_{r}}v\leq
\widetilde{E}%
( x,r) \limfunc{esup}_{B_{r}}f.
\]
The function $w:=u-v$ is harmonic in $B_{r}$ by Lemma \ref{LemG-G} and
nonnegative by Theorem \ref{TptOm}(b). By Lemma \ref{Lemosc}%
, we obtain%
\[
\limfunc{eosc}_{B_{\rho}}w\leq2\biggl( \frac{\rho}{r}\biggr)
^{\theta}%
\limfunc{eosc}_{B_{r}}w\leq2\biggl( \frac{\rho}{r}\biggr) ^{\theta
}%
\limfunc{esup}_{B_{r}}u.
\]
Since $u=v+w$, (\ref{osc+}) follows by adding up the two previous lines.

For a signed function $f$, write $f=f_{+}-f_{-}$ and set%
\[
\overline{u}:=G^{\Omega}f_{+} \quad\mbox{and}\quad\underline
{u}:=G^{\Omega
}f_{-}.
\]
Then $\overline{u}$ and $\underline{u}$ are nonnegative and
$u=\overline{u}-%
\underline{u}$, whence it follows that%
\[
\limfunc{eosc}u=\limfunc{eosc}( \overline{u}-\underline
{u}) \leq
\limfunc{eosc}\overline{u}+\limfunc{eosc}\underline{u}.
\]
Applying (\ref{osc+}) separately to $\overline{u}$ and $\underline
{u}$ and
adding up the inequalities, we obtain~(\ref{oscu}).
\end{pf}

%s5.2 ###
\subsection{Time derivative}

In this section, we assume only the basic hypotheses. If $f\in
L^{2}(
M) $, then, for any $t>0$, the function $u_{t}=P_{t}f$ is in
$\func{dom}%
( \mathcal{L}) $ and satisfies the heat equation%
%
%e5.8 ###
%
\begin{equation} \label{tL}
\partial_{t}u_{t}+\mathcal{L}u=0,
\end{equation}
where $\partial_{t}u_{t}$ is the strong derivative in $L^{2}(
M)
$ of the mapping $t\mapsto u_{t}$; cf. \cite{GrigHu} and~\cite
{Grigbook}, Section 4.3. The argument in the next lemma is well known
in the context of
the semigroup theory (see \cite{Davbooko,Davbook}), and we reproduce
it here for the sake of completeness.
\begin{lemma}
For any $f\in L^{2}( M) $ and all $0\leq s<t$, we have%
%
%e5.9 ###
%
\begin{equation}\label{dudt}
\Vert\partial_{t}u_{t}\Vert_{2}\leq\frac{1}{t-s}\Vert u_{s}\Vert_{2},
\end{equation}
where $u_{t}=\mathcal{P}_{t}f$.
\end{lemma}
\begin{pf}
Let $\{ E_{\lambda}\} _{\lambda\geq0}$ be spectral resolution
in $L^{2}( M) $ of the operator~$\mathcal{L}$. Then we
have%
\begin{eqnarray*}
u_{t} &=&e^{-t\mathcal{L}}f=\int_{0}^{\infty}e^{-t\lambda
}\,dE_{\lambda}f,\\
\partial_{t}u_{t} &=&-\mathcal{L}e^{-t\mathcal{L}}f=\int
_{0}^{\infty
}( -\lambda) e^{-t\lambda}\,dE_{\lambda}f,
\end{eqnarray*}
whence%
\begin{eqnarray*}
\Vert u_{t}\Vert_{2}^{2} &=&\int_{0}^{\infty}e^{-2t\lambda}\,d\Vert
E_{\lambda}f\Vert^{2}, \\
\Vert\partial_{t}u_{t}\Vert_{2}^{2} &=&\int_{0}^{\infty}\lambda
^{2}e^{-2t\lambda}\,d\Vert E_{\lambda}f\Vert^{2}.
\end{eqnarray*}
Since
\[
\lambda^{2}e^{-2t\lambda}=\bigl( \lambda e^{-( t-s)
\lambda
}\bigr) ^{2}e^{-2s\lambda}\leq\frac{1}{( t-s) ^{2}}%
e^{-2s\lambda},
\]
we obtain%
\[
\Vert\partial_{t}u_{t}\Vert_{2}^{2}\leq\frac{1}{( t-s
) ^{2}}%
\int_{0}^{\infty}e^{-2s\lambda}\,d\Vert E_{\lambda}f\Vert^{2}=\frac
{1}{%
( t-s) ^{2}}\Vert u_{s}\Vert_{2}^{2},
\]
which was to be proved.
\end{pf}

In the rest of this section, assume that $p_{t}( x,y) $ is the
heat kernel with domain~$D$.
\begin{corollary}
\label{CdtptL2}For any $t>0$ and $y\in D$, the function $t\mapsto
p_{t}( \cdot,y) $ is strongly differentiable in
$L^{2}(
M) $ and, for all $0<s<t$,%
\[
\Vert{}\partial_{t}p_{t}( \cdot,y) \Vert_{2}\leq\frac
{1}{t-s}%
\sqrt{p_{2s}( y,y) }.
\]
\end{corollary}
\begin{pf}
Setting $f=p_{\varepsilon}( \cdot,y) $ for some
$\varepsilon>0$
and using (\ref{ptPt}), we obtain that the function%
\[
u_{t}=\mathcal{P}_{t}f=p_{t+\varepsilon}( \cdot,y)
\]
satisfies (\ref{dudt}), that is,%
\[
\Vert{}\partial_{t}p_{t+\varepsilon}( \cdot,y) \Vert
_{2}\leq
\frac{1}{t-s}\Vert p_{s+\varepsilon}( \cdot,y) \Vert
_{2}=\frac{%
1}{t-s}\sqrt{p_{2( s+\varepsilon) }( y,y) }.
\]
Renaming $t+\varepsilon$ by $t$ and $s+\varepsilon$ by $s$, we finish the
proof.
\end{pf}

Set%
\[
p_{t}^{\prime}( x,y) \equiv\partial_{t}p_{t}(
x,y) ,
\]
where the strong derivative $\partial_{t}$ is taken with respect to the
first variable $x$. Hence, for any $y\in D$ and $t>0$, $p_{t}^{\prime
}( x,y) $ is an $L^{2}$-function of $x$.
\begin{lemma}
\label{LemqtPt}For all $0<s<t$ and $y\in D$, we have%
%
%e5.10 ###
%
\begin{equation} \label{p'}
p_{t}^{\prime}( \cdot,y) =\mathcal
{P}_{s}p_{t-s}^{\prime
}( \cdot,y) .
\end{equation}
\end{lemma}
\begin{pf}
Indeed, we have by (\ref{ptPt})%
\[
p_{t}( \cdot,y) =\mathcal{P}_{s}p_{t-s}( \cdot
,y) .\vadjust{\goodbreak}
\]
Since $\mathcal{P}_{s}$ is a bounded operator in $L^{2}$, it commutes with
the operator $\partial_{t}$ of strong derivation. Applying the latter to
the both sides of the above identity, we obtain~(\ref{p'}).
\end{pf}
\begin{corollary}
\label{Cdtpt}For all $t>0$, $y\in D$ and $\mu$-a.a.
$x\in D$,
%
%e5.11 ###
%
\begin{equation} \label{dtpt}
\vert\partial_{t}p_{t}( x,y) \vert\leq
\frac{2}{t}%
\sqrt{p_{t/2}( x,x) p_{t/2}( y,y) }.
\end{equation}
\end{corollary}
\begin{pf}
By Lemma \ref{LemqtPt}, we have, for all $y\in D$ and $\mu$-a.a.
$x\in D$,%
\[
p_{t}^{\prime}( x,y) =( p_{s}( x,\cdot)
,p_{t-s}^{\prime}( \cdot,y) ) ,
\]
whence by Corollary \ref{CdtptL2},
\[
\vert p_{t}^{\prime}( x,y) \vert\leq\Vert
{}p_{s}( x,\cdot) \Vert_{2}\Vert p_{t-s}^{\prime}
( \cdot
,y) \Vert_{2}\leq\frac{1}{t-s-r}\sqrt{p_{2s}(
x,x)
p_{2r}( y,y) }
\]
for any $0<r<t-s$. Choosing $s=r=t/4$, we finish the proof of (\ref{dtpt}).
\end{pf}
\begin{remark}
It follows easily from the identity%
\[
p_{t}( x,y) =( p_{s}( x,\cdot)
,p_{t-s}(
\cdot,y) ) ,
\]
that the function $t\mapsto p_{t}( x,y) $ is
differentiable for
all fixed $x,y\in D$ and%
\[
\frac{\partial}{\partial t}p_{t}( x,y) =(
p_{s}(
x,\cdot) ,\partial_{t}p_{t-s}( \cdot,y) )
=(
p_{s}( x,\cdot) ,q_{t-s}( \cdot,y) ) .
\]
Arguing as in the previous proof, we obtain%
\[
\biggl\vert\frac{\partial}{\partial t}p_{t}( x,y)
\biggr\vert
\leq\frac{2}{t}\sqrt{p_{t/2}( x,x) p_{t/2}(
y,y) }
\]
for all $t>0$ and $x,y\in D$. However, for applications we need estimate
(\ref{dtpt}) for the strong derivative $\partial_{t}p_{t}$ rather
than for
the partial derivative $\frac{\partial}{\partial t}p_{t}(
x,y) $.
\end{remark}
\begin{lemma}
\label{LemGu}If $\Omega$ is an open subset of $M$ and if $\widetilde
{E}%
( \Omega) <\infty$, then, for all $t>0$ and $z\in M$, the
function $u_{t}:=p_{t}^{\Omega}( \cdot,z) $ satisfies in
$%
( 0,+\infty) \times\Omega$ the equation%
\[
G^{\Omega}( \partial_{t}u_{t}) +u_{t}=0.
\]
\end{lemma}
\begin{pf}
By Lemma \ref{LG1-1}, the Green operator $G^{\Omega}$ is a bounded operator
in $L^{2}( \Omega) $, and $G^{\Omega}$ is the inverse operator
to $\mathcal{L}^{\Omega}$. Since the function $u_{t}$ satisfies the
equation $\partial_{t}u_{t}+\mathcal{L}^{\Omega}u_{t}=0$, applying $%
G^{\Omega}$ proves the claim.
\end{pf}

%s5.3 ###
\subsection{\texorpdfstring{The H\"{o}lder continuity}{The Holder continuity}}

In this subsection we use the hypotheses $\mbox{(\ref{VD})} +\mbox{(\ref{condH})}
+({E}_{F}\mbox{$\leq$}) $. As it follows from Theorem \ref{TG=>FK} and
Proposition \ref{Pdirhk}, under these hypotheses the heat semigroup
$\{
P_{t}\} $ is locally ultracontractive. Hence, by Theorem~\ref{TptOm},
for any open set $\Omega\subset M$, the heat kernel $p_{t}^{\Omega}$
exists with the domain $\Omega\setminus\mathcal{N}$ where $\mathcal
{N}%
\subset M$ is a fixed properly exceptional set; cf. the beginning of the
proof of Theorem \ref{TDUE}.\vadjust{\goodbreak}
\begin{lemma}
\label{LemoscpB}Let the hypotheses $\mbox{(\ref{VD})} +\mbox{(\ref{condH})}
+({E}_{F}\mbox{$\leq$}) $ be satisfied, and let $\Omega$ be an
open subset
of $M$. Fix $t>0$, $0<\rho\leq\mathcal{R}( t), $ and set
%
%e5.12 ###
%
\begin{equation} \label{r=}
r=( \mathcal{R}( t) ^{\beta}\rho^{\theta}
) ^{{1%
}/({\beta+\theta})},
\end{equation}
where $\beta$ is the exponent from (\ref{Fb}), and $\theta$
is the
constant from Lemma \ref{Lemosc}. Fix also a point $x\in
\Omega
\setminus\mathcal{N}$ and assume that the ball $B( x,r)
$ is
precompact, and its closure is contained in $\Omega$. Then%
%
%e5.13 ###
%
\begin{equation} \label{oscpB}
\limfunc{osc}_{y\in B( x,\rho) \setminus\mathcal{N}%
}p_{t}^{\Omega}( x,y) \leq C\biggl( \frac{\rho
}{\mathcal{R}%
( t) }\biggr) ^{\Theta}\sup_{y\in B( x,r)
\setminus
\mathcal{N}}p_{t/2}^{\Omega}( y,y) ,
\end{equation}
where $\Theta=\frac{\beta\theta}{\beta+\theta}$ and $C$ depends
on the
constants in $({E}_{F}\mbox{$\leq$}) $ and (\ref{Fb}).
\end{lemma}
\begin{pf}
By construction in the proof of Theorem \ref{TptOm}, the heat kernel~$%
p_{t}^{\Omega}$ is obtained as a monotone increasing limit of $p_{t}^{U_{n}}$
as $n\rightarrow\infty$ where $\{ U_{n}\} $ is an
exhaustion of
$\Omega$ by sets $U_{n}$ that are finite union of balls from a countable
base and the convergence is pointwise in $\Omega\setminus\mathcal{N}$.
Suppose for a~moment that we have proved (\ref{oscpB}) for $U_{n}$ instead
of $\Omega$, that is,%
%
%e5.14 ###
%
\begin{eqnarray} \label{ptUn}
&&\sup_{y\in B( x,\rho) \setminus\mathcal
{N}}p_{t}^{U_{n}}(
x,y)\nonumber\\[-8pt]\\[-8pt]
&&\qquad \leq\inf_{y\in B( x,\rho) \setminus
\mathcal{N}%
}p_{t}^{U_{n}}( x,y) +C\biggl( \frac{\rho}{\mathcal
{R}(
t) }\biggr) ^{\Theta}\sup_{y\in B( x,r)
\setminus\mathcal{%
N}}p_{t/2}^{U_{n}}( y,y)\nonumber
\end{eqnarray}
[note that if $n$ is large enough, then $B( x,r) \Subset U_{n}$].
Replacing on the right-hand side of (\ref{ptUn}) $p_{t}^{U_{n}}$ by a larger
value $p_{t}^{\Omega}$ and letting $n\rightarrow\infty$ on the left-hand
side, we obtain (\ref{oscpB}).

To prove (\ref{ptUn}), rename for simplicity $U_{n}$ into $U$ and recall
that, by construction in the proof of Theorem \ref{TptOm}, the domain
of $%
p_{t}^{U}$ is $U\setminus\mathcal{N}_{U}$ where~$\mathcal{N}_{U}$ is a
truly exceptional set in $U$, that is contained in $\mathcal{N}$. It follows
from Corollary \ref{Cesup=sup1}, that, for any $x\in U\setminus
\mathcal{N}$
\[
\sup_{y\in B( x,\rho) \setminus\mathcal
{N}}p_{t}^{U}(
x,y) =\limfunc{esup}_{y\in B( x,\rho)
}p_{t}^{U}(
x,y)
\]
and a similar identity for $\inf$ and $\limfunc{einf}$. Hence, it suffices
to prove that%
\[
\limfunc{eosc}_{y\in B( x,\rho) }p_{t}^{U}(
x,y) \leq
C\biggl( \frac{\rho}{\mathcal{R}( t) }\biggr) ^{\Theta}A,
\]
where%
\[
A=\sup_{y\in B( x,r) \setminus\mathcal
{N}}p_{t/2}^{U}(
y,y) .
\]
Set%
\[
u( y) =p_{t}^{U}( x,y) \quad\mbox{and}\quad
f(
y) =\partial_{t}p_{t}^{U}( x,y) ,
\]
where $\partial_{t}$ is the strong derivative in $L^{2}( U
) $
with respect to the variable~$y$. Applying Corollary \ref{Cdtpt} to
the heat
kernel $p_{t}^{U}$, we obtain, for $\mu$-a.a. $y\in
B(
x,r) $,%
\[
\vert f( y) \vert\leq\frac{2}{t}\sqrt{%
p_{t/2}^{U}( x,x) p_{t/2}^{U}( y,y) }\leq
\frac{2}{t}A.
\]
By Lemma \ref{LemGu}, we have $u=-G^{U}f$. Since for all $y\in B(
x,r) \setminus\mathcal{N}$%
\[
u( y) \leq\sqrt{p_{t/2}^{U}( x,x)
p_{t/2}^{U}(
y,y) }\leq A
\]
and $\rho\leq r$, we obtain by Proposition \ref{Posc} and $(
{E}_{F}\mbox{$\leq$}
) $ that%
\begin{eqnarray*}
\limfunc{eosc}_{B( x,\rho) }u &\leq&2\widetilde
{E}(
x,r) \limfunc{esup}_{B( x,r) }\vert f
\vert
+4\biggl( \frac{\rho}{r}\biggr) ^{\theta}\limfunc{esup}_{B(
x,r)
}\vert u\vert\\
&\leq&CF( r) \frac{A}{t}+4\biggl( \frac{\rho}{r}
\biggr) ^{\theta
}A.
\end{eqnarray*}
Since $r\leq\mathcal{R}( t) $, we have by (\ref{Fb})%
\[
F( r) \leq C\biggl( \frac{r}{\mathcal{R}( t
) }\biggr)
^{\beta}F( \mathcal{R}( t) ) =C\biggl(
\frac{r}{%
\mathcal{R}( t) }\biggr) ^{\beta}t,
\]
whence it follows that%
\[
\limfunc{eosc}_{B( x,\rho) }u\leq C\biggl( \biggl( \frac
{r}{%
\mathcal{R}( t) }\biggr) ^{\beta}+\biggl( \frac{\rho
}{r}\biggr)
^{\theta}\biggr) A.
\]
Note that this inequality is true for any $r$ such that $B(
x,r)
\Subset U$ and $\rho\leq r\leq\mathcal{R}( t) $.
Choosing $%
r=( \mathcal{R}( t) ^{\beta}\rho^{\theta}
) ^{{1%
}/({\beta+\theta})}$, we obtain (\ref{oscpB}).
\end{pf}
\begin{theorem}
\label{THolder}Let the hypotheses $\mbox{(\ref{VD})} +
\mbox{(\ref{condH})} + \mbox{(\ref{EFF})}$ be satisfied. Then, for any open
set $\Omega\subset
M$, the
heat kernel $p_{t}^{\Omega}( x,y) $ is H\"{o}lder
continuous in $%
x$ and $y$ in $\Omega\setminus\mathcal{N}$.
\end{theorem}
\begin{pf}
Fix $x\in\Omega\setminus\mathcal{N}$, $t>0$, and choose $\rho>0$
so small
that $B( x,r) \Subset\Omega$ where $r=r( t,\rho
) $
is defined by (\ref{r=}). Using Theorem \ref{TDUE}, (\ref{VD}), and
(\ref{Fb}), we obtain that, for any $y\in B( x,r)
\setminus
\mathcal{N}$,%
\begin{eqnarray*}
p_{t/2}^{\Omega}( y,y) &\leq&p_{t/2}( y,y)
\leq
\frac{C}{V( y,\mathcal{R}( t) ) } \\
&=&\frac{C}{V( x,\mathcal{R}( t) ) }\frac
{V( x,%
\mathcal{R}( t) ) }{V( y,\mathcal{R}(
t)
) } \\
&\leq&\frac{C}{V( x,\mathcal{R}( t) )
}\biggl( 1+\frac{%
d( x,y) }{\mathcal{R}( t) }\biggr) ^{\alpha
} \\
&\leq&\frac{C}{V( x,\mathcal{R}( t) ) },
\end{eqnarray*}
where we have used that $d( x,y) <r\leq\mathcal{R}(
t) $. Therefore, by Lemma \ref{LemoscpB},%
%
%e5.15 ###
%
\begin{equation} \label{oscptOm}
\limfunc{osc}_{y\in B( x,\rho) \setminus\mathcal{N}%
}p_{t}^{\Omega}( x,y) \leq\biggl( \frac{\rho}{\mathcal
{R}(
t) }\biggr) ^{\Theta}\frac{C}{V( x,\mathcal{R}(
t)
) }.
\end{equation}
In particular, if $y\in\Omega\setminus\mathcal{N}$ is close enough
to $x$,
then we have%
\[
\vert p_{t}^{\Omega}( x,x) -p_{t}^{\Omega}(
x,y) \vert\leq\biggl( \frac{d( x,y)
}{\mathcal{R}%
( t) }\biggr) ^{\Theta}\frac{C}{V( x,\mathcal
{R}(
t) ) },
\]
which means that $p_{t}^{\Omega}( x,\cdot) $ is H\"{o}lder
continuous in $\Omega\setminus\mathcal{N}$.
\end{pf}
\begin{corollary}
\label{CorHolder}Let the hypotheses $\mbox{(\ref{VD})} +
\mbox{(\ref{condH})}
+ \mbox{(\ref{EFF})}$ be satisfied, and let $B( x,R) $
be a
precompact ball, such that $x\in M\setminus\mathcal{N}$. Then for all
$%
\rho$ and $t$, such that
%
%e5.16 ###
%
\begin{equation} \label{troR}
0<\rho\leq\mathcal{R}( t) <R,
\end{equation}
and for all $y\in B( x,\rho) \setminus\mathcal{N}$, the
following estimate holds:
%
%e5.17 ###
%
\begin{equation} \label{pt-pt}
\bigl\vert p_{t}^{B( x,R) }( x,x)
-p_{t}^{B(
x,R) }( x,y) \bigr\vert\leq\biggl( \frac{\rho
}{\mathcal{R%
}( t) }\biggr) ^{\Theta}\frac{C}{V( x,\mathcal
{R}(
t) ) }.
\end{equation}
\end{corollary}
\begin{pf}
Set $\Omega=B( x,R) $. Then the condition $B(
x,r)
\Subset\Omega$ from the previous proof is satisfied because $r\leq
\mathcal{R}( t) <R$ by (\ref{r=}) and (\ref{troR}).
Hence,~(\ref{pt-pt}) follows from (\ref{oscptOm}).
\end{pf}

%s5.4 ###
\subsection{Proof of the lower bounds}
\label{Seclow}

\begin{lemma}
\label{LemDLE}Assume that $\mbox{(\ref{VD})} + \mbox{(\ref{EFF})}$ are
satisfied. Then there exists $\varepsilon>0$ such that, for all precompact
balls $B( x,R) $ with $x\in M\setminus\mathcal{N}$ and
for all $%
0<t\leq\varepsilon F( R) $,%
%
%e5.18 ###
%
\begin{equation} \label{pB>}
p_{t}^{B( x,R) }( x,x) \geq\frac{c}{V
( x,%
\mathcal{R}( t) ) }.
\end{equation}
\end{lemma}
\begin{pf}
Choose $r$ from the condition $t=\varepsilon F( r) $ which
implies $R\geq r$ and, hence, $p_{t}^{B( x,R) }\geq
p_{t}^{B( x,r) }$. Hence, it suffices to prove that%
\[
p_{t}^{B( x,r) }( x,x) \geq\frac{c}{V
( x,%
\mathcal{R}( t) ) }.
\]
Setting $B=B( x,r) $, we have by (\ref{PtOm})%
\[
\int_{B\setminus\mathcal{N}}p_{t}^{B}( x,y) \,d\mu
( y)
=\mathcal{P}_{t}^{B}1=\mathbb{P}_{x}( t<\tau_{B})
=1-\mathbb{P}%
_{x}( \tau_{B}\leq t) .
\]
By Corollary \ref{CEF}, we obtain%
\[
\mathbb{P}_{x}( \tau_{B}\leq t) \leq C\exp\biggl(
-c\biggl( \frac{%
F( r) }{t}\biggr) ^{{1}/({\beta^{\prime}-1})}
\biggr) =C\exp
\bigl( -c\varepsilon^{-{1}/({\beta^{\prime}-1})}\bigr) ,\vadjust{\goodbreak}
\]
whence it follows that, for small enough $\varepsilon\in(
0,1) $%
,%
\[
\int_{B}p_{t}^{B}( x,y) \,d\mu( y) \geq\frac{1}{2}.
\]
Therefore,%
\begin{eqnarray*}
p_{2t}^{B}( x,x) &=&\int_{B\setminus\mathcal
{N}}p_{t}^{B}(
x,y) ^{2}\,d\mu( y) \\
&\geq&\frac{1}{\mu( B) }\biggl( \int_{B\setminus
\mathcal{N}%
}p_{t}^{B}( x,y) \,d\mu( y) \biggr) ^{2} \\
&\geq&\frac{1/4}{V( x,\mathcal{R}( t/\varepsilon
) ) }%
,
\end{eqnarray*}
where we have used that $r=\mathcal{R}( t/\varepsilon) $.
Finally, using (\ref{Rb}) and (\ref{Va}) we obtain (\ref{pB>}).
\end{pf}
\begin{theorem}
\label{TNLE}Assume that the hypotheses $\mbox{(\ref{VD})} +
\mbox{(\ref{condH})}
+ \mbox{(\ref{EFF})}$ are satisfied, and all metric balls in $(
M,d) $ are precompact. Then there exist~$\varepsilon,\eta>0$ such
that\label{condLLE}%
%
%e5.19 ###
%
\begin{equation} \label{LLE}
p_{t}^{B( x,R) }( x,y) \geq\frac{c}{V
( x,%
\mathcal{R}( t) ) }
\end{equation}
for all $R>0$, $0<t\leq\varepsilon F( R) $ and $x,y\in
M\setminus\mathcal{N}$, provided
%
%e5.20 ###
%
\begin{equation} \label{dd}
d( x,y) \leq\eta\mathcal{R}( t)
\end{equation}
(see Figure \ref{pic9}). Consequently, the following
%
%f6 ###
%
\begin{figure}

\includegraphics{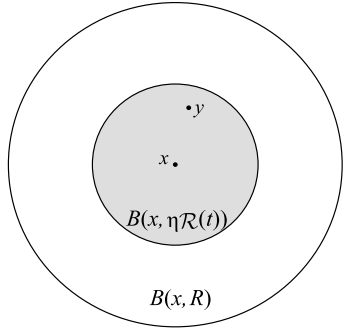}

\caption{Illustration to Theorem \protect\ref{TNLE}.}\label{pic9}
\end{figure}
inequality:\label{condNLE}%
{\renewcommand{\theequation}{\textit{NLE}}
\begin{equation}
p_{t}( x,y) \geq\frac{c}{V( x,\mathcal{R}(
t)
) }
\end{equation}}

\vspace*{-8pt}

\noindent
holds for all $t>0$ and $x,y\in M\setminus\mathcal{N}$ satisfying
(\ref{dd}).
\end{theorem}
\begin{pf}
Obviously, (\ref{NLE}) follows from (\ref{LLE})
by letting
$R\rightarrow\infty$, so that it suffices to prove (\ref{LLE}). Let
$\rho$
and $t$ be such that%
\[
0<\rho\leq\eta\mathcal{R}( t) \quad\mbox{and}\quad t\leq
\varepsilon F( R),
\]
where $\varepsilon\in( 0,1) $ is the constant from Lemma
\ref{LemDLE}, and $\eta\in( 0,1) $ is to be defined below.
Then the
hypotheses of Lemma \ref{LemDLE} and Corollary \ref{CorHolder} are
satisfied. Writing for simplicity $B=B( x,R) $, we obtain
by (\ref{LLE}) and (\ref{pt-pt}) that, for any $y\in B( x,\rho)
\setminus\mathcal{N}$,%
\begin{eqnarray*}
p_{t}^{B}( x,y) &\geq&p_{t}^{B}( x,x) -
\vert
p_{t}^{B}( x,x) -p_{t}^{B}( x,y) \vert
\\
&\geq&\frac{c}{V( x,\mathcal{R}( t) )
}-\biggl( \frac{%
\rho}{\mathcal{R}( t) }\biggr) ^{\Theta}\frac
{C}{V( x,%
\mathcal{R}( t) ) } \\
&\geq&\frac{c-C\eta^{\Theta}}{V( x,\mathcal{R}(
t)
) }.
\end{eqnarray*}
Choosing $\eta$ sufficiently small, we obtain (\ref{LLE}).
\end{pf}

Combining Theorems \ref{TUE}, \ref{THolder} and \ref{TNLE}, we
obtain the
main result:
\begin{theorem}
\label{Tmain}If the hypotheses $\mbox{(\ref{VD})} +\mbox{(\ref{condH})}
+ \mbox{(\ref{EFF})}$ are satisfied and all metric balls are precompact,
then the
heat kernel exists, is H\"{o}lder continuous in $x,y$, and satisfies
(\ref{UEE}) and (\ref{NLE}).
\end{theorem}
\begin{example}
Under the hypotheses of Theorem \ref{Tmain}, assume that the volume
function $V( x,r) $ satisfies the uniform estimate%
\[
V( x,r) \simeq r^{\alpha}
\]
with some $\alpha>0$, and function $F$ be as follows:%
%
%e5.21 ###
%
\setcounter{equation}{20}
\begin{equation} \label{F12}
F( r) =\cases{
r^{\beta_{1}}, &\quad$r<1$, \cr
r^{\beta_{2}}, &\quad$r\geq1$,}
\end{equation}
where $\beta_{1}>\beta_{2}>1$.\label{remwhysuchexampleexists?} Then
\[
\mathcal{R}( t) =
\cases{
t^{1/\beta_{1}}, &\quad$t<1$, \cr
t^{1/\beta_{2}}, &\quad$t\geq1$,}
\]
and the heat kernel satisfies the estimate%
\[
p_{t}( x,y) \leq\frac{C}{V( x,\mathcal{R}(
t)
) }\simeq C
\cases{
t^{-\alpha/\beta_{1}}, &\quad$t<1$,\cr
t^{-\alpha/\beta_{2}}, &\quad$t\geq1$.}
\]
It follows that
%
%e5.22 ###
%
\begin{equation} \label{pupbad}
p_{t}( x,y) \leq Ct^{-\alpha/\beta}
\end{equation}
for any $\beta$ from the interval $\beta_{2}<\beta<\beta_{1}$. Let us
verify that the following upper bound fails:%
%
%e5.23 ###
%
\begin{equation} \label{ptupbad}
p_{t}( x,y) \leq Ct^{-\alpha/\beta}\exp\biggl( -\biggl(
\frac{%
r^{\beta}}{Ct}\biggr) ^{{1}/({\beta-1})}\biggr) ,
\end{equation}
where $r=d( x,y) $. Indeed, by (\ref{NLE}) we
have%
\[
p_{t}( x,y) \geq\frac{c}{V( x,\mathcal{R}(
t)
) }
\]
provided $r\leq\eta\mathcal{R}( t) $. Assuming that
$t<1$ and
setting $r=\eta\mathcal{R}( t) =\eta t^{1/\beta_{1}}$
we obtain%
\footnote{%
The existence of a couple $x,y$ with a prescribed distance $r=d(
x,y) $ can be guaranteed, provided the space $( M,d
) $ is
connected.}
%
%e5.24 ###
%
\begin{equation} \label{plowab}
p_{t}( x,y) \geq\frac{c}{t^{\alpha/\beta_{1}}},
\end{equation}
while it follows from (\ref{ptupbad}) that%
%
%e5.25 ###
%
\begin{equation}\label{pupab}
p_{t}( x,y) \leq\frac{C}{t^{\alpha/\beta}}\exp\bigl(
-c(
t^{{\beta}/{\beta_{1}}-1}) ^{{1}/({\beta-1})}\bigr) .
\end{equation}
Since $\beta/\beta_{1}<1$, the exponent of $t$ under the exponential is
negative so that the right-hand side of (\ref{pupab}) becomes as $%
t\rightarrow0$ much smaller than that of~(\ref{plowab}), which is a
contradiction.

Another way to see a contradiction is to observe that (\ref{ptupbad})
implies (\ref{EFF}) with function $F( r)
\simeq
r^{\beta}$ (cf. \cite{KigamiNash,GrigHuUpper}), which is
incompatible with (\ref{EFF}) with function (\ref{F12}) [although
this argument requires the conservativeness of $( \mathcal
{E},\mathcal{F%
}) $].

The conclusion is that in general (\ref{pupbad}) does not imply (\ref
{ptupbad}). For comparison, let us note that if $\beta=2$ and the
underlying space is a Riemannian manifold, then (\ref{pupbad}) does
imply (%
\ref{ptupbad}); cf. \cite{GrigSuper}.
\end{example}

%s6 ###
\section{Matching upper and lower bounds}
\label{SecTwo}

%s6.1 ###
\subsection{\texorpdfstring{Distance $d_{\varepsilon}$}{Distance d_epsilon}}
\label{Secde}

\begin{definition}
We say that a sequence $\{ x_{i}\} _{i=0}^{N}$ of
points in $%
M $ is an $\varepsilon$-\textit{chain} between points $x,y\in M$ if%
\[
x_{0}=x,\qquad x_{N}=y\quad \mbox{and}\quad d(x_{i},x_{i-1})<\varepsilon
\qquad\mbox{for all }i=1,2,\ldots,N.
\]
\end{definition}

One can view an $\varepsilon$-chain as a sequence of \textit{chained}
balls $%
\{ B_{i}\} _{i=0}^{N}$ of radii~$\varepsilon$, that
connect $x$
and~$y$; that is, the center of $B_{0}$ is $x$, the center of $B_{N}$
is~$y$%
, and the center of $B_{i}$ is contained in $B_{i-1}$ for any
$i=1,\ldots,N$
%
%f7 ###
%
\begin{figure}

\includegraphics{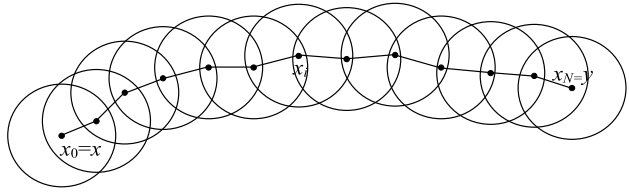}

\caption{An $\protect\varepsilon$-chain
connecting $x$ and $y$.}\label{pic8}
\end{figure}
(see Figure~\ref{pic8}).
\begin{definition}
For any $\varepsilon>0$ and all $x,y\in M$, define%
%
%e6.1 ###
%
\begin{equation}\label{dedef}
d_{\varepsilon}( x,y) =\inf_{\{ x_{i}\}
\ \mathrm{is}\ \varepsilon\mbox{{-}}\mathrm{chain}}\sum_{i=1}^{N}d(
x_{i},x_{i-1}),
\end{equation}
where the infimum is taken over all $\varepsilon$-chains $\{
x_{i}\} _{i=0}^{N}$ between $x,y$ with arbitrary~$N$.
\end{definition}

It is obvious that $d_{\varepsilon}( x,y) $ is a decreasing
left-continuous function of $\varepsilon$ and
%
%e6.2 ###
%
\begin{equation}\label{de>}
d_{\varepsilon}( x,y) \geq d( x,y) .
\end{equation}
Furthermore,
%
%e6.3 ###
%
\begin{equation} \label{eps>}
\varepsilon>d( x,y) \quad\Rightarrow\quad d_{\varepsilon}(
x,y) =d( x,y) .
\end{equation}
It is clear that $d_{\varepsilon}$ is an extended metric in the sense
that $%
d_{\varepsilon}$ satisfies all the axioms of a metric except for
finiteness. If an $\varepsilon$-chain exists for any two points $x,y$,
then $%
d_{\varepsilon}( x,y) <\infty$, and hence
$d_{\varepsilon}$ is
a true metric.
\begin{lemma}
\label{LemN}If $0<d_{\varepsilon}( x,y) <\infty$ for
some $%
x,y\in M$ and $\varepsilon>0$, then there exists an $\varepsilon
$-chain $%
\{ x_{i}\} _{i=0}^{N}$ between $x,y$ such that
%
%e6.4 ###
%
\begin{equation}\label{N<}
N\leq9\biggl\lceil\frac{d_{\varepsilon}( x,y) }{\varepsilon
}\biggr\rceil.
\end{equation}
\end{lemma}

Here \mbox{$\lceil\cdot\rceil$} stands for the least integer upper bound of the
argument. It follows from (\ref{dedef}) by the triangle inequality that
always
\[
N\geq\biggl\lceil\frac{d_{\varepsilon}( x,y) }{\varepsilon
}\biggr\rceil.
\]
Hence, denoting by $N_{\varepsilon}( x,y) $ the minimal
value of
$N$ for which there exists an $\varepsilon$ chain $\{
x_{i}\}
_{i=0}^{N}$ between $x$ and $y$, we obtain%
%
%e6.5 ###
%
\begin{equation} \label{Nede}
N_{\varepsilon}( x,y) \simeq\biggl\lceil\frac{d_{\varepsilon
}(
x,y) }{\varepsilon}\biggr\rceil.
\end{equation}
The number $N_{\varepsilon}( x,y) $ can be also viewed as the
minimal number in a~sequence of chained balls of radii $\varepsilon$
connecting $x$ and $y$.\vadjust{\goodbreak}
\begin{pf*}{Proof of Lemma \protect\ref{LemN}}
If $d_{\varepsilon}( x,y) <\varepsilon$, then also
$d(
x,y) <\varepsilon$, and hence $\{ x,y\} $ is an $%
\varepsilon$-chain with $N=1$. Assume $d_{\varepsilon}(
x,y)
\geq\varepsilon$, and let $\{ x_{i}\} _{i=0}^{n}$ be a~$%
\varepsilon$-chain between $x,y$, such that%
%
%e6.6 ###
%
\begin{equation} \label{xixi+1}
\sum_{i=1}^{n}d( x_{i},x_{i-1}) \leq2d_{\varepsilon
}(
x,y) ,
\end{equation}
which exists by hypothesis. Set $r_{i}=d( x_{i},x_{i-1})
$. Then (%
\ref{xixi+1}) implies%
\[
\#\{ i\dvtx r_{i}\geq\varepsilon/2\} \leq\frac
{4d_{\varepsilon
}( x,y) }{\varepsilon},
\]
whence%
\[
\#\{ i\dvtx r_{i}<\varepsilon/2\} \geq n-\frac
{4d_{\varepsilon
}( x,y) }{\varepsilon}.
\]
If $n>9\lceil\frac{d_{\varepsilon}( x,y) }{\varepsilon
}\rceil$,
then $n>9$ and $n>9\frac{d_{\varepsilon}( x,y)
}{\varepsilon}$,
whence it follows that%
\[
\#\{ i\dvtx r_{i}<\varepsilon/2\} >\frac{5n}{9}>\frac{n+1}{2}.
\]
Hence, there is an index $i$ such that both $r_{i-1}$ and $r_{i}$ are
smaller than $\varepsilon/2$. This implies that $d(
x_{i-1},x_{i+1}) <\varepsilon$ so that by removing the point $x_{i}$
from the chain we still have an $\varepsilon$-chain. Continuing this way,
we finally obtain an $\varepsilon$-chain satisfying~(\ref{N<}).
\end{pf*}

%s6.2 ###
\subsection{Two-sided estimates of the heat kernel}

If $x\neq y$, then it follows from~(\ref{de>}) and (\ref{Fb}) that%
%
%e6.7 ###
%
\begin{equation}\label{Fei}
\frac{F( \varepsilon) }{\varepsilon}d_{\varepsilon
}(x,y) \rightarrow\infty\qquad\mbox{as }\varepsilon\rightarrow
\infty.
\end{equation}
In this section, we make an additional assumption that, for all $x,y\in
M$,%
%
%e6.8 ###
%
\begin{equation}\label{Fe}
\frac{F( \varepsilon) }{\varepsilon}d_{\varepsilon
}(x,y) \rightarrow0 \qquad\mbox{as }\varepsilon\rightarrow0.
\end{equation}
In particular, (\ref{Fe}) implies the finiteness of $d_{\varepsilon
}$ for
all $\varepsilon>0$. Define the function $\varepsilon(
t,x,y) $
as follows:%
%
%e6.9 ###
%
\begin{equation} \label{epsdef}
\varepsilon( t,x,y) =\sup\biggl\{ \varepsilon>0\dvtx\frac
{F(
\varepsilon) }{\varepsilon}d_{\varepsilon}( x,y)
\leq t\biggr\}.
\end{equation}
If $x=y$, then $\varepsilon( t,x,x) =\infty$. If $x\neq
y$, then
it follows from (\ref{Fe}) and (\ref{Fei}) that $0<\varepsilon
(
t,x,y) <\infty$.
\begin{lemma}
\label{LemFe}If (\ref{Fe}) is satisfied, then the function $
\varepsilon( t,x,y) $ satisfies the identity%
%
%e6.10 ###
%
\begin{equation} \label{te}
\frac{F( \varepsilon) }{\varepsilon}d_{\varepsilon
}(
x,y) =t
\end{equation}
for all distinct $x,y\in M$ and $t>0$.
\end{lemma}
\begin{pf}
Since the function $F( \varepsilon) $ is continuous and $%
d_{\varepsilon}( x,y) $ is left-contin\-uous in
$\varepsilon$, we
have
\[
\frac{F( \varepsilon) }{\varepsilon}d_{\varepsilon
}(
x,y) \leq t.
\]
Assume from the contrary that
\[
\frac{F( \varepsilon) }{\varepsilon}d_{\varepsilon
}(
x,y) <t,
\]
and note that, for any $\varepsilon^{\prime}>\varepsilon$, we have by
definition of $\varepsilon$ that%
%
%e6.11 ###
%
\begin{equation} \label{eps}
\frac{F( \varepsilon^{\prime}) }{\varepsilon^{\prime}}
d_{\varepsilon^{\prime}}( x,y) >t.
\end{equation}
On the other hand, $d_{\varepsilon^{\prime}}( x,y) \leq
d_{\varepsilon}( x,y) $ and
\[
\frac{F( \varepsilon^{\prime}) }{\varepsilon^{\prime}}
\rightarrow\frac{F( \varepsilon) }{\varepsilon} \qquad\mbox
{as }%
\varepsilon^{\prime}\rightarrow\varepsilon+.
\]
Hence,%
\[
\limsup_{\varepsilon^{\prime}\rightarrow\varepsilon+}\frac{F(
\varepsilon^{\prime}) }{\varepsilon^{\prime}}\leq\frac
{F(
\varepsilon) }{\varepsilon}d_{\varepsilon}( x,y) <t,
\]
which contradicts (\ref{eps}).
\end{pf}
\begin{theorem}
\label{Ttwo}Let all metric balls be precompact. Let the hypotheses
$\mbox{(\ref{VD})} + \mbox{(\ref{EFF})} +\mbox{(\ref{condH})} $ and (\ref
{Fe}) be
satisfied, and let $\varepsilon( t,x,y) $ be the function from
(\ref{epsdef}). Then, for all $x,y\in M\setminus\mathcal{N}$
and $%
t>0 $, we have
%
%e6.12 ###
%
\begin{equation}\label{two}
p_{t}( x,y) \asymp\frac{C}{V(x,\mathcal{R}(
t) )}\exp
( -c\Phi( cd_{\varepsilon}( x,y) ,t)) ,
\end{equation}
where $\varepsilon=\varepsilon( \kappa t,x,y) $ and
$\kappa=8$
for the upper bound in (\ref{two}) while $\kappa$ is a~small enough
positive constant for the lower bound.
\end{theorem}

The proof of Theorem \ref{Ttwo} is preceded by a lemma.
\begin{lemma}
\label{Lemdee}For all distinct $x,y\in M$ and $t>0$, we have%
%
%e6.13 ###
%
\begin{equation}\label{dF}
\Phi( cd_{\varepsilon}( x,y) ,t) \leq\frac
{%
d_{\varepsilon}( x,y) }{\varepsilon}\leq\Phi(
Cd_{\varepsilon}( x,y) ,t) ,
\end{equation}
where $\varepsilon=\varepsilon( t,x,y) $.
\end{lemma}
\begin{pf}
Let us first show that, for all $\varepsilon>0$ and some $c\in(
0,1) $,
%
%e6.14 ###
%
\begin{equation} \label{FiF}
\Phi\biggl( c\frac{\varepsilon}{F( \varepsilon)
}\biggr) \leq
\frac{1}{F( \varepsilon) }\leq\Phi\biggl( 2\frac
{\varepsilon}{%
F( \varepsilon) }\biggr) .
\end{equation}
By (\ref{Fidef1}), we have, for all $r>0$,%
\[
\Phi\biggl( \frac{2\varepsilon}{F( \varepsilon)
}\biggr) \geq
\frac{2\varepsilon}{F( \varepsilon) r}-\frac{1}{F(
r) }.
\]
Choosing $r=\varepsilon$ we obtain the right-hand side inequality in
(\ref{FiF}). By (\ref{Fidef1}), the left-hand side inequality in (\ref
{FiF}) is
equivalent to%
\[
\frac{c\varepsilon}{F( \varepsilon) r}-\frac{1}{F(
r) }\leq\frac{1}{F( \varepsilon) } \qquad\mbox{for
all }r>0,
\]
that is, to%
\[
\frac{F( \varepsilon) }{F( r) }\geq\frac{%
c\varepsilon}{r}-1.
\]
If $r\geq\varepsilon$, then this is trivially satisfied provided
$c\leq1$.
If $r<\varepsilon$, then by (\ref{Fb}) we have%
\[
\frac{F( \varepsilon) }{F( r) }\geq c
\biggl( \frac{%
\varepsilon}{r}\biggr) ^{\beta}\geq c\frac{\varepsilon}{r},
\]
which proves the previous inequality and, hence, (\ref{FiF}).

Putting in (\ref{FiF}) $\varepsilon=\varepsilon( t,x,y)
$ and
using $\frac{\varepsilon}{F( \varepsilon) }=\frac{%
d_{\varepsilon}( x,y) }{t}$, which is true by Lemma \ref{LemFe}%
, we obtain%
\[
\Phi\biggl( c\frac{d_{\varepsilon}( x,y) }{t}\biggr)
\leq\frac{%
d_{\varepsilon}( x,y) }{\varepsilon t}\leq\Phi\biggl(
2\frac{%
d_{\varepsilon}( x,y) }{t}\biggr) ,
\]
whence (\ref{dF}) follows.
\end{pf}
\begin{pf*}{Proof of Theorem \protect\ref{Ttwo}}
If $x=y$, then $d_{\varepsilon}( x,y) =0$, and (\ref{two}) follows
from Theorems \ref{TUE} and \ref{TNLE}. Assume in the sequel that
$x\neq y$.
Let us first prove the lower bound in (\ref{two}), that is,%
%
%e6.15 ###
%
\begin{equation}\label{twol}
p_{t}( x,y) \geq\frac{c}{V(x,\mathcal{R}( t
) )}\exp
( -C\Phi( Cd_{\varepsilon}( x,y) ,t)) .
\end{equation}
By Theorem \ref{TNLE}, we have%
{\renewcommand{\theequation}{\textit{NLE}}
\begin{equation}
p_{t}( x,y) \geq\frac{c}{V( x,\mathcal{R}(
t)
) }
\end{equation}}

\vspace*{-8pt}

\noindent
provided
%
%e6.16 ###
%
\setcounter{equation}{15}
\begin{equation} \label{detR}
d( x,y) \leq\eta\mathcal{R}( t)
\end{equation}
for some $\eta>0$. Set $\varepsilon=\varepsilon( \kappa
t,x,y)
$ where $\kappa\in( 0,1) $ will be chosen later.

Consider first the case $\varepsilon\geq d_{\varepsilon}(
x,y) $%
. By (\ref{dF}), we have%
\[
\Phi( cd_{\varepsilon}( x,y) ,\kappa t)
\leq1.
\]
Applying (\ref{tFi}) with $R=cd_{\varepsilon}( x,y) $,
we obtain%
\[
F( cd_{\varepsilon}( x,y) ) \leq C\kappa t,\vadjust{\goodbreak}
\]
whence by (\ref{Rb})%
\[
d_{\varepsilon}( x,y) \leq c^{-1}\mathcal{R}(
C\kappa
t) \leq\eta\mathcal{R}( t) ,
\]
provided $\kappa$ is small enough. Since $d( x,y) \leq
d_{\varepsilon}( x,y) $, we see that the condition~(\ref{detR})
is satisfied and, hence, (\ref{twol}) follows from (\ref{NLE}).

Assume now that $\varepsilon<d_{\varepsilon}( x,y) $. By
Lemma %
\ref{LemN}, there is an $\varepsilon$-chain $\{ x_{i}\}
_{i=1}^{N}$ connecting $x$ and $y$ and such that%
%
%e6.17 ###
%
\begin{equation}\label{Ne}
N\simeq\frac{d_{\varepsilon}( x,y) }{\varepsilon}.
\end{equation}
By (\ref{semi}), we have%
%
%e6.18 ###
%
\begin{eqnarray} \label{chain}
p_{t}( x,y) &=&\int_{M}\cdots\int
_{M}p_{{t/N}%
}(x,z_{1})p_{{t/N}}(z_{1},z_{2})\cdots\nonumber\\
&&\hspace*{41.5pt}{}\times p_{{t/N}%
}(z_{N-1},y)\,dz_{1}\cdots dz_{N-1} \nonumber\\[-8pt]\\[-8pt]
&\geq&\int_{B_{1}}\cdots\int_{B_{N-1}}p_{{t/N}
}(z_{0},z_{1})p_{{t/N}}(z_{1},z_{2})\cdots\nonumber\\
&&\hspace*{54.8pt}{}\times  p_{{t/N}%
}(z_{N-1},z_{N})\,dz_{1}\cdots dz_{N-1},\nonumber
\end{eqnarray}
where $z_{0}=x$, $z_{N}=y$, $B_{i}=B( x_{i},\varepsilon)
$. We
will estimate $p_{t/N}( z_{i},z_{i+1}) $ from below by
means of (\ref{NLE}). For that, we need to verify the condition%
\[
d( z_{i},z_{i+1}) \leq\eta\mathcal{R}( t/N
) .
\]
By (\ref{te}), we have
%
%e6.19 ###
%
\begin{equation} \label{Tze}
\frac{F( \varepsilon) }{\varepsilon}=\frac{\kappa t}{%
d_{\varepsilon}( x,y) }.
\end{equation}
It follows from (\ref{Tze}) and (\ref{Ne}) that%
\[
F( \varepsilon) \simeq\frac{\kappa t}{N},
\]
whence by (\ref{Rb})
\[
\varepsilon\simeq\mathcal{R}\biggl( \frac{\kappa t}{N}\biggr) .
\]
Clearly, if $\kappa$ is small enough, then%
%
%e6.20 ###
%
\begin{equation}\label{ee}
\varepsilon\leq\frac{\eta}{3}\mathcal{R}\biggl( \frac{t}{N}\biggr) .
\end{equation}
Since in (\ref{chain}) $z_{i}\in B( x_{i},\varepsilon) $
and $%
d( x_{i},x_{i+1}) \leq\varepsilon$, it follows from
(\ref{ee})
that%
\[
d( z_{i},z_{i+1}) \leq3\varepsilon\leq\eta\mathcal
{R}(
t/N) .
\]
Hence, by (\ref{NLE}) and (\ref{Va}),%
\[
p_{{t/N}}( z_{i},z_{i+1}) \geq\frac{c}{V(
z_{i},%
\mathcal{R}( t/N) ) }\geq\frac{c}{V(
x_{i},\mathcal{R}%
( t/N) ) }\geq\frac{c}{V( x_{i},\varepsilon
) }.\vadjust{\goodbreak}
\]
Therefore, (\ref{chain}) implies%
%
%e6.21 ###
%
\begin{eqnarray}\label{4}
p_{t}( x,y) &\geq&\frac{c}{V( x,\mathcal{R}(
t/N) ) }\prod_{i=1}^{N-1}\frac{c}{V(
x_{i},\varepsilon
) }V( x_{i},\varepsilon) \nonumber\vadjust{\goodbreak}\\
&\geq&\frac{c^{-N}}{V( x,\mathcal{R}( t/N)
) } \nonumber\\[-8pt]\\[-8pt]
&\geq&\frac{\exp( -CN) }{V( x,\mathcal{R}(
t)
) } \nonumber\\
&\geq&\frac{\exp( -C({d_{\varepsilon}( x,y)
})/{%
\varepsilon}) }{V( x,\mathcal{R}( t)
) }.\nonumber
\end{eqnarray}
By Lemma \ref{Lemdee}, we have%
\[
\frac{d_{\varepsilon}( x,y) }{\varepsilon}\leq\Phi
(
Cd_{\varepsilon}( x,y) ,\kappa t) =\kappa\Phi
\biggl(
\frac{C}{\kappa}d_{\varepsilon}( x,y) ,t\biggr) .
\]
Substituting into (\ref{4}), we obtain (\ref{twol}).

To prove the upper bound in (\ref{two}), we basically repeat the proof of
Theorem~\ref{TUE} with $d$ being replaced by $d_{\varepsilon}$ for an
appropriate $\varepsilon$. Fix some $\varepsilon>0$ and denote by $%
B_{\varepsilon}( x,r) $ the ball in the metric
$d_{\varepsilon}$%
. It follows from (\ref{eps>}) that
%
%e6.22 ###
%
\begin{equation} \label{B=B}
B_{\varepsilon}( x,r) =B( x,r) \qquad\mbox{for all }%
r\leq\varepsilon,
\end{equation}
which allows to modify Lemma \ref{pf<exp} as follows: for all $x\in
M\setminus\mathcal{N}_{0}$ and $R>0$%
%
%e6.23 ###
%
\begin{equation} \label{Rr}
\mathbb{E}_{x}\bigl( e^{-\lambda\tau_{B_{\varepsilon}(x,R)}}
\bigr) \leq
C\exp\biggl( -c\frac{R}{r}\biggr) ,
\end{equation}
provided the values of parameters $r$ and $\lambda$ satisfy the
conditions%
%
%e6.24 ###
%
\begin{equation}\label{rl}
0<r\leq\varepsilon\quad\mbox{and}\quad\lambda\geq\frac{\sigma}{
F( r) }.
\end{equation}
Indeed, (\ref{Rr}) is analogous to estimate (\ref{n}) from the proof
of Lemma \ref{pf<exp} for $d$-balls, which was proved using $\lambda
\geq
\frac{\sigma}{F( r) }$. To repeat\vspace*{1pt} the proof for the
metric $%
d_{\varepsilon}$, we need the precompactness of $d_{\varepsilon}$-balls,
that follows from (\ref{de>}), and the condition (\ref{EFF}) for $%
d_{\varepsilon}$-balls of radii $\leq r$, that follows from (\ref{B=B}),
provided $r\leq\varepsilon$.

Consequently, the statement of Theorem \ref{TEF} is modified as
follows: for
all $x\in M\setminus\mathcal{N}_{0}$ and $R,t>0$,
%
%e6.25 ###
%
\begin{equation}\label{PG}
\mathbb{P}_{x}\bigl( \tau_{B_{\varepsilon}( x,R) }\leq
t\bigr)
\leq C\exp( -c\Phi_{\varepsilon}( cR,t) ),
\end{equation}
where
%
%e6.26 ###
%
\begin{equation} \label{G}
\Phi_{\varepsilon}( R,t) :=\sup_{0<r\leq\varepsilon
}\biggl\{
\frac{R}{r}-\frac{t}{F( r) }\biggr\}.
\end{equation}
Indeed, arguing as in (\ref{3a}) and using (\ref{Rr}) we obtain
that, under the assumptions of (\ref{rl}),%
\[
\mathbb{P}_{x}\bigl( \tau_{B_{\varepsilon}( x,R) }\leq
t\bigr)
\leq C\exp\biggl( -c\frac{R}{r}+\lambda t\biggr) .
\]
Setting here $\lambda=\sigma/F( r) $ yields%
%
%e6.27 ###
%
\begin{equation} \label{g}
\mathbb{P}_{x}\bigl( \tau_{B_{\varepsilon}( x,R) }\leq
t\bigr)
\leq C\exp\biggl( -\biggl( c\frac{R}{r}-\frac{\sigma t}{F(
r) }%
\biggr) \biggr) .
\end{equation}
Finally, minimizing the right-hand side of (\ref{g}) in $r\leq
\varepsilon$%
, we obtain (\ref{PG}).

Let us show that if
%
%e6.28 ###
%
\begin{equation}\label{te2}
t\leq\frac{1}{2}\frac{F( \varepsilon) }{\varepsilon}R,
\end{equation}
then%
%
%e6.29 ###
%
\begin{equation}\label{p}
\Phi( R,t) \leq2\Phi_{\varepsilon}( R,t) .
\end{equation}
We have%
\[
\sup_{r>\varepsilon}\biggl\{ \frac{R}{r}-\frac{t}{F( r
) }\biggr\}
\leq\frac{R}{\varepsilon},
\]
whereas%
\[
\sup_{0<r\leq\varepsilon}\biggl\{ \frac{R}{r}-\frac{t}{F(
r) }%
\biggr\} \geq\frac{R}{\varepsilon}-\frac{t}{F( \varepsilon
) }%
\geq\frac{R}{\varepsilon}-\frac{1}{2}\frac{R}{\varepsilon}=\frac
{1}{2}%
\frac{R}{\varepsilon}.
\]
It follows that%
\[
\Phi( R,t) =\sup_{r>0}\biggl\{ \frac{R}{r}-\frac
{t}{F(
r) }\biggr\} \leq2\sup_{0<r\leq\varepsilon}\biggl\{ \frac
{R}{r}-%
\frac{t}{F( r) }\biggr\},
\]
which proves (\ref{p}). Hence, we can rewrite (\ref{PG}) in the form
%
%e6.30 ###
%
\begin{equation} \label{PG1}
\mathbb{P}_{x}\bigl( \tau_{B_{\varepsilon}( x,R) }\leq
t\bigr)
\leq C\exp( -c\Phi( cR,t) ),
\end{equation}
provided the relation (\ref{te2}) between $\varepsilon,t,R$ is satisfied.

As in the last part of the proof of Theorem \ref{TUE}, we apply (\ref{PG1})
with $R=\frac{1}{4}d_{\varepsilon}( x,y) $ for fixed
$x,y\in
M\setminus\mathcal{N}$. Note that in (\ref{VRR}) $d( x,z
) $ can
be replaced by a larger value $d_{\varepsilon}( x,z) $.
The rest
of the argument goes through unchanged, and we obtain
%
%e6.31 ###
%
\begin{equation} \label{ptupe}
p_{t}( x,y) \leq\frac{C}{V( x,\mathcal{R}(
t)
) }\exp( -c\Phi( cd_{\varepsilon}(
x,y)
,t) ) ,
\end{equation}
provided
%
%e6.32 ###
%
\begin{equation} \label{18}
t\leq\frac{1}{8}\frac{F( \varepsilon) }{\varepsilon}%
d_{\varepsilon}( x,y) .
\end{equation}
By (\ref{te}), condition (\ref{18}) can be satisfied by setting $%
\varepsilon=\varepsilon( 8t,x,y) $.
\end{pf*}
\begin{corollary}
\label{Ctwo}Under the hypotheses of Theorem \ref{Ttwo}, we have
%
%e6.34 ###
%e6.33 ###
%
\begin{eqnarray}
\label{twod}
p_{t}( x,y) &\asymp&\frac{C}{V(x,\mathcal{R}(
t) )}%
\exp\biggl( -c\frac{d_{\varepsilon}( x,y) }{\varepsilon
}\biggr)
\\
\label{twoN}
&\asymp&\frac{C}{V(x,\mathcal{R}( t) )}\exp(
-cN_{\varepsilon}( x,y) ) ,
\end{eqnarray}
where $\varepsilon=\varepsilon( \kappa t,x,y) $ and
$\kappa$
is a large enough constant for the upper bound and a small enough positive
constant for the lower bound.
\end{corollary}
\begin{pf}
The lower bound in (\ref{twod}) follows from (\ref{4}), while the upper
bound follows from (\ref{ptupe}) and
\[
\frac{d_{\varepsilon}( x,y) }{\varepsilon}\leq\Phi
(
Cd_{\varepsilon}( x,y) ,\kappa t) =\kappa\Phi
\biggl(
\frac{C}{\kappa}d_{\varepsilon}( x,y) ,t\biggr) ,
\]
provided $\kappa$ is chosen large enough to ensure $C/\kappa\leq c$
where $%
c$ is the constant from (\ref{ptupe}). Estimate (\ref{twoN}) follows
then from (\ref{Nede}).
\end{pf}
\begin{remark}
\label{ExHK}A good example to illustrate Theorem \ref{Ttwo} and
Corollary~\ref{Ctwo} is the class of post critically finite (p.c.f.) fractals, where
the heat kernel estimate (\ref{twoN}) was proved by Hambly and Kumagai
\cite{HamblyKum}. Without going into the details of \cite{HamblyKum},
let us
mention that $d( x,y) $ is the resistance metric on such a
fractal $M$, and $\mu$ is the Hausdorff measure of $M$ of dimension
$\alpha
:=\dim_{H}M$. One has in this setting $V( x,r) \simeq
r^{\alpha
} $, in particular,~(\ref{VD}) is satisfied. Hambly and Kumagai
proved that (\ref{EFF}) is satisfied with $F(
r)
=r^{\beta}$ where $\beta=\alpha+1$; cf. \cite{HamblyKum}, Theorem 3.8.
Condition (\ref{Fe}) follows from the estimate%
%
%e6.35 ###
%
\begin{equation} \label{Ne<}
N_{\varepsilon}( x,y) \leq C\biggl( \frac{d(
x,y) }{%
\varepsilon}\biggr) ^{\gamma},
\end{equation}
proved in \cite{HamblyKum}, Lemma 3.3, with $\gamma=\beta/2$, as
(\ref{Ne<}%
) implies that%
\[
\frac{F( \varepsilon) }{\varepsilon}d_{\varepsilon
}(
x,y) \leq C\varepsilon^{\beta}N_{\varepsilon}(
x,y) \leq
Cd( x,y) ^{\gamma}\varepsilon^{\beta-\gamma
}\rightarrow0 %
\qquad\mbox{as }\varepsilon\rightarrow0.
\]
The Harnack inequality~(\ref{condH}) on p.c.f. fractals was
proved by
Kigami~\cite{Kigamibook}. Hence, Corollary~\ref{Ctwo} applies and
gives on
unbounded p.c.f. fractals estimate~(\ref{twoN}). The same estimate was
proved in~\cite{HamblyKum}, Theorem 1.1, for bounded p.c.f. fractals
using a
different method.

Note that (\ref{VD}) implies (\ref{Ne<}) with $\gamma
=\alpha$
[where $\alpha$ comes from (\ref{Va})], provided all balls in $M$ are
connected. Indeed, (\ref{VD}) implies by the classical ball
covering argument that any ball of radius $r$ can be covered by at
most~$C( \frac{r}{\varepsilon}) ^{\alpha}$ balls of radii $\varepsilon\in( 0,r)$.
Consequently, any point $y\in B(x,r) $ can be connected to $x$ by a chain
of $\varepsilon$-balls containing at most $C( \frac{r}{\varepsilon}) ^{\alpha}$ balls.
Taking $r\simeq d( x,y) $ we obtain (\ref{Ne<}) with
$\gamma
=\alpha$. Therefore, hypothesis~(\ref{Fe}) is satisfied automatically
for $F( r) =r^{\beta}$ with $\beta>\alpha$.\vadjust{\goodbreak}

Estimate (\ref{twoN}) means that the diffusion process goes from $x$
to $%
y$ in time~$t$ in the following way. The process first
``computes'' the value\footnote{%
For example, in the above setting, when (\ref{Ne<}) is satisfied with $
\gamma<\beta$, we obtain from~(\ref{te}) $\varepsilon^{\beta
}N_{\varepsilon}\simeq t$ whence%
\[
\varepsilon^{\beta}\biggl( \frac{d( x,y) }{\varepsilon
}\biggr)
^{\gamma}\geq ct
\]
and%
\begin{eqnarray*}
\varepsilon\geq c\biggl( \frac{t}{d( x,y) ^{\gamma
}}\biggr) ^{%
{1}/({\beta-\gamma})}.\\[-25pt]
\end{eqnarray*}}
of $\varepsilon$ as a~function of $t,x,y$, then
``detects'' a shortest chain of $\varepsilon$-balls
connecting $x$ and $y$ and finally goes along that chain (see Figure
\ref{pic10}).

%
%f8 ###
%
\begin{figure}

\includegraphics{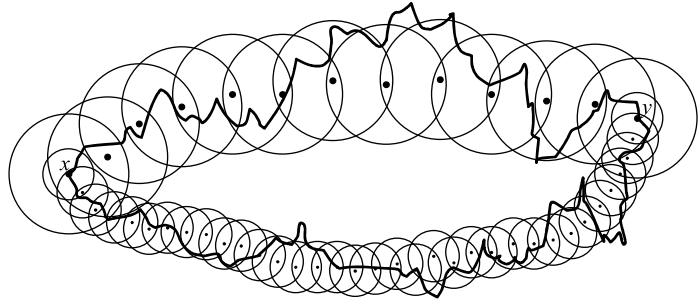}

\caption{Two shortest chains of
$\protect\varepsilon$-balls for two distinct values of $\protect
\varepsilon
$ provide different routes for the diffusion from $x$ to $y$ for two
distinct values of time $t$.}\label{pic10}
\end{figure}

This phenomenon was first observed by Hambly and Kumagai on p.c.f. fractals,
but it seems to be generic. Hence, to obtain matching upper and lower
bounds, one needs, in addition to the usual hypotheses, also the information
encoded in the function $N_{\varepsilon}( x,y) $, namely, the
graph distance between~$x$ and $y$ on any $\varepsilon$-net approximation
of $M$.
\end{remark}

%s6.3 ###
\subsection{Chain condition}
\label{SecChain}

The statement of Theorem \ref{Ttwo} can be simplified
if the
metric space $( M,d) $ possesses an additional property as
follows.
\begin{definition}
We say that a metric space $( M,d) $ satisfies the
\textit{%
chain condition} if there exists a constant $C\geq1$ such that, for any
positive integer~$n$ and for all $x,y\in M$, there is a sequence $
\{
x_{k}\} _{k=0}^{n}$ of points in $M$ such that $x_{0}=x$,
$x_{n}=y$ and%
\[
d( x_{k-1},x_{k}) \leq C\frac{d( x,y) }{n}
\qquad\mbox{for
all }k=1,\ldots,n.
\]
\end{definition}

For example, any geodesic metric satisfies the chain condition.
\begin{lemma}
\label{Lemde}If $( M,d) $ satisfies the chain condition,
then $%
d_{\varepsilon}\leq Cd$ for any $\varepsilon>0$.
\end{lemma}
\begin{pf}
Indeed, fix $\varepsilon>0$ and two distinct points $x,y\in M$, and
choose $%
n$ so big that $C\frac{d( x,y) }{n}<\varepsilon$. Let
$\{
x_{k}\} _{k=0}^{n}$ be a sequence from the chain condition. Then
it is
also an $\varepsilon$-chain, whence%
\[
d_{\varepsilon}( x,y) \leq\sum_{k=1}^{n}d(
x_{k-1},x_{k}) \leq Cd( x,y) ,
\]
which was to be proved.
\end{pf}
\begin{corollary}
\label{Cortwo}Let the metric space $( M,d) $ satisfy the chain
condition, and let all metric balls be precompact. If the hypotheses
$\mbox{(\ref{VD})} + \mbox{(\ref{EFF})} +\mbox{(\ref{condH})}$ are satisfied,
then, for
all $x,y\in M\setminus\mathcal{N}$ and $t>0$,
%
%e6.36 ###
%
\begin{equation}\label{twosided}
p_{t}( x,y) \asymp\frac{C}{V(x,\mathcal{R}(
t) )}\exp
\biggl( -ct\Phi\biggl( c\frac{d( x,y) }{t}\biggr)\biggr) ,
\end{equation}
where $\Phi( s) $ is defined by (\ref{Fidef1}).
\end{corollary}
\begin{pf}
Since by Lemma \ref{Lemde} $d_{\varepsilon}\leq Cd$, condition (\ref
{Fe}) is obviously satisfied. Since $d_{\varepsilon}\simeq d$, we can
replace in (\ref{two}) $d_{\varepsilon}$ by $d$, which together with~(\ref{FiFi}) yields (\ref{twosided}).
\end{pf}
\begin{remark}
\label{RemNLE}
$\!\!\!$Obviously, estimate (\ref{twosided}) (should it be
true) implies~(\ref{UEE}). We claim that (\ref{twosided}) implies
also (\ref{NLE}); moreover, the parameter $\eta$ in~(\ref{NLE})
can be chosen to be arbitrarily large, say $\eta>1$. Indeed,
we need to show that if%
\[
d( x,y) \leq\eta\mathcal{R}( t),
\]
where $\eta$ is a (large) given constant, then%
\[
t\Phi\biggl( c\frac{d( x,y) }{t}\biggr) \leq\func{const},
\]
which amounts to%
%
%e6.37 ###
%
\begin{equation} \label{txi}
\Phi\biggl( \eta\frac{\mathcal{R}( t) }{t}\biggr)
\leq\frac{%
\func{const}}{t},
\end{equation}
where we have renamed $c\eta$ to $\eta$. Indeed, by (\ref{Fidef1})
we have%
\[
\Phi( s) =\sup_{\xi>0}\biggl\{ \frac{s}{\mathcal
{R}( \xi
) }-\frac{1}{\xi}\biggr\}
\]
so that (\ref{txi}) is equivalent to%
%
%e6.38 ###
%
\begin{equation} \label{t/xi}
\frac{\eta\mathcal{R}( t) }{\mathcal{R}( \xi
) }\leq
\frac{t}{\xi}+\func{const.}\vadjust{\goodbreak}
\end{equation}
If $\xi\leq t$, then by (\ref{Rb})%
\[
\frac{\mathcal{R}( t) }{\mathcal{R}( \xi)
}\leq
C\biggl( \frac{t}{\xi}\biggr) ^{1/\beta}.
\]
If the ratio $\frac{t}{\xi}$ is large enough then, using $1/\beta
<1$, we
obtain that
\[
\eta C\biggl( \frac{t}{\xi}\biggr) ^{1/\beta}\leq\frac{t}{\xi},
\]
whence (\ref{t/xi}) follows. If $\frac{t}{\xi}$ is bounded by a constant,
say $\frac{t}{\xi}\leq C^{\prime}$ (which includes also the case
$\xi>t$),
then by (\ref{Rb})%
\[
\frac{\eta\mathcal{R}( t) }{\mathcal{R}( \xi
) }\leq
\eta\frac{\mathcal{R}( C^{\prime}\xi) }{\mathcal
{R}( \xi
) }\leq\func{const},
\]
whence (\ref{txi}) follows again.
\end{remark}

%s7 ###
\section{Consequences of heat kernel bounds}
\label{Secconv}

%s7.1 ###
\subsection{Harmonic function and the Dirichlet problem}
\label{SecHarm}

We assume only basic hypotheses in this subsection. Moreover,
we use neither the locality of $( \mathcal{E},\mathcal{F}
) $ nor
the existence of the process $\{ X_{t}\} $. We state and prove
some basic properties of the Dirichlet problem in the abstract setting, that
will be used in the proof of Theorem~\ref{Tconv}.

Fix an open set $\Omega\subset M$ such that $\lambda_{\min}(
\Omega
) >0$, and consider the following weak Dirichlet problem in
$\Omega$:
given a function $f\in\mathcal{F}$, find a function $u\in\mathcal
{F}$ such
that%
%
%e7.1 ###
%
\begin{equation} \label{D}
\cases{
u\mbox{ is harmonic in }\Omega, \cr
u=f \func{mod}\mathcal{F}( \Omega) ,}%
\end{equation}
where the second condition is a weak boundary condition and means that $
u-f\in\mathcal{F}( \Omega) $.
\begin{lemma}
\label{LemDir1}
\textup{(a)} For any $f\in\mathcal{F}$, problem
(\ref{D}) has a unique solution $u$.

\textup{(b)} If $u$ solves (\ref{D}) and $w\in
\mathcal{F}$ is
another function such that $w=f$ $\func{mod}\mathcal{F}( \Omega
),
$ then $\mathcal{E}( u) \leq\mathcal{E}( w
) $.
Moreover, the identity $\mathcal{E}( u) =\mathcal
{E}(
w) $ holds if and only if $u=w$.
\end{lemma}
\begin{pf}
(a) The condition $\lambda_{\min}( \Omega
) >0$
implies that%
\[
\mathcal{E}( \varphi) \simeq\mathcal{E}( \varphi
)
+\Vert\varphi\Vert_{2}^{2}
\]
for all $\varphi\in\mathcal{F}( \Omega) $. Hence,
$\mathcal{F}%
( \Omega) $ is a Hilbert space also with respect to the inner
product $\mathcal{E}( \varphi,\psi) $. The harmonicity
of $u$
in (\ref{D}) means that
%
%e7.2 ###
%
\begin{equation} \label{ufi}
\mathcal{E}( u,\varphi) =0 \qquad\mbox{for all }\varphi\in
\mathcal{F%
}( \Omega) .
\end{equation}
Equivalently, this means for the function $v=f-u$ that%
%
%e7.3 ###
%
\begin{equation}\label{vfi}
\mathcal{E}( v,\varphi) =\mathcal{E}( f,\varphi
) %
\qquad\mbox{for all }\varphi\in\mathcal{F}( \Omega) .\vadjust{\goodbreak}
\end{equation}
Since $\mathcal{E}( f,\varphi) \leq\mathcal{E}(
f)
^{1/2}\mathcal{E}( \varphi) ^{1/2}$, the functional
$\varphi
\mapsto\mathcal{E}( f,\varphi) $ is a bounded linear functional
in $\mathcal{F}( \Omega) $, and equation (\ref{vfi})
has a
unique solution $v\in\mathcal{F}( \Omega) $ by the Riesz
representation theorem. Then $u=f-v$ is a unique solution of~(\ref{D}).

(b) Setting $\varphi=w-u$ and noticing that $\varphi
\in
\mathcal{F}( \Omega) $, we obtain using (\ref{ufi})%
\[
\mathcal{E}( w) =\mathcal{E}( u+\varphi)
=\mathcal{E}%
( u) +2\mathcal{E}( u,\varphi) +\mathcal
{E}(
\varphi) =\mathcal{E}( u) +\mathcal{E}(
\varphi
) .
\]
Hence, $\mathcal{E}( u) \leq\mathcal{E}( w
) $, and the
equality is attained when $\mathcal{E}( \varphi) =0$,
that is,
when $\varphi=0$.
\end{pf}

In what follows, denote by $R$ the resolvent operator of (\ref{D}),
that is,
$u=Rf$. Obviously, $R$ is a linear operator in $\mathcal{F}$. Since by
Lemma %
\ref{LemDir1} $\mathcal{E}( Rf) \leq\mathcal{E}(
f) $%
, we see that the norm of the operator $R$ in $\mathcal{F}$ is bounded
by $1$.
\begin{lemma}
\label{LemDir2}
\textup{(a)} If $f\leq g$, then $Rf\leq Rg$. In
particular, if $f\geq0$, then $Rf\geq0$.

\textup{(b)} If $0\leq f\leq1$, then also $0\leq Rf\leq1$.

\textup{(c)} If $\{ f_{n}\} _{n=1}^{\infty}$ is an
increasing sequence from $\mathcal{F}$ and $f_{n}\stackrel{\mathcal
{F}}{%
\rightarrow}f$ as $n\rightarrow\infty$, then $Rf_{n}\rightarrow Rf$
a.e. in $\Omega$ as $n\rightarrow\infty$.
\end{lemma}
\begin{pf}
(a) The function $u=Rf-Rg$ is harmonic in $B$ and satisfies
the boundary condition $u\leq0$ $\func{mod}\mathcal{F}( \Omega
)
$. By \cite{GrigHu}, Lemma 4.4, the latter condition implies $u_{+}\in
\mathcal{F}( \Omega) $. Substituting $\varphi=u_{+}$
into (\ref{ufi}), we obtain $\mathcal{E}( u,u_{+}) =0$. On the
other hand,
by \cite{GrigHu}, Lemma 4.3, $\mathcal{E}( u,u_{+}) \geq
\mathcal{E%
}( u_{+}) $, whence it follows that $\mathcal{E}(
u_{+}) =0$ and, hence, $u_{+}=0$. Consequently, $u\leq0$ and
$Rf\leq
Rg$.

(b) Set $u=Rf$ and $w=u\wedge1$ so that $u,w\in
\mathcal{F}$
and $\mathcal{E}( w) \leq\mathcal{E}( u)
$. Setting $%
\varphi=u-f$ and $\psi=w-f$, we see that $\varphi\in\mathcal
{F}(
\Omega) $, $\psi\in\mathcal{F}$ and $\psi\leq
\varphi$. By~\cite{GrigHu}, Lemma 4.4, we conclude that $\psi_{+}\in\mathcal
{F}%
( \Omega) $. On the other hand, we have $\psi
_{-}=\varphi
_{-}\in\mathcal{F}( \Omega) $ whence $\psi\in\mathcal
{F}%
( \Omega) $. It follows that $w=f \func{mod}\mathcal
{F}(
\Omega) $. By Lemma~\ref{LemDir1} we conclude that $\mathcal
{E}(
u) \leq\mathcal{E}( w) $. Since the opposite
inequality is
true by the definition of a Dirichlet form, we see that $\mathcal
{E}(
w) =\mathcal{E}( u) $. It follows from Lemma \ref{LemDir1}
that $w=u$, which implies $u\leq1$.

(c) Since\vspace*{-1pt} $R$ is a bounded operator in $\mathcal{F}$,
we see
that $Rf_{n}\stackrel{\mathcal{F}}{\rightarrow}Rf$ as $n\rightarrow
\infty$%
. It follows that also $Rf_{n}\stackrel{L^{2}( \Omega) }{
\rightarrow}Rf$. Then\vspace*{1pt} there is a subsequence of $\{ Rf_{n}
\} $
that converges to $Rf$ almost everywhere in $\Omega$. Finally, since the
sequence $\{ Rf_{n}\} $ is monotone increasing, the entire
sequence $\{ Rf_{n}\} $ also converges to~$Rf$ almost everywhere
in $\Omega$.
\end{pf}

%s7.2 ###
\subsection{Some consequences of the main hypotheses}
\label{SecRVD}

The next lemma states useful consequences of the main
hypotheses and motivates the statement of Theorem \ref{Tconv} below.
It is
also used in the proof of Corollary \ref{CorUE}.\vadjust{\goodbreak}
\begin{lemma}
\label{Lemconv}Let all metric balls be precompact. Then the following
implications are true:

\begin{longlist}[(a)]
\item[(a)] (\ref{condH}) implies that the
metric space $%
( M,d) $ is connected;

\item[(b)] (\ref{EFF}) implies that the Dirichlet
form $( \mathcal{E},\mathcal{F}) $ is conservative;

\item[(c)] (\ref{EFF}) implies that
$\func{diam}%
M=\infty$.
\end{longlist}
\end{lemma}
\begin{pf}
(a) Assume that $( M,d) $ is
disconnected, and let $%
\Omega$ be a nonempty open subset of $M$ such that $\Omega^{c}$ is also
nonempty and open. There is a big enough ball $B\subset M$ such that the
intersections of $\delta B$ both with $\Omega$ and $\Omega^{c}$ are
nonempty, where $\delta$ is the parameter from (\ref{condH}).
Since $%
\overline{B}\cap\Omega$ is a compact set, there is a cutoff function $u$
of $\overline{B}\cap\Omega$ in $\Omega$; that is, $u\in\mathcal
{F}\cap
C_{0}( \Omega) $ and $u\equiv1$ in a neighborhood of
$\overline{%
B}\cap\Omega$. Obviously, $u\equiv0$ in $\Omega^{c}$. We claim that $u$
is harmonic in $B$. Indeed, for every function $v\in\mathcal{F}\cap
C_{0}( B) $, we have $uv\in\mathcal{F}\cap C_{0}(
B\cap
\Omega) $ and%
\[
\mathcal{E}( u,v) =\mathcal{E}( u,uv)
+\mathcal{E}%
( u,v-uv) .
\]
Since $\limfunc{supp}( uv) \subset\overline{B}\cap
\Omega$ and,
hence, $u\equiv1$ in a neighborhood of $\limfunc{supp}(
uv) $,
we obtain by the strong locality of $( \mathcal{E},\mathcal
{F}) $
that $\mathcal{E}( u,uv) =0$. Since
\[
\limfunc{supp}\bigl( v( 1-u) \bigr) \subset\overline
{B}\cap
( \overline{B}\cap\Omega) ^{c}=\overline{B}\cap\Omega^{c}
\]
and $u=0$ in $\Omega^{c}$, it follows that $\mathcal{E}(
u,v-uv)
=0$. Hence $\mathcal{E}( u,v) =0$ and $u$ is a nonnegative
harmonic function in $B$. However, the function $u$ does not satisfy
(\ref{condH}) because $u$ takes in $\delta B$ the values $1$ and $0$.

(b) By Corollary \ref{CEF}, (\ref{EFF}) implies
that%
\[
\mathbb{P}_{x}\bigl( \tau_{B( x,R) }\leq t\bigr) \leq
C\exp
\biggl( -c\biggl( \frac{F( R) }{t}\biggr) ^{
{1}/({\beta^{\prime
}-1})}\biggr)
\]
for any $x\in M\setminus\mathcal{N}_{0}$, $R>0$, $t>0$. Using this estimate
and (\ref{PtOm}), we obtain%
\begin{eqnarray*}
\mathcal{P}_{t}1( x) &\geq&\mathcal{P}_{t}^{B(
x,R)
}1( x) \\
&=&\mathbb{P}_{x}\bigl( \tau_{B( x,R) }>t\bigr) \\
&\geq&1-C\exp\biggl( -c\biggl( \frac{F( R) }{t}\biggr)
^{{1}/({%
\beta^{\prime}-1})}\biggr) .
\end{eqnarray*}
As $R\rightarrow\infty$, we see that $\mathcal{P}_{t}1(
x) \geq
1$, which proves the stochastic completeness.

(c) If $\func{diam}M=R<\infty, $ then $M=B_{R}$ so
that the
exit time from $B_{R}$ is~$\infty$ and $({E}_{F}\mbox{$\leq$}) $ fails.
\end{pf}

%s7.3 ###
\subsection{The converse theorem}

In the next statement, we use weaker versions of (\ref{UEE})
and (\ref{NLE}) that will be denoted by $(\mathit{UE}_{\mathrm{weak}}
) $ and $%
(\mathit{NLE}_{\mathrm{weak}}) $. Namely, in each of these conditions we assume
that the heat kernel exists as a measurable integral kernel of the heat
semigroup $\{ P_{t}\} $ and satisfies the estimates (\ref{UEE}) and (\ref{NLE})\vadjust{\goodbreak}
for all $t>0$ and for \textit
{almost} all
$x,y\in M$. Note that unlike the conditions (\ref{UEE}) and
(\ref{NLE}), their weak versions do not use the diffusion process
$\{
X_{t}\} $.

\begin{theorem}
\label{Tconv}Assume that all metric balls are precompact and\break $\func
{diam}%
M=\infty$. Then the following sets of conditions are equivalent:

\begin{longlist}
\item $\mbox{(\ref{VD})} + \mbox{(\ref{condH})}
+ \mbox{(\ref{EFF})}$;

\item $\mbox{(\ref{VD})} + \mbox{(\ref{UEE})}
+\mbox{(\ref{NLE})}$, and the heat kernel is H\"{o}lder continuous
outside a
properly exceptional set;

\item $\mbox{(\ref{VD})} +(
\mathit{UE}_{\mathrm{weak}})
+(\mathit{NLE}_{\mathrm{weak}}) $.
\end{longlist}
\end{theorem}

Note that, by Lemma \ref{Lemconv}, (i) implies that
$\func{%
diam}M=\infty$. However, neither of conditions (ii)
or (iii) implies that $M$ is unbounded because (ii)
is satisfied on any compact Riemannian manifold.
\begin{pf*}{Proof of Theorem \ref{Tconv}}
The implication $\mbox{(i)} \Rightarrow\mbox{(ii)}$ is contained
in Theorem \ref{Tmain}, and the implication $\mbox{(ii)}
\Rightarrow
\mbox{(iii)}$ is trivial. In what follows we prove the
implication $\mbox{(iii)} \Rightarrow \mbox{(i)}$.

Assuming (iii), let us first show that $M$ is connected.
Indeed, let $M$ split into a disjoint union of two nonempty open sets $
\Omega_{1}$ and $\Omega_{2}$. By the continuity of the paths of
$\{
X_{t}\} $, we have $p_{t}( x,y) =0$ for all $t>0$
and $x\in
\Omega_{1}\setminus\mathcal{N}$, $y\in\Omega_{2}\setminus\mathcal{N}$,
whereas by (\ref{NLE}) we have $p_{t}( x,y) >0$ whenever
$t>\eta^{-1}d( x,y) $. This contradictions proves the
connectedness of $M$. By \cite{GrigHuUpper}, Corollary 5.3,~(\ref{VD}),
the connectedness, and the unboundedness of $M$ imply the
\textit{%
reverse volume doubling} (\ref{RVD}); that is, the
following
inequality holds:\label{condRVD}%
{\renewcommand{\theequation}{\textit{RVD}}
\begin{equation}\label{RVD}
\frac{V( x,R) }{V( x,r) }\geq c\biggl( \frac
{R}{r}%
\biggr) ^{\alpha^{\prime}},
\end{equation}}

\vspace*{-8pt}

\noindent
which holds for all $x\in M$, $0<r\leq R$, with some positive constants
$%
c,\alpha^{\prime}$. By~\cite{GrigHuUpper}, Theorem 2.2 and Section 6.4 (see
also \cite{KigamiNash}), $\mbox{(\ref{VD})} + \mbox{(\ref{RVD})}
+(\mathit{UE}_{\mathrm{weak}}) $ imply $({E}_{F}\mbox{$\leq$}) $.\footnote{%
Note that (\ref{RVD}) is essential for $({E}_{F}\mbox{$\leq$}
) $
(see \cite{GrigHuUpper}, Theorem 2.2). In fact, it was shown in \cite
{GrigHuUpper} and \cite{KigamiNash} that $\mbox{(\ref{VD})} + \mbox
{(\ref{RVD})}
+(\mathit{UE}_{\mathrm{weak}}) $ imply also $({E}_{F}\mbox{$\geq$}
) $
provided the Dirichlet form is conservative. In our setting the
conservativeness of the Dirichlet form can also be proved but a direct proof
of $({E}_{F}\mbox{$\geq$}) $ is shorter.}

Let us now prove $({E}_{F}\mbox{$\geq$}) $, that is,
%
%e7.4 ###
%
\setcounter{equation}{3}
\begin{equation} \label{PF}
\int_{0}^{\infty}\mathcal{P}_{t}^{B( x,R) }1(
x)
\,dt\geq cF( R)
\end{equation}
for all $x\in M\setminus\mathcal{N}$ and $R>0$, where $\mathcal{N}$
is a
properly exceptional set. It suffices to prove that there is a constant
$%
\zeta>0$ such that, for any ball $B=B( x_{0},R) $,%
%
%e7.5 ###
%
\begin{equation} \label{PBB}
\int_{0}^{\infty}P_{t}^{B}1\,dt\geq cF( R)
\qquad\mbox{a.e. in }\zeta B.
\end{equation}
Indeed, the function%
\[
u=\int_{0}^{\infty}\mathcal{P}_{t}^{B}1\,dt=G^{B}1
\]
is quasi-continuous by \cite{FOT}, Theorem 4.2.3.\label{remcheckthisinFOT}
By \cite{GrigHuUpper}, Proposition 6.1, if $u( x)
\geq a$
for almost all $x\in\Omega$, where $a$ is a constant and $\Omega$ is an
open set, then $u( x) \geq a$ for all $x\in\Omega
\setminus
\mathcal{N}$ where $\mathcal{N}$ is a properly exceptional set.
Hence,~(\ref{PBB}) implies that%
%
%e7.6 ###
%
\begin{equation} \label{PFN}
\int_{0}^{\infty}\mathcal{P}_{t}^{B}1( x) \,dt\geq
cF(
R) \qquad\mbox{for all }x\in\zeta B\setminus\mathcal{N}
\end{equation}
for some properly exceptional set $\mathcal{N}=\mathcal{N}_{B}$.
Taking the
union of such sets $\mathcal{N}_{B}$ where $B$ varies over a countable
family $S$ of all balls with rational radii and whose centers form a dense
subset of $M$, we obtain a properly exceptional set $\mathcal{N}$ such
that (%
\ref{PFN}) holds for any ball $B\in S$. Approximating any ball $B$ from
inside by balls of the family $S$, we obtain (\ref{PFN}) for all balls,
which implies (\ref{PF}).

Now let us prove (\ref{PBB}). By the comparison principle of
\cite{GrigHuUpper}, Proposition~4.7 (see also \cite{GrigHu}, Lemma
4.18), we have,
for any nonnegative function $f\in L^{2}\cap L^{\infty}(
B) $,
%
%e7.7 ###
%
\begin{equation}\label{eqs-2}
P_{t}f( x) \leq P_{t}^{B}f( x) +\sup_{s\in
(0,t]}%
\limfunc{esup}_{y\in B\setminus({1/2})B}P_{s}f( y)
\end{equation}
for almost all $x\in B$. Let $\zeta$ be a small positive constant to be
specified below, and set $f=\mathbf{1}_{\zeta B}$. It follows from
$(\mathit{NLE}_{\mathrm{weak}}) $ and (\ref{Va}) that%
%
%e7.8 ###
%
\begin{equation} \label{infp}
p_{t}( x,z) \geq\frac{c}{V( x_{0},\mathcal
{R}(
t) ) }\qquad\mbox{for a.a. }x,z\in B\biggl(x_{0},\frac
{1}{2}\eta
\mathcal{R}( t) \biggr),
\end{equation}
provided $0<t\leq\varepsilon F( R) $. The initial value
of $%
\varepsilon$ is given by the condition $(\mathit{NLE}_{\mathrm{weak}}) $
but we
are going to further reduce this value of $\varepsilon$ in the course of
the proof. Assume that $t$ varies in the following interval:%
%
%e7.9 ###
%
\begin{equation} \label{tR}
\tfrac{1}{2}\varepsilon F( R) \leq t\leq\varepsilon
F(R) .
\end{equation}
The left-hand side inequality in (\ref{tR}) implies by (\ref{Fb}) that
%
%e7.10 ###
%
\begin{equation} \label{Rt}
R\leq C\biggl( \frac{1}{\varepsilon}\biggr) ^{1/\beta}\mathcal
{R}(
t) .
\end{equation}
Chose $\zeta$ from the identity%
%
%e7.11 ###
%
\begin{equation}\label{zeta}
\zeta C\biggl( \frac{1}{\varepsilon}\biggr) ^{1/\beta}=\frac
{1}{2}\eta
\end{equation}
so that (\ref{Rt}) implies%
\[
B( x_{0},\zeta R) \subset B\bigl(x_{0},\tfrac{1}{2}\eta
\mathcal{R}%
( t) \bigr).
\]
Integrating (\ref{infp}) over $B( x_{0},\zeta R) $ and
using (\ref{VD}) and (\ref{zeta}), we obtain%
%
%e7.12 ###
%
\begin{eqnarray} \label{eqs-4}
P_{t}f( x) &=&\int_{B( x_{0},\zeta R)
}p_{t}(
x,z) \,d\mu( z) \nonumber\\
&\geq&\frac{cV( x_{0},\zeta R) }{V( x_{0},\mathcal
{R}(
t) ) } \nonumber\\[-8pt]\\[-8pt]
&\geq&c\zeta^{\alpha} \nonumber\\
&=&c^{\prime}\varepsilon^{\alpha/\beta}\nonumber
\end{eqnarray}
for almost all $x\in B( x_{0},\zeta R) $. On the other
hand, for
almost all $y\in B\setminus\frac{1}{2}B$, we have by $(
\mathit{UE}_{\mathrm{weak}}) $ and Lemma \ref{LemtiFi}
\begin{eqnarray*}
P_{s}f( y) &=&\int_{B( x_{0},\zeta R)
}p_{s}(
y,z) \,d\mu( z) \\
&\leq&C\frac{V( x_{0},R) }{V( y,\mathcal{R}
( s)
) }\exp\biggl( -c\biggl( \frac{F( R) }{s}\biggr)
^{{1}/({%
\beta^{\prime}-1})}\biggr) ,
\end{eqnarray*}
where we have used that $d( y,z) \simeq R$ and $s\leq
t<F(
R) $. Using (\ref{Va}) and~(\ref{Fb}) we obtain%
\[
\frac{V( x_{0},R) }{V( y,\mathcal{R}(
s) ) }%
\leq C\biggl( \frac{R}{\mathcal{R}( s) }\biggr)
^{\alpha}\leq
C^{\prime}\biggl( \frac{F( R) }{s}\biggr) ^{\alpha
/\beta}.
\]
Finally, it follows from (\ref{tR}) and $s\leq t$ that $\frac{F(
R) }{s}\geq\frac{1}{\varepsilon}$ whence%
%
%e7.13 ###
%
\begin{equation} \label{eqs-3}
P_{s}f( y) \leq C\biggl( \frac{1}{\varepsilon}\biggr)
^{\alpha
/\beta}\exp\biggl( -c\biggl( \frac{1}{\varepsilon}\biggr) ^{
{1}/({\beta
^{\prime}-1})}\biggr)
\end{equation}
for almost all $y\in B\setminus\frac{1}{2}B$. Combining (\ref{eqs-2}),
(\ref{eqs-4}) and (\ref{eqs-3}), we obtain, for almost all $x\in B(
x_{0},\zeta R) $,%
\begin{eqnarray*}
P_{t}^{B}f( x) &\geq&P_{t}f( x) -\sup_{s\in
(0,t]}%
\limfunc{esup}_{B\setminus K}P_{s}f \\
&\geq&c^{\prime}\varepsilon^{\alpha/\beta}-C\biggl( \frac
{1}{\varepsilon
}\biggr) ^{\alpha/\beta}\exp\biggl( -c\biggl( \frac{1}{\varepsilon
}\biggr)
^{{1}/({\beta^{\prime}-1})}\biggr) \\
&\geq&\frac{1}{2}c{^{\prime}}\varepsilon^{\alpha/\beta},
\end{eqnarray*}
provided $\varepsilon$ is chosen small enough. The path $t\mapsto
P_{t}^{B}f $ is a continuous path in $L^{2}( B) $ and,
hence, can
be integrated in $t$. It follows from the previous inequality that%
\[
\int_{0}^{\infty}P_{t}^{B}1\,dt\geq\int_{({1/2})\varepsilon
F(
R) }^{\varepsilon F( R) }P_{t}^{B}f \,dt\geq
c\varepsilon
^{\alpha/\beta+1}F( R) ,
\]
which finishes the proof of $({E}_{F}\mbox{$\geq$}) $.\vadjust{\goodbreak}

We are left to prove that $\mbox{(iii)} \Rightarrow\mbox{(\ref{condH})}$.
By \cite{BarGrigKumHar}, Theorem 3.1
(see also \cite{FS} and \cite{HebSChar}, Theorem 5.3), $\mbox{(\ref
{VD})} +(
\mathit{UE}_{\mathrm{weak}}) +(\mathit{NLE}_{\mathrm{weak}}) $ imply the \textit
{parabolic} Harnack inequality for
bounded caloric function and, hence, the Harnack inequality (\ref{condH})
for bounded harmonic functions (note that this result uses the
precompactness of the balls). We still have to obtain (\ref{condH}) for
all nonnegative harmonic functions. Note that by \cite{GrigHuUpper},
Theorem 2.1,
\[
\mbox{(\ref{VD})} + \mbox{(\ref{RVD})} +(\mathit{UE}_{\mathrm{weak}})
\Rightarrow
\mbox{(\ref{FK})}.
\]
In particular, for any ball $B$, we have $\lambda_{\min}(
B) >0$%
. Given a function $u\in\mathcal{F}$ that is nonnegative and
harmonic in a
ball $B\subset M$, set $f_{n}=u\wedge n$ for any $n\in\mathbb{N}$, and
denote by $u_{n}$ the solution of the Dirichlet problem%
\[
\cases{
u_{n}\mbox{ is harmonic in }B, \cr
u_{n}=f_{n} \func{mod}\mathcal{F}( B);}
\]
cf. Section \ref{SecHarm}. Since $0\leq f_{n}\leq n$, we have also
$0\leq
u_{n}\leq n$. Since the sequence~$\{ f_{n}\} $ increases
and $%
f_{n}\stackrel{\mathcal{F}}{\rightarrow}u$ (cf. \cite{FOT}, Theorem
1.4.2), it
follows by Lemma~\ref{LemDir2} that $u_{n}\rightarrow u$ almost everywhere
in $B$. Each function $u_{n}$ is bounded and, hence, satisfies the Harnack
inequality in $B$, that is,%
\[
\limfunc{esup}_{\delta B}u_{n}\leq C\limfunc{einf}_{\delta B}u_{n}.
\]
Replacing in the right-hand side $u_{n}$ by a larger function $u$ and
passing to the limit in the left-hand side as $n\rightarrow\infty$, we
obtain the same inequality for~$u$, which was to be proved.
\end{pf*}
\begin{corollary}
\label{CorUE}Assume that all metric balls are precompact, $\func{diam}
M=\infty$, and the Dirichlet form $( \mathcal{E},\mathcal
{F}) $
is conservative. Then the following sets of conditions are equivalent:

\begin{longlist}
\item $\mbox{(\ref{VD})} +\mbox{(\ref{condH})}
+(\mathit{UE}
_{\mathrm{weak}}) $;

\item $\mbox{(\ref{VD})} +\mbox{(\ref{UEE})}
+ \mbox{(\ref{NLE})}$.
\end{longlist}
\end{corollary}
\begin{pf}
In the view of Theorem \ref{Tconv}, it suffices to prove that $\mbox{(i)}
\Rightarrow\mbox{(\ref{EFF})} $. By Lemma \ref{Lemconv},
(\ref{condH}) implies the connectedness of $M$. By \cite{GrigHuUpper},
Corollary~5.3,
$\mbox{(\ref{VD})} \Rightarrow\mbox{(\ref{RVD})}$ provided $M$ is connected
and unbounded, which is the case now. By \cite{GrigHuUpper}, Theorem
2.2, the
conservativeness and $\mbox{(\ref{VD})} + \mbox{(\ref{RVD})} +(
\mathit{UE}_{\mathrm{weak}}) $ imply (\ref{EFF}).
\end{pf}

Many equivalent conditions for $(\mathit{UE}_{\mathrm{weak}}) $ were proved in
\cite{GrigHuUpper} under the standing assumptions $\mbox{(\ref{VD})}
+ \mbox{(\ref{RVD})}$ and the conservativeness of $( \mathcal
{E},\mathcal{F}%
) $. Of course, each of these conditions can replace $(
\mathit{UE}_{\mathrm{weak}}) $ in the statement of Corollary \ref{CorUE}.
\begin{corollary}
\label{Corunbound}Assume that all metric balls are precompact, $\func
{diam}%
M=\infty$, and $( M,d) $ satisfies the\vadjust{\goodbreak} chain condition.
Then the
following two sets of conditions are equivalent:

\begin{longlist}
\item $\mbox{(\ref{VD})} +\mbox{(\ref{condH})}
+ \mbox{(\ref{EFF})}$;

\item The heat kernel exists and satisfies the
two-sided estimate (\ref{twosided}).
\end{longlist}
\end{corollary}
\begin{pf}
The implication $\mbox{(i)} \Rightarrow\mbox{(ii)} $ is contained
in Corollary \ref{Cortwo}. Let us prove the implication
$\mbox{(ii)}\Rightarrow\mbox{(i)} $. Estimate (\ref{twosided}) implies (\ref{UEE}) as well
as~(\ref{NLE}) with any value
of $\eta$%
, in particular, $\eta>1$; cf. Remark \ref{RemNLE}. By
\cite{GrigHuLauD}, Lemma 4.1, (\ref{NLE}) with $\eta>1$ implies (\ref{VD}).
Finally, by Theorem \ref{Tconv}, we obtain $\mbox{(\ref{condH})}
+ \mbox{(\ref{EFF})}$.
\end{pf}

\begin{appendix}
\section*{Appendix: List of conditions}

We briefly list the lettered conditions used in this paper with references
to the appropriate places in the main body.

\begin{longlist}[$(\mathrm{RVD}) $]
\item[$(H)$] $\mathop{\limfunc{esup}}_{B(
x,\delta r) }u\leq C \mathop{\func{einf}}_{B( x,\delta
r) }u$
(Section \ref{SecHarnack});

\item[$(\mathit{VD})$] $V( x,2r) \leq CV(
x,r) $
(Section \ref{SecHarnack});

\item[$(E_{F})$] $\mathbb{E}_{x}\tau_{B(
x,r)
}\simeq F( r) $ (Section \ref{SecEF});

\item[$(\mathit{FK})$] $\lambda_{\min}( \Omega)
\geq\frac{%
c}{F( R) }( \frac{\mu( B) }{\mu
( \Omega
) }) ^{\nu}$ (Section \ref{SecEF});

\item[$(\mathit{UE})$] $p_{t}( x,y) \leq\frac
{C}{V( x,%
\mathcal{R}( t) ) }\exp( -\frac{1}{2}\Phi
(
cd( x,y) ,t) ) $ (Section \ref{SecDUE});

\item[$(\mathit{NLE})$] $p_{t}( x,y) \geq\frac{c}{V(
x,\mathcal{R}%
( t) ) } $ provided $d( x,y) \leq
\eta
\mathcal{R}( t) $ (Section \ref{Seclow});

\item[$(\mathit{RVD}) $] $\frac{V( x,R) }{V(
x,r) }%
\geq c( \frac{R}{r}) ^{\alpha^{\prime}}$ (Section
\ref{SecRVD}).
\end{longlist}
\end{appendix}

\makeatletter\write@toc@ignorecontentsline\makeatother
\section*{Acknowledgments}
The authors thank Martin Barlow, Jiaxin Hu, Jun Kigami, Takashi Kumagai
and Alexander Teplyaev for valuable conversations on the topics of this
paper. The authors are grateful to the unnamed referees for careful reading
of the manuscript and for the useful remarks.
\makeatletter\write@toc@restorecontentsline\makeatother

%suskaldyti doi

% imsref loaded by lrinkeviciute, 2011-04-12 10:08:21
% imsref loaded by lrinkeviciute, 2011-04-12 10:28:35
%
% imsref loaded by lrinkeviciute, 2011-04-12 10:34:31
%
% imsref loaded by lrinkeviciute, 2011-04-12 10:40:12
%

%
\printaddresses

\end{document}